%% file: main.tex
\documentclass[nonblindrev]{informs3}
\OneAndAHalfSpacedXI % current default line spacing
\usepackage{microtype}
\usepackage{graphicx}
\usepackage{algorithm}
\usepackage{algorithmic}
\usepackage{subfigure}
\usepackage{booktabs} % for professional tables
\usepackage{multirow}
\usepackage{changepage}
\usepackage{svg}

\usepackage[colorinlistoftodos,textsize=tiny%,disable
]{todonotes}

\usepackage{hyperref}
\usepackage{natbib}
 \bibpunct[, ]{(}{)}{,}{a}{}{,}%
 \def\bibfont{\small}%

\usepackage{comment}

\OneAndAHalfSpacedXI

\TheoremsNumberedThrough     % Preferred (Theorem 1, Lemma 1, Theorem 2)
\ECRepeatTheorems

\usepackage{amsmath}
\usepackage{amssymb}
\usepackage{mathtools}

\usepackage[capitalize,noabbrev]{cleveref}
\crefname{assumption}{Assumption}{assumptions}
\theoremstyle{plain}
\newcommand{\out}{\text{out}}
\newcommand{\inn}{\text{in}}

\newcommand{\revision}[1]{#1}

\newcommand{\minor}[1]{#1}

\begin{document}

 \RUNAUTHOR{Qu, Galichon, Gao, Ugander}

 \RUNTITLE{Sinkhorn's Algorithm and Choice Modeling}

% Enter the full title:
\TITLE{On Sinkhorn's Algorithm and Choice Modeling}

% List of affiliations: The first argument should be a (short)
% identifier you will use later to specify author affiliations
% Academic affiliations should list Department, University, City, Region, Country
% Industry affiliations should list Company, City, Region, Country

% You can specify symbols, otherwise they are numbered in order.
% Ideally, you should not use this facility. Affiliations will be numbered
% in order of appearance and this is the preferred way.
\ARTICLEAUTHORS{%
\AUTHOR{Zhaonan Qu}
\AFF{Data Science Institute, Columbia University, New York, USA \EMAIL{zq2236@columbia.edu}} 
\AUTHOR{Alfred Galichon}
\AFF{Department of Mathematics and Department of Economics, New York University, New York, USA\\ Department of Economics, Sciences Po, Paris, France \EMAIL{alfred.galichon@nyu.edu}}
\AUTHOR{Wenzhi Gao}
\AFF{Institute for Computational and Mathematical Engineering, Stanford University, California, USA \EMAIL{gwz@stanford.edu}}
\AUTHOR{Johan Ugander}
\AFF{Department of Management Science \& Engineering, Stanford University, California, USA \EMAIL{jugander@stanford.edu}}
% Enter all authors
} % end of the block

\ABSTRACT{For a broad class of models widely used in practice for choice and ranking data based on Luce's choice axiom, including the Bradley--Terry--Luce and Plackett--Luce models, we show that the associated maximum likelihood estimation problems are equivalent to a classic matrix balancing problem with target row and column sums. This perspective opens doors between two seemingly unrelated research areas, and allows us to unify existing algorithms in the choice modeling literature as special instances or analogs of Sinkhorn's celebrated algorithm for matrix balancing. We draw inspirations from these connections and resolve some open problems on the study of Sinkhorn's algorithm. We establish the global linear convergence of Sinkhorn's algorithm for non-negative matrices whenever finite scaling matrices exist, and characterize its linear convergence rate in terms of the algebraic connectivity of a weighted bipartite graph. We further derive the sharp asymptotic rate of linear convergence, which generalizes a classic result of Knight (2008). To our knowledge, these are the first quantitative linear convergence results for Sinkhorn's algorithm for general non-negative matrices and positive marginals. Our results highlight the importance of connectivity and orthogonality structures in matrix balancing and Sinkhorn's algorithm, which could be of independent interest. More broadly, the connections we establish in this paper between matrix balancing and choice modeling could also help motivate further transmission of ideas and lead to interesting results in both disciplines.
}
\maketitle

\input{intro}

\input{literature}

\input{preliminaries}
\input{equivalence}
\input{convergence-analyses}
\input{conclusion}

\bibliography{sinkhorn-choice}
\bibliographystyle{informs2014}

%%%%%%%%%%%%%%%%%%%%%%%%%%%%%%%%%%%%%%%%%%%%%%%%%%%%%%%%%%%%%%%%%%%%%%%%%%%%%%%
%%%%%%%%%%%%%%%%%%%%%%%%%%%%%%%%%%%%%%%%%%%%%%%%%%%%%%%%%%%%%%%%%%%%%%%%%%%%%%%
% APPENDIX
%%%%%%%%%%%%%%%%%%%%%%%%%%%%%%%%%%%%%%%%%%%%%%%%%%%%%%%%%%%%%%%%%%%%%%%%%%%%%%%
%%%%%%%%%%%%%%%%%%%%%%%%%%%%%%%%%%%%%%%%%%%%%%%%%%%%%%%%%%%%%%%%%%%%%%%%%%%%%%%
\newpage
\renewcommand{\theHchapter}{A\arabic{chapter}}

\begin{APPENDICES}
\crefalias{section}{appendix}
\input{connection}
\input{regularization}
\input{empirical}
\input{proofs}
 \end{APPENDICES}
%%%%%%%%%%%%%%%%%%%%%%%%%%%%%%%%%%%%%%%%%%%%%%%%%%%%%%%%%%%%%%%%%%%%%%%%%%%%%%%
%%%%%%%%%%%%%%%%%%%%%%%%%%%%%%%%%%%%%%%%%%%%%%%%%%%%%%%%%%%%%%%%%%%%%%%%%%%%%%%

\end{document}

%% file: intro.tex
\section{Introduction}
The modeling of choice and ranking data is an important topic across many disciplines. Given a collection of $m$ objects, a universal problem is to aggregate choice or partial ranking data over them to arrive at a reasonable description of either the behavior of decision makers,  the intrinsic qualities of the objects, or both. Work on such problems dates back over a century at least to the work of Landau, who considered $m$ chess players and a record of their match results against one another, aiming to aggregate the pairwise comparisons to arrive at a
global ranking of all players \citep{landau1895relativen,elo1978rating}. More generally, comparison data can result from choices from subsets of \emph{varying} sizes, from partial or complete rankings of objects, or from mixtures of different data types.

The modern rigorous study of comparisons primarily builds on the foundational works of \citet{thurstone1927method} and \citet{zermelo1929berechnung}. Both proposed models based on a numerical ``score'' for each item (e.g., chess player), but with different specifications of choice probabilities.
\citet{zermelo1929berechnung} builds on the intuition that choice probability should be proportional to the score, and proposes a iterative algorithm to estimate the scores from pairwise comparison data. As one of the foundational works in this direction, \citet{luce2012individual}
formalized the multinomial logit model
of discrete choice starting from the axiom of \emph{independence of irrelevant alternatives}
(IIA). It states that the relative likelihood of choosing an item
$j$ over another item $k$ is independent of 
the presence of other alternatives. In other words, if $S$ and $S'$ are two subsets of the $m$
alternatives, both containing $j$ and $k$, and $\Pr(j, S)$ denotes
the probability of choosing item $j$ from $S$, then for $\Pr(k, S),\Pr(k, S')>0$,
\begin{align*}
\frac{\Pr(j, S)}{\Pr(k, S)} & =\frac{\Pr(j, S')}{\Pr(k, S')}.
\end{align*}
This invariance property, together with a natural condition for zero probability alternatives, are often referred to as Luce's choice axioms. They guarantee that each alternative can be summarized
by a non-negative score $s_{j}$ such that the probability of choice can be parameterized by
\begin{align}
\label{eq:Luce}
\Pr(j, S) & =\frac{s_{j}}{\sum_{k\in S}s_{k}},
\end{align}
 for any set $S$ that contains $j$. The parameters $s_{j}$ reflect
the ``intrinsic'' value of item $j$, and are unique up to a normalization,
which can be set to $\sum_{j}s_{j}=1$. This general choice model includes as a special case the Bradley--Terry--Luce model (BTL) for pairwise comparisons \citep{bradley1952rank}, and also applies to ranking data when each $k$-way ranking is broken down into $k-1$ choice observations, where an item is chosen over the set of items ranked lower \citep{critchlow1991probability,plackett1975analysis,hausman1987specifying}. The many subsequent works that build on Luce's choice axioms speak to its fundamental importance in choice modeling. Other works have also sought to address the limitations of Luce choice models. Prominent among them are probit models \citep{thurstone1927method,berkson1944application}, random utility (RUM) models \citep{mcfadden2000mixed}, context-dependent (CDM) models \citep{batsell1985new,seshadri2020learning}, and behavioral models from psychology \citep{tversky1972elimination}.

Matrix balancing (or scaling), meanwhile, is a seemingly unrelated mathematical problem with an equally long history. In its most common form that we study in this paper, the problem seeks positive diagonal matrices $D^0,D^1$ of a given (entry-wise) non-negative matrix $A\geq0$, such that the scaled matrix $D^0AD^1$ has row and column sums equal to some prescribed positive marginals $p,q$: 
\begin{align}
\label{eq:matrix-balancing}
\begin{split}
    (D^0AD^1)\mathbf{1} &= p,\\
    (D^0AD^1)^T\mathbf{1} &= q.
\end{split}
\end{align}
Over the years, numerous applications and problems across different domains, including statistics \citep{yule1912methods,deming1940least,ireland1968contingency}, 
economics \citep{stone1962multiple,bacharach1970biproportional,galichon2021matching}, transportation networks \citep{kruithof1937telefoonverkeersrekening,lamond1981bregman,chang2024inferring}, optimization \citep{bregman1967relaxation,ruiz2001scaling}, and machine learning \citep{cuturi2013sinkhorn,peyre2019computational}, have found themselves essentially solving a new incarnation of the old matrix balancing problem, which attests to its universality and importance. 

A major appeal of the matrix balancing problem lies in the simplicity and elegance of its popular solution method, widely known as Sinkhorn's algorithm \citep{sinkhorn1964relationship}. \revision{Observe that it is easy to scale the rows or columns of $A$ such that the resulting matrix satisfies one of the two marginal constraints in \eqref{eq:matrix-balancing}. However, it is more difficult to construct scalings $D^0,D^1$ that \emph{simultaneously} satisfy both constraints. Sinkhorn's algorithm (\cref{alg:scaling}) simply alternates between updating the scalings $D^0$ and $D^1$ to satisfy one of the two marginal conditions in the hope of converging to a solution, leading to lightweight implementations that have proved effective for practical problems of massive size.} In particular, Sinkhorn's algorithm has gained much popularity in the recent decade thanks to its empirical success at approximating optimal transport (OT) distances \citep{cuturi2013sinkhorn,altschuler2017near}, which are bedrocks of important recent topics in operations research, such as
Wasserstein distributionally robust optimization \citep{mohajerin2018data,gao2023distributionally,kuhn2019wasserstein,blanchet2019robust,blanchet2022distributionally}.

 \revision{Despite the widespread popularity of Sinkhorn's algorithm, its convergence behavior is yet to be 
fully understood. In particular, while there have been extensive studies of convergence, many focus on the setting when the matrix $A>0$, i.e., entry-wise strictly positive, which includes most OT problems. In contrast, other applications of matrix balancing, particularly those with network structures, have $A\geq0$ with zero elements, and is therefore potentially sparse. In this setting, quantitative analyses are less common and more fragmented, employing different assumptions whose connections and distinctions remain less clear. On one hand, works such as \citet{kalantari2008complexity}, \citet{chakrabarty2021better}, and \citet{leger2021gradient} have established \emph{global} (that is, true for all iterations $t\geq1$) sub-linear convergence results, i.e., convergence to an $\varepsilon$ accuracy solution requires a total number of iterations that is polynomial in $1/\varepsilon$. On the other hand, \citet{knight2008sinkhorn} establishes local (and more specifically \emph{asymptotic}) linear convergence for square matrix $A\geq 0$ and uniform marginals $p,q$. In other words, as $t\rightarrow \infty$, solution accuracy at iteration $t+1$ improves over that at iteration $t$ with a constant factor. Furthermore, a general result in \citet{luo1992convergence} implies global linear convergence of Sinkhorn's algorithm, i.e., convergence to an $\varepsilon$ accuracy solution requires iterations polynomial in $\log(1/\varepsilon)$. However, their result is implicit and does not characterize the dependence on problem parameter and structure as those in the sub-linear results.
}

\revision{
These results leave open several questions on the convergence of Sinkhorn's algorithm. First, when does a \emph{quantitative} global linear convergence result exist for $A\geq 0$, and how to characterize the global convergence rate in terms of the problem primitives? Second, how to characterize the sharp, i.e., best possible, asymptotic linear convergence rate $\lambda$ that is applicable to general non-negative $A$ and non-uniform $p,q$? Third, how to reconcile and clarify the linear vs. sub-linear convergence results under different assumptions on the problem structure? Given that many applications of matrix balancing with network structures correspond to the setting with sparse $A\geq0$, such as transportation and trade, it is therefore important to better understand the convergence of Sinkhorn's algorithm in this setting.}

In this paper, we provide answers to these open questions in matrix balancing. Surprisingly, the inspirations for our solutions come from results in the seemingly unrelated topic of choice modeling. Our main contributions are summarized below.

 \revision{Our first set of contributions, which we detail in \cref{sec:equivalence}, is recognizing Luce choice models as yet another instance where a central problem reduces to that of matrix balancing. More precisely, we formally establish the equivalence between the maximum likelihood estimation of Luce choice models and matrix balancing problems with an $A\geq 0$ with binary elements (\cref{prop:mle-scaling}).} We also clarify the relations and distinctions between problem assumptions in the two literatures (\cref{thm:necessary-and-sufficient}). More importantly, we demonstrate that classic and new algorithms from the choice literature, including those of \citet{zermelo1929berechnung,dykstra1956note,ford1957solution,hunter2004mm,maystre2017choicerank,agarwal2018accelerated}, can be viewed as special cases or analogs of Sinkhorn's algorithm when applied to various problems in the choice setting (\cref{thm:algorithm-equivalence}). These intimate mathematical and algorithmic connections allow us to provide a unifying perspective on works from both areas. More broadly, they enable researchers to import insights and tools from one domain to the other. In particular, recent works on choice modeling \citep{shah2015estimation,seshadri2020learning,vojnovic2020convergence} have highlighted the importance of \emph{algebraic connectivity} \citep{fiedler1973algebraic,spielman2007spectral} of the data structure for efficient parameter learning, which motivates us to also consider this quantity in the convergence analysis of Sinkhorn's algorithm.

  \revision{Our next set of contributions, detailed in \cref{sec:linear-convergence}, is establishing novel convergence bounds on Sinkhorn's algorithm, drawing from the connections to choice modeling that we establish. First, we provide a global linear convergence bound for Sinkhorn's algorithm, whenever the matrix balancing problem has a finite solution pair $D^0,D^1$ (\cref{thm:global-convergence}). We characterize the global convergence rate in terms of the algebraic connectivity of the weighted bipartite graph whose \emph{biadjacency} matrix is precisely $A$. To our knowledge, this result is the first to highlight the fundamental role of algebraic connectivity in the study of matrix balancing with sparse matrices.} %This result complements the linear convergence result of \citet{franklin1989scaling} for positive matrices, and the sub-linear bound of \citet{leger2021gradient} for non-negative matrices under weaker conditions. 
In addition, we characterize the asymptotic linear rate of convergence in terms of the scaled matrix $D^0AD^1$ with target marginals $p,q$, generalizing a result of \citet{knight2008sinkhorn} for uniform marginals and square matrices (\cref{thm:convergence}). This result employs a more explicit analysis that exploits an intrinsic orthogonality structure of Sinkhorn's algorithm. We also clarify the convergence behavior of Sinkhorn's algorithm under two regimes: when a finite scaling pair $D^0,D^1$ exists, Sinkhorn's algorithm converges linearly; otherwise, it only converges sub-linearly under the minimal conditions required for convergence (\cref{thm:lower-bound}). %Taken together, our results in \cref{sec:linear-convergence} provide solutions to the open questions concerning Sinkhorn's algorithm. They highlight the importance of spectral properties when $A$ is sparse.

%The challenges of ill-defined matrix balancing problems and non-convergence of Sinkhorn's algorithm in practice also motivate us to propose a regularized version of Sinkhorn's algorithm. It is inspired by the regularization of Luce choice models using Gamma priors, and 
 {Besides the contributions above, we further discuss connections between Sinkhorn's algorithm and topics in optimization and choice modeling in \cref{sec:connections}. For example, interpreting Sinkhorn's algorithm as a distributed optimization algorithm on a bipartite graph helps explain the importance of the spectral properties of the graph on its convergence. Inspired by Bayesian regularization of Luce choice models using gamma priors, we also design a regularized Sinkhorn's algorithm in \cref{subsec:regularization}, that is guaranteed to converge even when the standard algorithm does not, which is not uncommon when the data is very sparse and there are measurement errors.}

\revision{We believe that the connections we establish in this paper between choice modeling and matrix balancing can lead to further interesting results in both disciplines, and are therefore relevant to researchers working on related topics. In particular, the fundamental role of algebraic connectivity in the study of matrix balancing for sparse matrices goes beyond quantifying the algorithmic efficiency of Sinkhorn's algorithm. See, for example, \citet{chang2024inferring}, which quantifies the statistical efficiency of a network traffic model using algebraic connectivity.
}

%% file: literature.tex
\section{Related Work}
\label{sec:related-works}
This section includes an extensive review of related works in choice modeling and matrix balancing. Well-versed readers may skip ahead to the mathematical preliminaries (\cref{sec:formulations}) and our core results (\cref{sec:equivalence,sec:linear-convergence}).
\subsection{Choice Modeling}
Methods for aggregating choice and comparison data usually take one of two closely related approaches: maximum likelihood estimation of a statistical model or ranking according to the stationary distributions of a random walk on a Markov chain. Recent connections between maximum likelihood and spectral methods have put these two classes of approaches in increasingly close conversation with each other.

\textbf{Spectral Methods.}
The most well-known spectral method for rank aggregation is perhaps
the PageRank algorithm \citep{page1999pagerank}, which ranks web pages
based on the stationary distribution of a random walk on a
hyperlink graph. The use of stationary distributions also features in the work of \citet{dwork2001rank}, the Rank Centrality (RC) algorithm \citep{negahban2012iterative,negahban2016rank}, which generates consistent estimates for the Bradley--Terry--Luce
pairwise comparison model under assumptions on the sampling frame, and the Luce Spectral Ranking (LSR) and iterative LSR (I-LSR) algorithms of \citet{maystre2015fast} for choices from pairs as well as larger sets. Following that work, \citet{agarwal2018accelerated} proposed the Accelerated Spectral Ranking (ASR) algorithm with provably faster mixing times than RC and LSR, and better sample complexity bounds than \citet{negahban2016rank}.
\citet{knight2008sinkhorn} is an intriguing work partially motivated by \citet{page1999pagerank} that applies Sinkhorn's algorithm, which is central to the current work, to compute authority and hub scores similar to
those proposed by \citet{kleinberg1999authoritative} and \citet{tomlin2003new}, although the focus in \citet{knight2008sinkhorn} is on Markov chains rather than maximum likelihood estimation of choice models.
For ranking data, \citet{soufiani2013generalized} decompose rankings into pairwise comparisons and develop consistent estimators for Plackett--Luce models based on a generalized method of moments. 
Other notable works that make connections between Markov chains and 
choice modeling include \citet{blanchet2016markov} and \citet{ragain2016pairwise}.

\textbf{Maximum Likelihood Methods.}
Maximum likelihood estimation of the Bradley--Terry--Luce model 
dates back to \citet{zermelo1929berechnung}, \citet{dykstra1956note}, and \citet{ford1957solution}, which all give variants of the same iterative algorithm and prove its convergence to the MLE when the directed comparison graph is strongly connected. 
Much later, \citet{hunter2004mm} observed that their
algorithms are instances of a class of minorization-maximization or majorization-minimization (MM) algorithms, and develops MM algorithms for the Plackett--Luce model for ranking data, among others. \citet{vojnovic2020convergence,vojnovic2023accelerated} further investigate the convergence rate of MM algorithm for choice models, quantifying it in terms of the algebraic connectivity of the comparison graph. \citet{newman2023efficient} proposes an alternative to the classical iterative algorithm for pairwise comparisons based on a reformulated moment condition, achieving impressive empirical speedups. \citet{negahban2012iterative} is arguably the first work that connects maximum likelihood estimation to Markov chains, followed by \citet{maystre2015fast}, whose spectral method is based on a balance equation interpretation of the optimality condition.
\citet{kumar2015inverting} consider the problem of inverting the stationary distribution of a Markov
chain, and embed the maximum likelihood problem of the Luce choice model into this
framework, where the MLEs parameterize the desired transition matrix. \citet{maystre2017choicerank} consider the estimation of a network choice model with similarly parameterized random walks. Lastly, a vast literature in econometrics on discrete choice also considers different aspects of the ML estimation problem. In particular, the present paper is related to the Berry--Levinsohn--Pakes (BLP) framework of \citet{berry1995automobile}, well-known in econometrics. The matrix balancing interpretation of maximum likelihood estimation of choice models that we develop in this paper connects many of the aforementioned works.

\revision{Besides optimization problems related to maximum likelihood estimation, there have also been extensive studies on the
statistical properties of maximum likelihood estimates themselves \citep{hajek2014minimax,rajkumar2014statistical}. In particular, a line of recent works have highlighted the importance of {algebraic connectivity}---as quantified by the Fiedler eigenvalue \citep{fiedler1973algebraic,spielman2007spectral}---on the statistical efficiency of the MLEs. \citet{shah2015estimation} is the first to recognize this significance of data structure for the statistical efficiency of the BTL model, which they refer to as ``topology-dependence''. As a byproduct of analysis for a context-dependent generalization of the Luce choice model, \citet{seshadri2020learning} obtain tight expected risk and tail risk bounds for the MLEs of Luce choice models (which they call MNL) and Plackett--Luce ranking models in terms of the algebraic connectivity, extending and improving upon previous works by \citet{hajek2014minimax,shah2015estimation,vojnovic2016parameter}. Other works with tight risk bounds on the BTL model include \citet{hendrickx2020minimax} and \citet{bong2022generalized}, who also provide the first high probability guarantees for the existence of finite MLEs of the BTL model under conditions on a Fiedler eigenvalue. Interestingly, the statistical significance of algebraic connectivity has also been highlighted in models of networks in econometrics and machine learning by the works of \citet{de2017econometrics,jochmans2019fixed,chang2024inferring}, among others.} Our present work is primarily concerned with the optimization aspects of the maximum likelihood estimation of choice models. Nevertheless, the statistical importance of algebraic connectivity in the aforementioned works also provide motivations for us to leverage it in our convergence analysis of Sinkhorn's algorithm for matrix balancing.  

Lastly, a short note on terminology. Even though a choice model based on \eqref{eq:Luce} is technically a  ``multinomial logit model'' with only intercept terms \citep{mcfadden1973conditional}, there are subtle differences. 
When \eqref{eq:Luce} is applied to model
ranking and choice data with distinct items, each observation $i$ usually consists of a possibly \emph{different}
subset $S_i$ of the universe of all alternatives, so that there is a large number of different configurations of the choice menu in the dataset. On the other hand, common applications of multinomial logit models, such as classification models in statistics and machine learning \citep{bishop2006pattern} and discrete choice models in econometrics \citep{mcfadden1973conditional}, often deal with repeated observations consisting of the \emph{same} number of alternatives. However, these alternatives now possess ``characteristics'' that vary across observations, which are often mapped {parametrically} to the scores in \eqref{eq:Luce}. 
 In this paper, we primarily use the term \emph{Luce choice model} to refer to the model \eqref{eq:model}, although it is also called MNL (for multinomial logit) models in some works. We refrain from using the term MNL to avoid confusion with parametric models for featurized items used in ML and econometrics.

\subsection{Matrix Balancing}
Matrix balancing (or scaling) is an important topic in optimization and numerical linear algebra that underlies a diverse range of applications. The particular question of scaling rows and columns of a matrix $A$ so that the resulting matrix has target row and column norms $p,q$ has been studied as early as the 1930s, and continue to interest researchers from different disciplines today. The present paper only contains a partial survey of the vast literature on this topic. \cref{app:related-works} provides a summary of some popular applications to illustrate the ubiquity of the matrix balancing problem. \citet{schneider1990comparative} and \citet{idel2016review} also provide excellent discussions of many applications. 

The standard iterative algorithm for the matrix balancing problem we study in this paper
has been rediscovered independently quite a few times. As a result, it has domain-dependent names, including the iterative proportional fitting (IPF) procedure \citep{deming1940least}, biproportional fitting \citep{bacharach1965estimating} and the RAS algorithm \citep{stone1962multiple}, but is perhaps most widely known as \emph{Sinkhorn's algorithm} \citep{sinkhorn1964relationship}. A precise description can be found in \cref{alg:scaling}. Sinkhorn's algorithm is also closely related to relaxation and coordinate descent type methods for solving the dual of entropy optimization problems \citep{bregman1967relaxation,cottle1986lagrangean,tseng1987relaxation,luo1992convergence}, as well as message passing and belief propagation algorithms in distributed optimization \citep{balakrishnan2004polynomial,agarwal2018accelerated}. 

The convergence behavior of Sinkhorn's algorithm in different settings has been extensively studied by \citet{sinkhorn1964relationship,bregman1967proof,lamond1981bregman,franklin1989scaling,ruschendorf1995convergence,kalantari2008complexity,knight2008sinkhorn,pukelsheim2009iterative,altschuler2017near,dvurechensky2018computational,marino2020optimal,chakrabarty2021better,leger2021gradient,carlier2022linear}, among many others. For $A$ with strictly \emph{positive} entries, \citet{franklin1989scaling} establish the global linear convergence of Sinkhorn's algorithm in the Hilbert projective metric $d$ \citep{bushell1973hilbert}. More precisely, if $r^{(t)}$ denotes the row sum of the scaled matrix after $t$ iterations of Sinkhorn's algorithm that enforce column constraints, then 
 \begin{align}
 \label{eq:hilbert-contraction} 
   d(r^{(t)}, p) \leq  \lambda^{t} \cdot d(r^{(0)}, p)
 \end{align}
  for some $\lambda \in (0,1)$ dependent on $A$. \revision{On the other hand, works such as \citet{kalantari1993rate}, \citet{altschuler2017near}, and \citet{dvurechensky2018computational} develop complexity bounds on the number of iterations required for the $\ell^1$ distance $\|r^{(t)}- p\|_1 \leq \varepsilon$ for a given $\varepsilon>0$. Although these bounds imply a convergence that is sub-linear, i.e., $\|r^{(t)}- p\|_1 = \mathcal{O}(1/t)$,
  their focus is on optimal dependence on problem size and dimension. An important class of problems with $A>0$ is entropy regularized optimal transport \citep{cuturi2013sinkhorn}, where $A$ is of the form $A=\exp(-c/\gamma)$ with a finite cost function (matrix) $c$, i.e., $A$ strictly positive everywhere. In this setting, convergence of Sinkhorn's algorithm in discrete and continuous problems has been studied by \citet{altschuler2017near,marino2020optimal,leger2021gradient,ghosal2022convergence}, among others. 
  The linear convergence of Sinkhorn's algorithm for $A>0$ has also been extended to the multi-marginal continuous setting by \citet{carlier2022linear}, building on the work of \citet{marino2020optimal}.}
  
  However, the matter of convergence is more delicate when the matrix contains zero entries, and additional assumptions on the problem structure are required to guarantee the existence of scalings $D^0,D^1$ and the convergence of Sinkhorn's algorithm. For non-negative $A$, convergence is first established by \citet{sinkhorn1967concerning} in the special case of {square} $A\geq 0$ and uniform $p=q=\mathbf{1}_n=\mathbf{1}_m$. \revision{Their necessary and sufficient condition is that $A$
has \emph{support}, i.e., there exists a permutation $\sigma$ such that the ``diagonal'' $(A_{1\sigma(1)},A_{2\sigma(2)},\dots,A_{n\sigma(n)})$ is strictly positive. \citet{soules1991rate} and \citet{achilles1993implications} further show that the convergence is linear if and only if the stronger condition of 
\emph{total support} holds, i.e., any non-zero entry of $A$ must be in $(A_{1\sigma(1)},A_{2\sigma(2)},\dots,A_{n\sigma(n)})$
for some permutation $\sigma$. \citet{knight2008sinkhorn} provides a tight \emph{asymptotic} linear convergence rate in terms of the sub-dominant (second largest) singular value of the scaled doubly stochastic matrix $D^0AD^1$. The convergence in \citet{knight2008sinkhorn} is measured by some implicit distance to the target marginals $p,q$. However, no asymptotic linear convergence rate is previously known for non-square $A\geq0$ and non-uniform marginals.}

For general non-negative matrices and non-uniform marginals, the necessary and sufficient conditions on $A$ in the matrix balancing problem that generalize that of \citet{sinkhorn1967concerning} have been studied by \citet{thionet1964note,bacharach1965estimating,brualdi1968convex,menon1968matrix,djokovic1970note,sinkhorn1974diagonal,balakrishnan2004polynomial,pukelsheim2009iterative}, among others, and convergence of Sinkhorn's algorithm under these conditions is well-known. Connecting Sinkhorn's algorithm to dual coordinate descent for entropy optimization, \citet{luo1992convergence} show that the dual optimality gap, defined in \cref{eq:potential convergence}, converges linearly globally with some unknown rate $\lambda$ when finite scalings $D^0,D^1$ exist. However, their result is implicit and there are no results that quantify the global linear rate $\lambda$, even for special classes of non-negative matrices. When convergence results for $A>0$ in previous works are applied to non-negative matrices, the bounds often blow up or become degenerate as soon as $\min_{ij}A_{ij} \downarrow 0$. For example, in \eqref{eq:hilbert-contraction} the contraction factor $\lambda \rightarrow 1$ when $A$ contains zero entries. {When $\min_{ij}A_{ij} = 0$, complexity bounds on $\|r^{(t)}- p\|_1$ and $\|r^{(t)}- p\|_2$ have been established for Sinkhorn's algorithm, for example by the works of \citet{kalantari2008complexity} and \citet{chakrabarty2021better}, with polynomial dependence on $1/\varepsilon$, i.e., sub-linear convergence.} Under the minimal condition that guarantees the convergence of Sinkhorn's algorithm, \citet{leger2021gradient} gives a {quantitative} global sub-linear bound on the KL divergence between $r^{(t)}$ and $p$ in the continuous setting for general probability distributions, which include non-negative matrices $A\geq 0$. 

It therefore remains to reconcile the various results on Sinkhorn's algorithm for $A\geq 0$ and characterize the global and asymptotic linear convergence rates for non-negative $A$. Our results precisely fill these gaps left by previous works. The global linear convergence result in \cref{thm:global-convergence} establishes a contraction like \eqref{eq:hilbert-contraction} for the optimality gap whenever finite scalings $D^0,D^1$ exist, and characterize the convergence rate $\lambda$ in terms of the algebraic connectivity. Moreover, the asymptotic linear rate in \cref{thm:convergence} directly extends the result of \citet{knight2008sinkhorn}. See \cref{tab:convergence-summary} for a detailed summary and comparison of the convergence results in previous works and this paper. The dependence of Sinkhorn's convergence rate on spectral properties of graphs can be compared to convergence results in the  literature on decentralized optimization and gossip algorithms, where a spectral gap quantifies the convergence rate \citep{boyd2006randomized,xiao2007distributed}.

\revision{Lastly, we note that other algorithms with better complexities have been developed for the matrix balancing problem, utilizing e.g., the ellipsoid algorithm \citep{kalantari1996complexity} and geometric programming \citep{nemirovski1999complexity}, interior point algorithms \citep{cohen2017matrix,chen2022maximum}, or customized first/second order techniques \citep{linial1998deterministic,allen2017much}. However, despite having better theoretical complexities, most of these algorithms are yet to be implemented practically.  Sinkhorn's algorithm, on the other hand,  remains an attractive choice in practice due to its simplicity, robustness, and parallelization capabilities.}

%% file: preliminaries.tex
\section{Preliminaries on Choice Modeling and Matrix Balancing}
\label{sec:formulations}
We start by providing brief but self-contained introductions to the two main subjects of this paper,
choice modeling and matrix balancing, including their respective underlying mathematical problems and assumptions. Then, we formally establish their equivalence in \cref{sec:equivalence}.
\subsection{Maximum Likelihood Estimation of Luce Choice Models}
\label{subsec:Luce}
In the Luce choice modeling framework, we have $n$ observations $\{(j_i,S_i)\}_{i=1}^n$, each consisting of a choice set $S_{i}\subseteq \{1,\dots,m\}=[m]$ that is a subset
of the total $m$ alternatives/items/objects, and the alternative selected, denoted by $j_i \in S_i$.
 The choice probability is prescribed by Luce's axiom of choice given model parameter $s\in\mathbb{R}^m_{++}$
 in the interior of the probability simplex $\Delta_m$:
\begin{align*}
\Pr(j_i, S_i) & =\frac{s_{j_i}}{\sum_{k\in S_{i}}s_{k}},
\end{align*}
and the likelihood of the observed data is thus given by 
\begin{align}
\label{eq:model}
L(s;\{(j_i,S_i)\}_{i=1}^n):=& \prod_{i=1}^{n}\frac{s_{j_i}}{\sum_{k\in S_{i}}s_{k}}.
\end{align} 
A popular method to estimate $s = \{s_{1},\dots,s_m\}$ is the maximum likelihood estimation approach, which maximizes the log-likelihood 
\begin{align}
\label{eq:log-likelihood}
\ell(s):=\log L(s) & =\sum_{i=1}^{n}\left(\log s_{j_i}-\log\sum_{k\in S_{i}}s_{k}\right).
\end{align}
over the interior of the probability simplex. \revision{Note that the choice sets $S_i$ can vary across $i$. In other words, in each observation the choice is made from a potentially distinct set of alternatives. This feature of the problem turns out to be important for both the algorithmic efficiency of computing the maximizers to \eqref{eq:log-likelihood}, as well as the statistical efficiency of the resulting maximum likelihood estimators (MLEs), which can be quantified by a measure of \emph{connectivity} of the data structure. We will elaborate on these points shortly. For now, we focus on the existence and uniqueness of MLE.}  %Importantly, the fact that \eqref{eq:model} has potentially different choice sets $S_i$ across observations have structure-dependent implications for both the optimization and the statistical efficiencies of the model.

If we reparameterize $\exp(u_{j})=s_{j}$, it is obvious that \eqref{eq:log-likelihood} is concave in $u$. However, to ensure the log-likelihood  \eqref{eq:log-likelihood} has a unique maximizer in the \emph{interior} of the simplex, additional assumptions on the comparison structure of the dataset $\{(j_i,S_i)\}_{i=1}^n$ are needed. 
The following classic condition is necessary and sufficient for the maximum likelihood problem to be well-posed.
\begin{assumption}[\textbf{Strong Connectivity}]
\label{ass:strong-connected}
In any partition of $[m]$ into two nonempty subsets $S$ and its complement $S^C$,
some $j \in S$ is \emph{selected} at least once over some
$k \in S^C$. Equivalently, the {\it directed comparison graph}, with items as vertices and an edge $j\rightarrow k$ if and only if $k$ is selected in some $S_i$ for which $j,k\in S_i$, is strongly connected. 
\end{assumption}
\cref{ass:strong-connected} is standard in the literature \citep{hunter2004mm,noothigattu2020axioms} and appeared as early as the work of \citet{zermelo1929berechnung} and \citet{ford1957solution} for pairwise comparisons. \citet{hunter2004mm} shows that \cref{ass:strong-connected} is necessary and sufficient for the upper-compactness
of \eqref{eq:log-likelihood}, which guarantees the existence of a maximizer in the interior of the probability simplex. In fact, when an interior maximizer exists, it is also \emph{unique}, since \cref{ass:strong-connected} implies the following weaker condition, which guarantees the strict concavity of \eqref{eq:log-likelihood}.
\begin{assumption}[\textbf{Connectivity}]
\label{ass:weak-connected}
In any partition of $[m]$ into two nonempty subsets $S$ and $S^C$,
some $j\in S$ and some $k\in S^C$ \emph{appear} in the same choice set $S_i$ for some $i$. 
\end{assumption}
The intuitions provided by \citet{ford1957solution} and \citet{hunter2004mm} are helpful for understanding Assumptions \ref{ass:strong-connected} and \ref{ass:weak-connected}. If items from some $S\subsetneq [m]$ are never compared with those in $S^C$, i.e., never appeared together in any choice set $S_i$, it is impossible to rank across the two subsets. In this case, we can rescale the relative weights of $S$ and $S^C$ of an interior maximizer and obtain another maximizer. On the other hand, if items in $S$ are always preferred to those in $S^C$, we can increase the likelihood by scaling $s_j$ for items $j\in S^C$ towards $0$, and no maximizer in the \emph{interior} of the probability simplex exists. 
Nevertheless, a boundary solution can still exist. This case turns out to be important in the present work: in the equivalent matrix balancing problem, it corresponds to the  \emph{slowdown} regime of Sinkhorn's algorithm, where scalings diverge but the scaled matrix converges (\cref{sec:linear-convergence}). 

\cref{ass:weak-connected} also has a concise graph-theoretic interpretation. Define the weighted {undirected comparison graph} $G_c$ on $m$ vertices with adjacency matrix $A^c$ given by 
\begin{align}
\label{eq:comparison-graph-adjacency}
{A}^c_{jk} & =\begin{cases}
0 & j=k\\
|\{i\in [n]\mid j,k\in S_i\}| & j\neq k,
\end{cases}
\end{align}
In other words, there is an undirected edge between $j$ and $k$ if and only if they are both included in some choice set $S_i$, with the edge weight equal to the number of their co-occurrences, which could be zero. We can verify that \cref{ass:weak-connected} precisely requires $G_c$ to be connected.
\revision{
\begin{remark}[\textbf{Importance of Graph Connectivity}]
    Under the standard Assumptions \ref{ass:strong-connected} and \ref{ass:weak-connected}, previous works have studied the statistical efficiency of the MLE \citep{hajek2014minimax,shah2015estimation,seshadri2020learning} as well as the computational efficiency of the MM algorithm for computing the MLE \citep{vojnovic2020convergence}. In both cases, the \emph{algebraic connectivity} of $G_c$ \citep{fiedler1973algebraic}, quantified by the second smallest eigenvalue of the graph Laplacian of $G_c$, plays an important role. See \cref{subsec:graph-laplacian} for more details. The importance of spectral properties for parameter learning in data with graph or matrix structures has appeared as early as \citet{kendall1940method} and in the classic work of \citet{keener1993perron} on ranking sports teams, as well as in works in economics \citep{abowd1999computing,jochmans2019fixed}. These results, together with
    the connections we establish in this paper between choice modeling and matrix balancing, 
    inspire us to quantify the convergence of Sinkhorn's algorithm also using the algebraic connectivity of a \emph{bipartite} graph, defined in \eqref{eq:bipartite-adjacency}. 
\end{remark}
}
\subsection{The Canonical Matrix Balancing Problem}
\label{subsec:matrix-balancing}
Matrix balancing is a classic problem that shows up in a wide range of disciplines. See \cref{app:related-works} for a short survey on some applications. 
The underlying mathematical problem can be stated concisely in matrix form as:
\begin{quote} 
Given positive vectors $p \in \mathbb{R}_{++}^n,q \in \mathbb{R}_{++}^m$with $\sum_i p_i=\sum_j q_j=c>0$, which without loss of generality can be set to $c=1$, and a non-negative matrix $A\in \mathbb{R}_+^{n\times m}$, find positive diagonal matrices 
$D^{1}$, $D^{0}$ satisfying the conditions $D^{1}AD^{0}\cdot\mathbf{1}_m=p$
and $D^{0}A^{T}D^{1}\cdot\mathbf{1}_n=q$.\end{quote}
We henceforth refer to the above as the ``canonical'' matrix balancing problem. %It seeks positive row and column scalings of an (entry-wise) non-negative rectangular matrix $A$ such that the scaled matrix has positive target row and column sums $p$ and $q$. 
Other variants of the problem replace the row and column sums (the 1-norm) with other norms \citep{bauer1963optimally,ruiz2001scaling}. Note that for any $c>0$, $(D^0/c,cD^1)$ 
 is also a solution whenever $(D^0,D^1)$ is. A finite positive solution $(D^{0},D^{1})$ to the canonical matrix balancing problem is often called a \emph{direct scaling}.

 The structure of the matrix balancing problem suggests a simple iterative scheme: starting
from any initial positive diagonal $D^{0}$, invert $D^{1}AD^{0}\mathbf{1}_m = p$ using $p/(AD^{0}\mathbf{1}_m)$ 
 to update $D^{1}$. Then invert $D^{0}A^{T}D^{1}\mathbf{1}_n=q$ using  $q/(A^{T}D^{1}\mathbf{1}_n)$ to compute the new estimate of $D^{0}$, and repeat the procedure, leading to a solution if it converges. Here, divisions involving two vectors of the same length are \emph{entry-wise}. This simple iterative scheme is precisely Sinkhorn's algorithm, described in \cref{alg:scaling}, where vectors $d^0,d^1$ are the diagonal elements of $D^0,D^1$.

An important dichotomy occurs depending on whether the entries of $A$ are strictly positive. If $A$ contains no zero entries, then direct scalings and a unique scaled matrix $D^1AD^0$ always exist \citep{sinkhorn1964relationship}. Moreover, Sinkhorn's algorithm converges linearly \citep{franklin1989scaling}. 
On the other hand, when $A$ contains zero entries, the problem becomes more complicated. Additional conditions are needed to guarantee meaningful solutions, and the convergence behavior of Sinkhorn's algorithm is less clearly understood.  Well-posedness of the matrix balancing problem has been studied by \citet{brualdi1968convex,sinkhorn1974diagonal,pukelsheim2009iterative}, among others, who characterize the following equivalent existence conditions.
\begin{assumption}[\textbf{Strong Existence}]
\label{ass:matrix-existence}
\textbf{(a)} There exists a non-negative matrix $A'\in \mathbb{R}_+^{n\times m}$ with the same zero patterns as $A$ and with row and column sums $p$ and $q$. Or, equivalently,

\textbf{(b)} For every pair of sets of indices $N \subsetneq [n]$ and $M \subsetneq [m]$ such that $A_{ij}=0$ for $i\notin N$ and $j\in M$, $\sum_{i\in N}p_i \geq \sum_{j\in M}q_j$, with equality if and only if $A_{ij}=0$ for all $i \in N$ and $j \notin M$ as well.
\end{assumption} 
It is well-known in the matrix balancing literature that the above two conditions are equivalent, and that a positive finite solution $(D^{0},D^{1})$ to the canonical problem exists if and only if they hold. See, for example, Theorem 6 in \citet{pukelsheim2009iterative}. \revision{\cref{ass:matrix-existence} also guarantees the convergence of Sinkhorn's algorithm. However, it is not a necessary condition. In other words, Sinkhorn's algorithm could converge even if the matrix balancing problem does not admit a direct scaling. This phenomenon turns out to be important in characterizing the convergence rate, which we study in \cref{sec:linear-convergence}.}

Clearly, \cref{ass:matrix-existence}(a) is the minimal necessary condition when a solution to the matrix balancing problem exists, and trivially holds when $A>0$ (take for example $A'$ as the Kronecker product of $p,q$). \cref{ass:matrix-existence}(b) is closely connected to conditions for perfect matchings in bipartite graphs \citep{hall1935representatives,galichon2021matching}. In flow networks \citep{gale1957theorem,ford1956maximal,ford1957simple}, it is a capacity constraint that
guarantees the maximum flow on a weighted bipartite graph is equal to $\sum_i p_i=\sum_j q_j$ and with positive flow on every edge \citep{idel2016review}. \revision{The weighted bipartite graph, denoted by $G_b$, is important in this paper. Its adjacency matrix $A^b \in \mathbb{R}^{(n+m)\times(n+m)}$ can be represented concisely using $A$ as
\begin{align}
\label{eq:bipartite-adjacency}
    A^b := \begin{bmatrix}\mathbf{0} & {A}\\
{A}^{T} & \mathbf{0}
\end{bmatrix},	
\end{align}
and $A$ is sometimes called the \emph{biadjacency} matrix of $G_b$. See \cref{subsec:graph-laplacian} for more information. Just like in the choice setting, where the connectivity of the undirected comparison graph $G_c$ plays an important role, the connectivity of $G_b$ turns out to be important for the linear convergence rate of Sinkhorn's algorithm (see \cref{sec:linear-convergence}).}

Lastly, the necessary and sufficient condition for the uniqueness of finite scalings essentially requires that $A$ is not block-diagonal, and precisely guarantees that $G_b$ is \emph{connected}.
\begin{assumption}[\textbf{Uniqueness}]
\label{ass:matrix-uniqueness}
 $D^{0}$ and $D^{1}$ are unique modulo normalization if and only if $A$ is indecomposable, i.e., there does not exist permutation matrices $P,Q$ such that $PAQ$ is block diagonal.
\end{assumption} 
 \begin{algorithm}[tb]
\caption{Sinkhorn's Algorithm}
   \label{alg:scaling}
\begin{algorithmic}
   \STATE {\bfseries Input:}  $A, p, q,\epsilon_{\text{tol}}$.
   \STATE {\bfseries initialize} $d^{0}\in\mathbb{R}_{++}^{m}$
   \REPEAT
   \STATE $d^{1} \leftarrow  p/( A d^0)$ 
   \STATE $d^{0}\leftarrow  q/({A}^{T} d^{1})$
   \STATE 
$\epsilon\leftarrow$ maximal update in $(d^{0},d^1)$ or distance between $D^1Ad^0$ and $p$
\UNTIL{$\epsilon<\epsilon_{\text{tol}}$}
\end{algorithmic}
\end{algorithm}
With a proper introduction to both problems, we are now ready to establish the equivalence between Luce choice model estimation and matrix balancing. In \cref{sec:linear-convergence}, we return to Sinkhorn's algorithm for the matrix balancing problem, and provide answers to open problems concerning its linear convergence for non-negative $A$, by leveraging the connections we establish next. 
%In \cref{sec:connections}, we further discuss the extensive connections of matrix balancing and Sinkhorn's algorithm to choice modeling and optimization.

%% file: equivalence.tex
\section{Connecting Choice Modeling and Matrix Balancing}
\label{sec:equivalence}

In this section, we formally establish the connections between choice modeling and matrix balancing. We show that
maximizing the log-likelihood \eqref{eq:log-likelihood} is equivalent to solving a canonical matrix balancing problem. We also precisely describe the correspondence between the relevant conditions in the two problems. In view of this equivalence, we show that Sinkhorn's algorithm, when applied to estimate Luce choice models, is in fact a \emph{parallelized} generalization of the classic iterative algorithm for choice models, dating back to \citet{zermelo1929berechnung,dykstra1956note,ford1957solution}, and studied extensively also by \citet{hunter2004mm,vojnovic2020convergence,vojnovic2023accelerated}.

\subsection{Maximum Likelihood Estimation of Luce Choice Models as Matrix Balancing}
\label{subsec:reformulation}
The optimality conditions for maximizing the log-likelihood \eqref{eq:log-likelihood} for each $s_j$ are given by
\begin{align*}
\partial_{s_{j}}\ell(s)=\sum_{i\in[n] \mid (j,S_i)}\frac{1}{s_{j}}-\sum_{i\in[n]\mid j\in S_{i}}\frac{1}{\sum_{k\in S_{i}}s_{k}} & =0.
\end{align*}
 Multiplying by $s_{j}$ and dividing by $1/n$, we have 
\begin{align}
\label{eq:optimality-original}
\frac{W_j}{n} & = \frac{1}{n} \sum_{i\in [n]\mid j\in S_{i}}\frac{s_{j}}{\sum_{k\in S_{i}}s_{k}},
\end{align}
where $W_j:=|\{i\in[n]\mid (j,S_i)\}|$ is the number of observations where $j$ is selected.

Note that in the special case where $S_i\equiv [n]$, i.e., every choice set contains \emph{all} items, the MLE simply reduces to the familiar empirical frequencies $\hat s_j = {W_j}/{n}$. However, when the choice sets $S_i$ vary, no closed form solution to \eqref{eq:optimality-original} exists, which is the primary motivation behind the long line of works on the algorithmic problem of solving \eqref{eq:optimality-original}. With varying $S_i$, we can interpret the optimality condition as requiring the \emph{observed} frequency of $j$ being chosen (left hand side) be equal to the average or \emph{expected} probability of $j$ being selected (right hand side), which conditional on choice set $S_i$ is $\frac{s_{j}}{\sum_{k\in S_{i}}s_{k}}$ if $j\in S_i$ and 0 otherwise. In addition, note that since the optimality condition in \eqref{eq:optimality-original} only involves the \emph{frequency} of selection, distinct datasets could yield the same optimality conditions and hence the same MLEs. For example, suppose that two alternatives $j$ and $k$ both appear in choice sets
$S_{i}$ and $S_{i'}$, with $j$ selected in $S_{i}$ and  $k$ selected in $S_{i'}$. Then switching
the choices in $S_{i}$ and $S_{i'}$ does not alter the likelihood and optimality conditions. This feature holds more generally with longer cycles of items and choice sets, and can be viewed as a consequence of the context-independent nature of Luce's choice axiom, i.e., IIA. In some sense, it is also the underpinning of many works in economics that estimate choice models based on \emph{marginal} sufficient statistics. A prominent example is  \citet{berry1995automobile}, which estimates consumer preferences using data on {aggregate} market shares of products.

\begin{remark}[\revision{\textbf{Reduction to Unique Choice Sets}}]
    In practice, the choice sets of many observations may be identical to each other, i.e., $S_i\equiv S_{i'}$ for some $i,i'\in [n]$. Because \eqref{eq:optimality-original} only depends on the total ``winning'' counts of items, we may aggregate over observations with the same $S_i$:
\begin{align*}
\sum_{i\in[n]\mid j\in S_{i}}\frac{s_{j}}{\sum_{k\in S_{i}}s_{k}} & = \sum_{i'\in[n^\ast]\mid j\in S^\ast_{i'}} R_{i'} \cdot \frac{s_{j}}{\sum_{k\in S^\ast_{i'}}s_{k}},
\end{align*}
where each $S_{i'}^\ast$ is a unique choice set that appears in $R_{i'}\geq 1$ observations, for $i'=1,\dots,n^\ast \leq n$. By construction, $\sum_{i'=1}^{n^\ast}R_{i'} =n$. Note, however, that the selected item could vary across different appearances of $S_i^\ast$, yet the optimality conditions only involve each item's winning count $W_j$. From now on, we will assume this reduction and drop the $^\ast$ superscript. In other words, without loss of generality, we assume that we observe $n$ \emph{unique} choice sets, and choice set $S_i$ has \emph{multiplicity} $R_i$ \citep{shah2015estimation}. The resulting maximum likelihood problem has optimality conditions
\begin{align}
\label{eq:optimality}
{W_j} & = \sum_{i\in[n]\mid j\in S_{i}}{R_i} \cdot \frac{s_{j}}{\sum_{k\in S_{i}}s_{k}}.
\end{align}
\end{remark}

We are now ready to reformulate \eqref{eq:optimality} as a canonical matrix balancing problem. Define $p\in\mathbb{R}^{n}$ as $p_i=R_i$, i.e., the number of times choice set $S_i$ appears in the data. Define
$q\in\mathbb{R}^{m}$ as $q_j={W_j}$,
i.e., the number of times item $j$ was \emph{selected} in the data. By construction we have $\sum_i p_i=\sum_j q_j$, and $p_i,q_j>0$ whenever \cref{ass:strong-connected} holds.

\begin{figure}
    \centering
\includegraphics[width=0.35\linewidth]{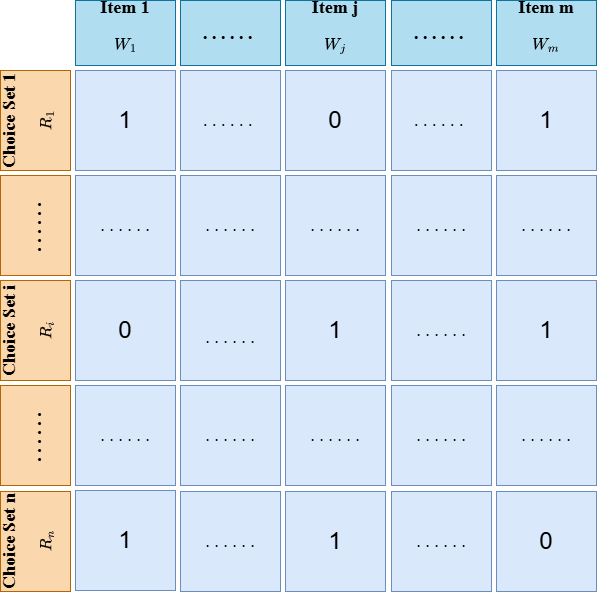}
    \caption{Representation of Luce choice data with participation matrix $A$, $R_i$ the frequency of appearances of choice set $i$, and $W_j$ the frequency of choices of item $j$ as a matrix balancing problem $(A,p,q)$ with target marginals $p_i=R_i$ and $q_j=W_j$.} \label{fig:matrix_balancing}
\end{figure}

Now define the $n\times m$ binary matrix $A$ by
$A_{ij}=1\{j\in S_{i}\}$, so the $i$-th row of $A$
is the indicator of which items appear in the (unique) choice set $S_i$, and the $j$-th column of $A$ is the indicator of which choice sets
item $j$ appears in. We refer to this $A$ constructed from a choice dataset as the \emph{participation matrix}. By construction, $A$ has distinct rows, but may still have identical columns. We may also remove repeated columns by ``merging'' items and their win counts. Their estimated scores can be computed from the score of the merged item proportional to their respective win counts. We do not require this reduction in our results. \cref{fig:matrix_balancing} provides an illustration of the matrix balancing representation of the luce choice modeling problem with $(A,p,q)$ defined as above.

Let $D^{0}\in\mathbb{R}^{m\times m}$ be the diagonal matrix with
$D_{j}^{0}=s_{j}$ and $D^{1}\in\mathbb{R}^{n\times n}$ be the
diagonal matrix with $D_{i}^{1}={R_i}/{\sum_{k\in S_{i}}s_{k}}$,
and define the scaled matrix
\begin{align}
\label{eq:scaled-matrix}
\hat{A} & :=D^{1}AD^{0}.
\end{align}
The matrices $D^{1}$ and $D^{0}$ are scalings of rows and columns
of $A$, respectively, and
\begin{align*}
    \hat{A}_{ij} = \frac{R_i}{\sum_{k\in S_{i}}s_{k}}\cdot1\{j\in S_{i}\}\cdot s_{j}.
\end{align*}
The key observation is that the optimality conditions \eqref{eq:optimality} can be rewritten as
\begin{align}
\label{eq:bridge}
\hat{A}^T \mathbf{1}_n & = q.
\end{align}
Moreover, by construction $\hat{A}$ also satisfies
\begin{align}
\label{eq:marginal}
\hat{A}\mathbf{1}_m & = p.
\end{align}
 Therefore, if $s_j$'s satisfy the optimality conditions for maximizing \eqref{eq:log-likelihood}, then $D^0,D^1$  defined above solve the matrix balancing problem in \cref{eq:scaled-matrix,eq:bridge,eq:marginal}. Moreover, the converse is also true, and we thus establish the equivalence between choice maximum likelihood estimation and matrix balancing. All omitted proofs appear in \Cref{app:proofs}.
\begin{theorem}[\revision{Equivalence of Problems}]
\label{prop:mle-scaling}
\revision{Let $(A,p,q)$ be constructed from a choice dataset as follows: $p\in\mathbb{R}^{n}$ with $p_{i}=R_i$, the multiplicity of choice set $S_i$; $q\in\mathbb{R}^{m}$ with $q_{j}=W_j$, the total number of times item $j$ is chosen in the choice dataset; 
$A\in\mathbb{R}^{n\times m}$
with $A_{ij}=1\{j\in S_{i}\}$, i.e., the $i$-th row of $A$ is a one-hot encoding of the choice set $S_i$.} 

Then  $D^{0},D^{1}>0$ with $\sum_j D_j^0=1$ solves the matrix balancing problem
\begin{align}
\label{eq:equation-system}
\begin{split}
D^{1}AD^{0} & =\hat{A}\\
\hat{A}\mathbf{1}_m & =p\\
\hat{A}^{T} \mathbf{1}_n & =q
\end{split}
\end{align}
if and only if $s \in \Delta_m$ with $s_j=D^0_j$ satisfies the optimality conditions \eqref{eq:optimality} of the maximum likelihood estimation problem of a Luce choice model given the choice dataset.
\end{theorem}
\cref{prop:mle-scaling} implies that \eqref{eq:log-likelihood} has a unique maximizer $s$ in the interior of the probability simplex if and only if \eqref{eq:equation-system} has a unique normalized solution $D^0$ as well. The next question, naturally, is then how \cref{ass:strong-connected} and \cref{ass:weak-connected} for choice modeling are connected to \cref{ass:matrix-existence} and \cref{ass:matrix-uniqueness} for matrix balancing.
\begin{proposition}
[\revision{Equivalence of Assumptions}]
\label{thm:necessary-and-sufficient}
Let (A,p,q) be constructed from the choice dataset as in \cref{prop:mle-scaling}, with $p,q>0$. \cref{ass:weak-connected} is equivalent to \cref{ass:matrix-uniqueness}. Furthermore, \cref{ass:strong-connected} holds if and only if $(A,p,q)$ satisfy \cref{ass:matrix-existence} and $A$ satisfies \cref{ass:matrix-uniqueness}.
\end{proposition} 
Thus when the choice maximum likelihood estimation problem is cast as a matrix balancing problem, \cref{ass:matrix-existence} exactly characterizes the \emph{gap} between \cref{ass:weak-connected} and \cref{ass:strong-connected}. 
 We provide some intuition for \cref{thm:necessary-and-sufficient}. When we construct a triplet $(A,p,q)$ from a choice dataset, with $p$ the numbers of appearances of unique choice sets and $q$ the winning counts of each item, \cref{ass:matrix-uniqueness} precludes the possibility of partitioning the items into two subsets that never get compared with each other, i.e., \cref{ass:weak-connected}. Then \cref{ass:matrix-existence} requires that whenever a strict subset $M\subsetneq [m]$ of objects only appear in a strict subset $N\subsetneq [n]$ of the observations, their total winning counts are \emph{strictly} smaller than the total number of these observations, i.e., there is some object $j\notin M$ that is selected in $S_i$ for some $i\in N$, which is required by \cref{ass:strong-connected}.

Interestingly, while \cref{ass:strong-connected} requires the directed comparison graph, defined by the $m\times m$ matrix of counts of item $j$ being chosen over item $k$, to be strongly connected, the corresponding conditions for the equivalent matrix balancing problem concern the $n\times m$ participation matrix $A$ and positive vectors $p,q$, which do not explicitly encode the specific \emph{choice} of each observation. This apparent discrepancy is due to the fact that $(A,p,q)$ form the \emph{sufficient statistics} of the Luce choice model. In other words, there can be more than one choice dataset with the same optimality condition \eqref{eq:optimality} and $(A,p,q)$ defining the equivalent matrix balancing problem (see \cref{fig:matrix_balancing}).

\begin{remark}[\revision{\textbf{Aggregate Data as Sufficient Statistics}}]
     The feature where ``marginal'' or aggregate quantities constitute the sufficient statistics of a parametric model is an important characteristic that underlies many works in economics and statistics \citep{kullback1997information,stone1962multiple,good1963maximum,birch1963maximum,theil1967economics,fienberg1970iterative,berry1995automobile,fofana2002balancing,maystre2017choicerank,chang2024inferring}. It makes the task of estimating a \emph{joint} model from marginal quantities feasible, and is very useful because in many applications, only marginal data is available due to high sampling cost or privacy-preserving considerations.
\end{remark}

Having formulated a particular matrix balancing problem from the estimation problem given choice data, we may ask how one can go in the other direction. In other words, when/how can we construct a ``choice dataset'' whose sufficient statistics is a given triplet $(A,p,q)$? First off, for $(A,p,q)$ to be valid sufficient statistics of a Luce choice model, $p,q$ need to be positive integers. Moreover, $A$ must be a binary matrix with unique rows, each containing at least two non-zero elements (valid choice sets have at least two items). Given such a $(A,p,q)$ satisfying Assumptions \ref{ass:matrix-existence} and \ref{ass:matrix-uniqueness}, a choice dataset can be constructed efficiently. Such a procedure is described, for example, in \citet{kumar2015inverting}, where $A$ is motivated
by random walks on a graph instead of matrix balancing (\cref{sec:connections}). Their construction relies on finding the max flow on the bipartite graph $G_b$. For rational $p,q$, this max flow can be found efficiently in polynomial time \citep{balakrishnan2004polynomial,idel2016review}. Moreover, the maximum flow implies a matrix $A'$ satisfying \cref{ass:matrix-existence}(a), thus providing a feasibility certificate for the matrix balancing problem as well. 

We have thus closed the loop and fully established the equivalence of the maximum likelihood estimation of Luce choice models and the canonical matrix balancing problem.
\begin{corollary}
    There is a one-to-one correspondence between classes of maximum likelihood estimation problems with the same optimality conditions \eqref{eq:optimality} and canonical matrix balancing problems with $(A,p,q)$, where $A$ is a valid binary participation matrix and $p,q>0$ have integer entries. 
\end{corollary}
\revision{Connections to discrete choice modeling have also been established for the related problem of regularized semi-discrete optimal transport \citep{tacskesen2023semi}, although the problem setting and results are distinct from the ones studied in this paper.} We next turn our attention to the algorithmic connections between choice modeling and matrix balancing. 
\subsection{Algorithmic Connections between Matrix Balancing and Choice Modeling}
\label{subsec:IPF}
 Given the equivalence between matrix balancing and choice modeling, we can naturally consider applying Sinkhorn's algorithm to maximize \eqref{eq:log-likelihood}. In this case, one can verify that the updates in each full iteration of \cref{alg:scaling} reduce algebraically to
\begin{align}
\label{eq:scaling-iteration}
s_{j}^{(t+1)} & =W_{j}/\sum_{i\in[n]\mid j\in S_{i}}\frac{R_i}{\sum_{k\in S_{i}}s_{k}^{(t)}}
\end{align}
 in the $t$-th iteration. 
 Comparing \eqref{eq:scaling-iteration} to the optimality condition in \eqref{eq:optimality}, which we recall is given by 
\begin{align*}
 W_j & =\sum_{i\in[n]\mid j\in S_{i}} R_i\frac{s_{j}}{\sum_{k\in S_{i}}s_{k}}= s_{j} \cdot \sum_{i\in[n]\mid j\in S_{i}}\frac{R_i}{\sum_{k\in S_{i}}s_{k}},
 \end{align*}
we can therefore interpret the iterations as simply dividing the winning count $W_j$ by the coefficient of $s_j$ on the right repeatedly, in the hope of converging to a \emph{fixed point}. A similar intuition was given by \citet{ford1957solution} in the special case of pairwise comparisons. Indeed, the algorithm proposed in that paper is a cyclic variant of \eqref{eq:scaling-iteration} applied to pairwise comparisons. However, this connection is mainly algebraic, as the optimality condition in \citet{ford1957solution} does not admit a reformulation as the matrix balancing problem in \eqref{eq:equation-system}.

In \cref{sec:connections}, we provide further discussions on the connections of Sinkhorn's algorithm to existing frameworks and algorithms in the choice modeling literature and the optimization literature. We demonstrate that many existing algorithms for Luce choice model estimation are in fact special cases or analogs of Sinkhorn's algorithm. These connections also illustrate the many interpretations of Sinkhorn's algorithm, e.g., as a distributed optimization algorithm as well as a minorization-maximization or majorization-minimization (MM) algorithm \citep{lange2016mm}. However, compared to most algorithms for choice modeling discussed in this work, Sinkhorn's algorithm is more general as it applies to non-binary $A$ and non-integer $p,q$, and has the additional advantage of being paralellized and distributed, hence more efficient in practice. \revision{To highlight the ubiquity of Sinkhorn's algorithm in the choice setting, we summarize these algorithmic connections below.
\begin{theorem}[\textbf{Equivalence of Algorithms}]
\label{thm:algorithm-equivalence}
    Sinkhorn's algorithm, when applied to matrix balancing formulations of various choice modeling problems, is equivalent to the following algorithms:
    \begin{itemize}
        \item The iterative algorithms of \citet{zermelo1929berechnung,dykstra1956note,ford1957solution} for the BTL model of pairwise comparison data;
        \item The MM algorithm of \citet{hunter2004mm} for the Plackett-Luce model of ranking data;
        \item The unregularized version of the ChoiceRank algorithm of \citet{maystre2017choicerank} for their proposed network choice model;
        \item The BLP algorithm of \citet{berry1995automobile} in a logit random utility model with only intercepts.
    \end{itemize}
    Moreover, Sinkhorn's algorithm can be viewed as a general MM algorithm as well as a message passing algorithm, and in the latter case it is a variant of the accelerated spectral ranking algorithm for Luce choice models of \cite{agarwal2018accelerated} based on a different moment condition.
\end{theorem}
}
The mathematical and algorithmic connections between matrix balancing and choice modeling we establish in this paper allow the transmission of ideas in both directions. For example, inspired by regularized maximum likelihood estimation \citep{maystre2017choicerank}, we propose a regularized version of Sinkhorn's algorithm in \cref{subsec:regularization}, which is guaranteed to converge even when the original Sinkhorn's algorithm does not converge. \revision{Leveraging optimization connections between maximum likelihood estimation of choice models and matrix balancing, \citet{chang2024inferring} propose a statistical framework for network traffic that justifies the popular use of Sinkhorn's algorithm to infer detailed dynamic networks from aggregate node-level activities.} In the rest of the main text, we focus on resolving some interesting open problems on the convergence of Sinkhorn's algorithm, motivated by results in choice modeling on the importance of algebraic connectivity in quantifying statistical and computational efficiencies. 

%% file: convergence-analyses.tex
\section{Linear Convergence of Sinkhorn's Algorithm for Non-negative Matrices}
\label{sec:linear-convergence}
In this section, we turn our attention on matrix balancing and study the global and asymptotic linear convergence of Sinkhorn's algorithm for general non-negative matrices $A\geq 0$ and positive marginals $p,q>0$. We first present the relevant optimization principles behind matrix balancing and discuss some existing results in \cref{subsec:convergence-background}, and then present the convergence results in Sections \ref{subsec:global-linear-convergence}--\ref{subsec:sharp-rate}.
Throughout, we use superscript $(t)$ to denote quantities after $t$ iterations of Sinkhorn's algorithm with normalized columns, described in \cref{alg:scaling}.  

\subsection{Preliminaries}
\label{subsec:convergence-background}
We start with the optimization principles associated with matrix balancing and Sinkhorn's algorithm. \revision{Given a matrix balancing problem with $A\geq 0,\sum_{ij}A_{ij}=1$ and target marginals $p,q$ with $\sum_ip_i=\sum_jq_j=1$, consider the following KL divergence (relative entropy) minimization problem 
     \begin{align}
\label{eq:relative-entropy-minimization}
\begin{split}
  \min_{\hat A\in \mathbb{R}^{n\times m}_+} & D_{\text{KL}}(\hat{A}\| A)\\
\hat{A}\mathbf{1}_m & =p\\
\hat{A}^{T}\mathbf{1}_n & =q.
\end{split}
\end{align}
It is well-known that when scalings $D^0,D^1$ solve the matrix balancing problem with $(A,p,q)$, the scaled matrix $\hat A=D^1AD^0$ is the unique minimizer of \eqref{eq:relative-entropy-minimization} \citep{bregman1967proof,ireland1968contingency}. Moreover, vector representations $d^0,d^1$ of the optimal scalings $D^0,D^1$ precisely minimize the following (negative) dual objective of \eqref{eq:relative-entropy-minimization}:
 \begin{align}
 \label{eq:log-barrier}
     g(d^0,d^1)	:=(d^1)^{T}Ad^0-\sum_{i=1}^{n}p_{i}\log d^1_{i}-\sum_{j=1}^{m}q_{j}\log d^0_{j},
\end{align} 
and Sinkhorn's algorithm is a block coordinate descent (BCD) type algorithm \citep{tseng2001convergence} applied to minimize \eqref{eq:log-barrier}. In \citet{luo1992convergence}, the authors study the linear convergence of block coordinate descent algorithms for a general class of objectives that includes \eqref{eq:log-barrier}. In particular, their result implies that the
convergence of Sinkhorn's algorithm, measured in terms of the optimality gap of $g$, is linear with some unknown rate,
 as long as finite positive scalings $D^0,D^1$ exist that satisfy \eqref{eq:matrix-balancing}.} The function $g$, sometimes referred to as the \emph{potential function} of Sinkhorn's algorithm, also turns out to be crucial in quantifying the global linear convergence rate in the present work. 
\begin{remark}[\textbf{Optimization Connections}]
   Interestingly, minimizing \eqref{eq:log-barrier} is in fact equivalent to maximizing the log-likelihood function $\ell(s)$ in \eqref{eq:log-likelihood} for valid $(A,p,q)$, because $\min_{d^1}g(d^0,d^1)=-\ell(d^0)+c$ for some $c>0$. Moreover, the optimality condition of minimizing $g$ with respect to $d^0$ reduces to the optimality condition \eqref{eq:optimality}. A detailed discussion can be found in \cref{sec:sinkhorn-MM}. This connection relates choice modeling and matrix balancing from an optimization perspective. 
\end{remark}

      \revision{Although convergence results on Sinkhorn's algorithm are abundant, they are often stated in different forms, developed under different assumptions, and apply to settings of varying degrees of generality. Next, we briefly discuss some existing works to clarify their connections and distinctions, which also helps motivate the technical results in this paper. They are summarized in \cref{tab:convergence-summary}.} 

      \revision{First and foremost, how to define and measure the convergence of Sinkhorn's algorithm is not entirely trivial. Because of the indeterminacy of scalings under the transformation $(D^0,D^1) \rightarrow (D^0/c,c\cdot D^1)$, most works define convergence using quantities that are \emph{invariant} under this transformation. Let $D^{0(t)},D^{1(t)}$ be the scalings obtained after $t$ iterations of Sinkhorn's algorithm, and let ${A}^{(t)}:=D^{1(t)}AD^{0(t)}$ be the scaled matrix based on these scalings. Since ${A}^{(t)}$ is invariant under the transformation, some earlier works such as \citet{franklin1989scaling,soules1991rate} measure convergence in terms of ${A}^{(t)}$ to the optimally scaled matrix $\hat A= D^{1}AD^{0}$ that has the target row and column sums $p,q$. Most later works focus instead on the convergence of the marginal quantities
      \begin{align*}
          r^{(t)}:={A}^{(t)}\mathbf{1}_m;\quad c^{(t)}:={A}^{(t)T}\mathbf{1}_n
      \end{align*}
to the target row and column sums $p,q$. Since after each iteration in \cref{alg:scaling}, the column constraint is always satisfied: ${A}^{(t)}\mathbf{1}_n=q$, it suffices to focus on the convergence of $r^{(t)}$ to $p$. For example,
      \citet{leger2021gradient} uses the KL divergence $D_{\text{KL}}(r^{(t)}\| p)$, while \citet{altschuler2017near,chakrabarty2021better} use the $\ell^1$ distance $\|r^{(t)}- p\|_1$, which is upper bounded by the KL divergence via Pinsker's inequality. Given the entropy optimization perspective of matrix balancing and Sinkhorn's algorithm in \eqref{eq:relative-entropy-minimization} and \eqref{eq:log-barrier}, it is also possible to measure convergence in terms of the dual optimality gap $g(d^{0(t)},d^{1(t)})-g(d^{0},d^{1})$, which in turn bounds KL divergences via
      \begin{align}
      \label{eq:optimality-gap-property}
g(d^{0(t)},d^{1(t)})-g(d^{0},d^{1}) & = D_{\text{KL}}(D^1AD^0\|A^{(t)})=\sum_{s=t}^\infty D_{\text{KL}}(p\|r^{(s)})+D_{\text{KL}}(q\|c^{(s)}).
\end{align} \citet{luo1992convergence} show that this dual optimality gap converges linearly, but how the convergence rate depends on problem structure has been an open question. Our global linear convergence result in \cref{thm:global-convergence} is the first to characterize this rate.}

\revision{Next, we note the distinction between global and local (particularly asymptotic) convergence results. Global results hold for all iterations $t>0$. For example, \citet{leger2021gradient} shows $D_{\text{KL}}(r^{(t)}\| p)\leq D_{\text{KL}}^\ast /t$ for all $t$, where $D_{\text{KL}}^\ast$ is the optimal value of the relative entropy minimization problem \eqref{eq:relative-entropy-minimization}. On the other hand, local convergence results pertain to the behavior of an algorithm in a neighborhood of the optimal solution, while asymptotic results only hold in the limit as $t\rightarrow \infty$. For example, \citet{knight2008sinkhorn} characterizes the asymptotic rate $\lim_{t\rightarrow \infty}\|r^{(t+1)}-p\|_{\ast}/\|r^{(t)}-p\|_{\ast}$ of linear convergence for Sinkhorn's algorithm, where $\|\cdot\|_\ast$ is an implicitly defined norm. When $A\geq 0$, \citet{knight2008sinkhorn} is the only work on the exact asymptotic convergence rate, for square $A$ and uniform $p,q$. Our  \cref{thm:convergence} is the first asymptotic result for general $A\geq0$ and non-uniform $p,q$, with an explicit norm ($\|\cdot\|_2$). While asymptotic results provide a more precise description of an algorithm's behavior near the optimal solution, global results are useful for obtaining complexity bounds on the number of iterations required to obtain $\varepsilon$ accuracy solutions. In fact, global results are often stated directly as complexity bounds. For example, the result in \citet{altschuler2017near} is that for $A>0$, $t=4(\varepsilon)^{-2}\log(\sum_{ij}A_{ij}/\min_{ij}A_{ij})$ iterations of Sinkhorn's algorithm guarantee $\|r^{(t)}- p\|_1\leq \varepsilon$.} 

\revision{Lastly, we note the distinction between (global) linear and sub-linear convergence results. Linear convergence is often understood as successive improvements of the convergence metric by a constant factor. For example, \citet{franklin1989scaling} show that for $A>0$ and the Hilbert metric $d$, $d(r^{(t+1)}, p) \leq \lambda \cdot d(r^{(t)}, p)$ for all $t>0$ for some $\lambda\in (0,1)$. As a result, $d(r^{(t)}, p) \leq \lambda^t \cdot d(r^{(0)}, p)$ decreases \emph{exponentially} in $t$, so that $d(r^{(t)}, p)\leq \varepsilon$ in $\mathcal{O}(\log(1/\varepsilon))$ iterations. In contrast, in sub-linear results, such as \citet{leger2021gradient}, the convergence metric $D_{\text{KL}}(r^{(t)}\| p)$ only decreases \emph{polynomially} in $t$, requiring $\mathcal{O}(1/\varepsilon)$ iterations to guarantee $D_{\text{KL}}(r^{(t)}\| p) \leq \varepsilon$. While sub-linear complexity bounds have worse (polynomial) dependence on $1/\varepsilon$, they often focus on optimizing the dependence on problem size and dimension. Our main focus in this paper is on understanding the linear convergence behavior of Sinkhorn's algorithm when $A\geq0$, i.e.,  $\mathcal{O}(\log(1/\varepsilon))$ iteration complexity. Nevertheless, we also provide refined complexity bounds in \cref{prop:iteration-complexity} that optimize dependence on problem constants.}

\begin{table}
\caption{Summary of some convergence results on Sinkhorn's algorithm. Throughout, assume that $\|p\|_1=\|q\|_1=1$. Define $r^{(t)}:={A}^{(t)}\mathbf{1}_m$ where ${A}^{(t)}$ is the scaled matrix after $t$ Sinkhorn iterations. In \citet{franklin1989scaling}, $\kappa(A)=\frac{\theta(A)^{1/2}-1}{\theta(A)^{1/2}+1}$,
where $\theta(A)$ is the diameter of $A$ in the Hilbert metric.
The norm in \citet{knight2008sinkhorn} is not explicitly specified, and $\sigma_{2}(\hat{A})$
denotes the second largest singular value of the scaled
doubly stochastic matrix $\hat{A}$. The bound in \citet{altschuler2017near} is
originally stated as the complexity bound that $\|r^{(t)}-p\|_{1}\protect\leq\varepsilon$
in $t=\mathcal{O}(\varepsilon^{-2}\log(\frac{\sum_{ij}A_{ij}}{\min_{ij}A_{ij}}))$
iterations, while the original result in \citet{chakrabarty2021better} is $\|r^{(t)}-p\|_{1}\protect\leq\varepsilon$
in $t=\mathcal{O}(\varepsilon^{-2}\log(\frac{\Delta \max_{ij}A_{ij}}{\min_{ij,A_{ij}>0}A_{ij}}))$
iterations, where $\Delta:=\max_j| i\in [n]:A_{ij}>0|$. The result in \citet{leger2021gradient} applies more generally to couplings
of probability distributions.  
In our asymptotic result, $\lambda_{2}(\tilde{A}^T\tilde{A})$ is the second largest eigenvalue of $\tilde{A} :=\mathcal{D}(1/\sqrt{p})\cdot\hat{A}\cdot\mathcal{D}(1/\sqrt{q})$.
In our global bound, $g^{(t)}=g(d^{0(t)},d^{1(t)})$ while $g^\ast$ is the minimum value of \eqref{eq:log-barrier}. $\lambda_{-2}(\mathcal{L})$ is the second \emph{smallest} eigenvalue of the Laplacian of the bipartite graph defined by $A$, $l=\min \{\max_j (A^T\mathbf{1}_n)_j, \max_i (A\mathbf{1}_m)_i\}$,
$c_B=\exp(-4B)$, and $B$ is a bound on the initial sub-level set, which is finite if and only if \cref{ass:matrix-existence} holds.
}
 \begin{adjustwidth}{-1.5cm}{}
\begin{centering}
\begin{tabular}{c|c|c|c|c}
 & convergence statement & $\lambda$ & $A$ & $p,q$\tabularnewline
\hline 
\citet{franklin1989scaling} & $d_{\text{Hilbert}}(r^{(t)},p)\leq\lambda^t d_{\text{Hilbert}}(r^{(0)},p)$ & $\kappa^{2}(A)$ & $A>0$, rectangular & uniform\tabularnewline
\hline 
\citet{luo1992convergence} & $g^{(t)}-g^\ast\leq\lambda^t (g^{(0)}-g^\ast)$ & \text{unknown} & $A\geq0$, rectangular & general\tabularnewline
\hline 
\citet{knight2008sinkhorn} & $\|r^{(t+1)}-p\|_{\ast}/\|r^{(t)}-p\|_{\ast}\rightarrow\lambda$ & $\sigma_{2}^{2}(\hat{A})$ & $A\geq0$, square  & uniform\tabularnewline
\hline 
%\citet{pukelsheim2009iterative} & $\|r^{(t)}-p\|_{1}\rightarrow0$ & no rate & $A\geq0$, rectangular & general\tabularnewline \hline 
\citet{altschuler2017near} & $\|r^{(t)}-p\|_{1}\leq 2 \sqrt{\frac{\lambda}{t}}$ & $\log(\frac{\sum_{ij}A_{ij}}{\min_{ij}A_{ij}})$ & $A>0$, rectangular & general\tabularnewline
\hline 
\citet{leger2021gradient} & $D_{\text{KL}}(r^{(t)}\| p) \leq\frac{\lambda}{t}$ & $D_{\text{KL}}(\hat{A}\| A)$ & $A\geq0$, continuous & general\tabularnewline
\hline 
current work, asymptotic & $\|\frac{r^{(t+1)}}{\sqrt{p}}-\sqrt{p}\|_{2}/\|\frac{r^{(t)}}{\sqrt{p}}-\sqrt{p}\|_{2}\rightarrow\lambda$ & $\lambda_{2}(\tilde{A}^T\tilde{A})$ & $A\geq0$, rectangular & general\tabularnewline
\hline 
current work, global & $g^{(t)}-g^\ast\leq\lambda^t (g^{(0)}-g^\ast)$ & $1-c_B\lambda_{-2}(\mathcal{L})/l$ & $A\geq0$, rectangular & general\tabularnewline
\end{tabular}
\par\end{centering}
\label{tab:convergence-summary}
\end{adjustwidth}
\end{table}

As discussed before, when $A\geq 0$, only sub-linear convergence results with explicit rates are known \citep{kalantari2008complexity,chakrabarty2021better,leger2021gradient}, while \citet{luo1992convergence} implies global linear convergence with an unknown rate. We now characterize this global linear rate of convergence in terms of the algebraic connectivity of the bipartite graph defined in \eqref{eq:bipartite-adjacency}.

\subsection{Global Linear Convergence}
\label{subsec:global-linear-convergence}
 Our analysis starts with the following change of variables to transform the potential function \eqref{eq:log-barrier}:
\begin{align}
\label{eq:change-of-variables}
    u:=\log d^0,\quad v:=-\log d^1,
\end{align}
resulting in the reparameterized potential function $g(u,v)$ of \eqref{eq:log-barrier}:
\begin{align}
\label{eq:transformed-potential}
  g(u,v)	:=\sum_{ij}A_{ij}e^{-v_{i}+u_{j}}+\sum_{i=1}^{n}p_{i}v_{i}-\sum_{j=1}^{m}q_{j}u_{j}.
\end{align}
Note first that $g(u,v)=g(u+a,v+a)$ for any constant $a\in \mathbb{R}$. We can verify that Sinkhorn's algorithm is equivalent to the alternating minimization algorithm \citep{bertsekas1997nonlinear} for \eqref{eq:transformed-potential}, which alternates between minimizing with respect to $u$ and $v$, holding the other block fixed:
\begin{align}
    u^{(t)}\leftarrow \arg\min_u g(u,v^{(t-1)}), \quad v^{(t)}\leftarrow \arg\min_v g(u^{(t)},v),
\end{align}
or written more explicitly element-wise,
\begin{align}
   \label{eq:alternating-minimization} u_j^{(t)}\leftarrow \log \frac{q_j}{\sum_i A_{ij}e^{-v^{(t-1)}_i}},\quad v_i^{(t)}\leftarrow  \log \frac{p_i}{\sum_j A_{ij}e^{u^{(t)}_j}}.
\end{align}

\revision{A main reason to focus on \eqref{eq:transformed-potential} instead of the log-barrier form \eqref{eq:log-barrier} is that \eqref{eq:transformed-potential} has a Hessian with desirable properties for proving linear convergence. The Hessian of $g(u,v)$ is  
\begin{align}
\nabla^2g(u,v)=\begin{bmatrix}\mathcal{D}(\hat{A}\mathbf{1}_{m}) & -\hat{A}\\
-\hat{A}^{T} & \mathcal{D}(\hat{A}^{T}\mathbf{1}_{n})
\end{bmatrix},
\end{align}
where $\mathcal{D}$ converts a vector to a diagonal matrix, and $\hat{A}=\mathcal{D}(d^1)A\mathcal{D}(d^0)=\mathcal{D}(\exp(-v))A\mathcal{D}(\exp(u))$ is the matrix scaled by $u,v$. Note that the Hessian $\nabla^2g(u,v)$ always has $\mathbf{1}_{m+n}$ in its null space.} On the surface, it may seem that standard linear convergence results for first-order methods, which require strong convexity (or the Polyak--{\L}ojasiewicz condition) of the objective function, do not apply to $g(u,v)$. However, we will show that whenever the matrix balancing problem has finite scaling solutions, $g(u,v)$ is in fact strongly convex when \emph{restricted} to bounded subsets of the subspace 
\begin{align}
\label{eq:orthogonal-subspace}
\mathbf{1}_{m+n}^\perp:= \{u\in \mathbb{R}^m,v\in \mathbb{R}^n:(u,v)^T\mathbf{1}_{m+n}=0\}.
\end{align}
\revision{Moreover, the invariance of $g(u,v)$ and its gradient and Hessian under constant translations of $(u,v)$ by $\mathbf{1}_{m+n}$ guarantees that the strong convexity constant of $g(u,v)$ on $\mathbf{1}_{m+n}^\perp$ in fact quantifies the linear convergence of Sinkhorn's algorithm even if the iterates $u^{(t)},v^{(t)}$ are not in $\mathbf{1}_{m+n}$. Similar types of ``restricted strong convexity'' properties have been studied by for example \citet{agarwal2010fast}. It also shares similarities with the exp-concavity property popular in online learning \citep{hazan2016introduction,orabona2019modern}, which implies the strong convexity of a function in the direction of the gradient, evaluated at \emph{any} point. In contrast, $g(u,v)$ is strongly convex along any direction orthogonal to $\mathbf{1}_{m+n}$, but its gradient evaluated at any $(u,v)$ is not necessarily orthogonal to $\mathbf{1}_{m+n}$. However, the key is that along the trajectory of iterates $(u^{(t)},v^{(t)})$ obtained by running Sinkhorn's algorithm, the gradients of $g$ evaluated at $(u^{(t)},v^{(t)})$ are indeed orthogonal to $\mathbf{1}_{m+n}$, which is sufficient to guarantee the linear convergence of Sinkhorn's algorithm.}

We now introduce the key quantities and definitions used in our result. Let Sinkhorn's algorithm initialize with a $u^{(0)}$, with $v^{(0)}$ given by \eqref{eq:alternating-minimization}. 
Define the constant $B$ as 
\begin{align}
\label{eq:B-bound}
\begin{split}
     B:&= \sup_{(u,v)} \|(u,v)\|_\infty\\  \text{subject to } (u,v)^{T}&\mathbf{1}_{m+n}=0,\\
   g(u,v)&\leq g(u^{(0)},v^{(0)})
\end{split}
\end{align}
In other words, $B$ is the \emph{diameter} of the initial normalized sub-level set. \revision{We will show that $B$ is finite and that it bounds normalized Sinkhorn iterates $\|(u^{(t)},v^{(t)})\|_\infty$, since under \cref{ass:matrix-existence} the function $g(u,v)$ is \emph{coercive} on the subspace $\mathbf{1}^\perp_{m+n}$. Coercivity is an important property, and we define it below, following \citet{bertsekas2016nonlinear}. 
\begin{quote}\textbf{Definition 1. }(Coercivity) 
A function $f(x): \mathbb{R}^d \rightarrow \mathbb{R}$ is coercive on a subspace $V\subseteq \mathbb{R}^d$ if
\begin{align}
\label{eq:coercivity}
    f(x^{(t)}) \rightarrow +\infty \quad \text{whenever } x^{(t)} \in V \quad \text{and}
    \quad \|x^{(t)}\|_\infty \rightarrow +\infty.
\end{align}
\end{quote}
}
% Define the normalized optimal solution pair 
% \begin{align}
% \label{eq:normalized-optimum}
% (u^\ast,v^\ast):=\arg \min_{(u,v)\in \mathbf{1}_{m+n}^\perp} g(u,v),
% \end{align}
% and 
Next, define the \emph{Laplcian} matrix $\mathcal{L}:=\mathcal{L}(A)$ of the bipartite graph $G_b$ (see \eqref{eq:bipartite-adjacency}) as 
\begin{align}
\label{eq:graph-laplacian}
\mathcal{L}:&=\begin{bmatrix}\mathcal{D}({A}\mathbf{1}_{m}) & -{A}\\
-{A}^{T} & \mathcal{D}({A}^{T}\mathbf{1}_{n})
\end{bmatrix}, 
\end{align}
and refer to the second smallest eigenvalue $\lambda_{-2}(\mathcal{L})$ as the Fiedler eigenvalue. $\lambda_{-2}(\mathcal{L})>0$ if and only if \cref{ass:matrix-uniqueness} holds, and it quantifies ``connectivity'' of the data structure \citep{spielman2007spectral}. Although an important quantity in the choice modeling literature, algebraic connectivity has not been previously used in the analysis of Sinkhorn's algorithm. For details on the graph Laplacian and the Fiedler eigenvalue, see \cref{subsec:graph-laplacian}.  Finally, define the smoothness parameters
\begin{align}
\label{eq:smoothness}
    l_0:= \max_j (A^T\mathbf{1}_n)_j, \quad l_1:= \max_i (A\mathbf{1}_m)_i,
\end{align}
which are used to quantify the smoothness of $g(u,v)$.

We can now state one of our main contributions to the study of Sinkhorn's algorithm. 
\begin{theorem}[\textbf{Global Linear Convergence}]
\label{thm:global-convergence}
Suppose \cref{ass:matrix-existence} and
\cref{ass:matrix-uniqueness} hold. \revision{Let $\mathcal{L}$ be the bipartite graph Laplacian defined in \eqref{eq:graph-laplacian} and $\lambda_{-2}(\mathcal{L})$ its second smallest eigenvalue. Let 
$l_0, l_1$ be the smoothness parameters defined in \eqref{eq:smoothness}. Let $(u^{(t)},v^{(t)})$ be Sinkhorn iterates at iteration $t$ defined in \eqref{eq:alternating-minimization} and $B$ the bound on $\|(u^{(t)},v^{(t)})\|_\infty$ defined in \eqref{eq:B-bound}.
Define $g^\ast:=\inf_{u,v}g(u,v)$. For all $t>0$, the optimality gap of the dual objective $g(u,v)$ defined in \eqref{eq:transformed-potential} satisfies
\begin{align}
\label{eq:potential convergence}
g(u^{(t+1)},v^{(t+1)}) - g^\ast \leq (1-e^{-4B} \frac{\lambda_{-2}(\mathcal{L})} {\min\{l_0,l_1\}})\left( g(u^{(t)},v^{(t)})-g^\ast \right).
    \end{align} 
    The ratio $\min\{l_{0},l_{1}\}/\lambda_{-2}(\mathcal{L})$ can be interpreted as a condition number of $\mathcal{L}$.}
\end{theorem}
\revision{The linear convergence rate of Sinkhorn's algorithm is therefore quantified by $\lambda_{-2}(\mathcal{L})/\min\{l_{0},l_{1}\}$, which is invariant under rescalings of $A\rightarrow c\cdot A$. Although the corresponding bipartite graph $G_b$ with biadjacency matrix $A$ is a natural object to consider in the study matrix balancing problems, to our knowledge, \cref{thm:global-convergence} is the first to highlight the precise role of its spectral property, described by $\lambda_{-2}(\mathcal{L})$, in the linear convergence of Sinkhorn's algorithm. It fills the gap left by \citet{luo1992convergence} who establish linear convergence with an implicit rate, and allows us to compute its dependence on problem parameters by applying established bounds on $\lambda_{-2}(\mathcal{L})$ from spectral graph theory \citep{spielman2007spectral}.}
\begin{remark}[\textbf{Importance of Assumptions for Linear Convergence}]
The importance of Assumptions \ref{ass:matrix-existence} and \ref{ass:matrix-uniqueness} are clearly reflected in the bound \eqref{eq:potential convergence}.
First, note that the Fiedler eigenvalue $\lambda_{-2}(\mathcal{L})>0$ if and only if \cref{ass:matrix-uniqueness} holds (see \cref{subsec:graph-laplacian}). On the other hand, \cref{ass:matrix-existence} guarantees the \emph{coercivity} of $g$ on $\mathbf{1}_{m+n}^\perp$ (see \eqref{eq:coercivity}). This property ensures that $B$ defined in \eqref{eq:B-bound} satisfies $B<\infty$, and consequently, that normalized iterates stay bounded by $B$. That \cref{ass:matrix-existence} guarantees  $g(u,v)$ is coercive should be compared to the observation by \citet{hunter2004mm} that \cref{ass:strong-connected} guarantees the upper compactness (a closely related concept) of the log-likelihood function \eqref{eq:log-likelihood}. When \cref{ass:matrix-existence} fails, $B$ may become infinite and $\|(u^{(t)},v^{(t)})\|_\infty \rightarrow \infty$.
\end{remark} 

\begin{remark}[\textbf{Self-normalizing Property of Sinkhorn}]
   The ability of Sinkhorn's algorithm to exploit the (subspace) strong convexity of $g(u,v)$ on $\mathbf{1}_{m+n}^\perp$ to achieve linear convergence relies critically on the invariance of the scaled matrix $\hat{A}=D^1AD^0$ and $g(u,v)$ under the transformation $(D^0,D^1) \rightarrow (D^0/c,c\cdot D^1)$. This is an intrinsic feature of the matrix balancing problem that has been well-known but not fully exploited in the convergence analysis so far. %Recall that $u=\log d^0$ and $v=-\log d^1$, where $d^0,d^1$ are the diagonals of the scaling $(D^0,D^1)$. 
   It guarantees that the translation $(u,v)\rightarrow(u-\log c,v-\log c)$ does not alter $g(u,v)$ and its derivatives in \eqref{eq:transformed-potential}. We can therefore impose the \emph{auxiliary} normalization $(u,v)^T\mathbf{1}_{m+n}=0$, or equivalently $\prod_j d^0_j = \prod_i d^1_i$, which is easily achieved by requiring that after every update in \cref{alg:scaling}, a normalization $(d^0/c,c d^1)$ is performed using the normalizing constant 
\begin{align}
\label{eq:normalization}
   c=\sqrt{\prod_j d^0_j /\prod_i d^1_i}.
\end{align}
See \cref{alg:scaling-normalized} for the normalized Sinkhorn's algorithm, which given \eqref{eq:normalization} results in a virtual sequence of $u^{(t)},v^{(t)}$ satisfying $(u^{(t)},v^{(t)})^T\mathbf{1}_{m+n}=0$. Moreover, the values of $g(u,v)$ on this virtual sequence is identical to those on the standard Sinkhorn iterates. As a result, the convergence result \eqref{eq:potential convergence} applies to the standard Sinkhorn's algorithm without normalization (or with any other normalization), due to the invariance of $g(u,v)$. Normalization of Sinkhorn's algorithm is also considered in the analyses in \citet{carlier2023sista}, although they use the asymmetric condition $u_0=0$, which does not guarantee that normalized Sinkhorn iterates stay in $\mathbf{1}_{m+n}^\perp$.
\end{remark} 

 With this auxiliary normalization procedure, the proof of \cref{thm:global-convergence} then relies on the observation that the Hessian of $g(u,v)$ is precisely the graph {Laplacian} $\mathcal{L}(u,v)$ of the bipartite graph with biadjacency matrix $\hat A = \mathcal{D}(\exp(-v))A\mathcal{D}(\exp(u))$. As $(u,v)$ are {bounded} on normalized Sinkhorn iterates thanks to the coercivity of $g$, the Fiedler eigenvalue of $\mathcal{L}=\mathcal{L}(0,0)$ quantifies the strong convexity on $\mathbf{1}_{m+n}^\perp$. Linear convergence then follows from results on block coordinate descent and alternating minimization methods for strongly convex and smooth functions \citep{beck2013convergence}. Typically, the leading eigenvalue of the Hessian quantifies the smoothness, which is bounded by $2\max\{l_0,l_1\}$ for $\mathcal{L}$. For alternating minimization methods, the better smoothness constant $\min\{l_0,l_1\}$ is available. {Thus the quantity $\min\{l_{0},l_{1}\}/\lambda_{-2}(\mathcal{L})$ in \eqref{eq:potential convergence} can be interpreted as a type of  ``condition number'' of the graph Laplacian $\mathcal{L}$.} When $A$ is positive (not just non-negative), then the strong existence and uniqueness conditions are trivially satisfied, and our results continue to hold with the rate quantified by $\min\{l_{0},l_{1}\}/\lambda_{-2}(\mathcal{L})$. In this case, both $\min \{l_0,l_1\}$ and  $\lambda_{-2}(\mathcal{L})$ are $\Theta(n)$ where $n$ is  problem dimension, so $\min\{l_{0},l_{1}\}/\lambda_{-2}(\mathcal{L})$ does not increase with problem dimension.
\revision{
\begin{remark}[\textbf{Significance of \cref{thm:global-convergence}}]
    A main innovation in our paper is in introducing the concept of algebraic connectivity when quantifying the global convergence of Sinkhorn's algorithm for non-negative matrices. In this respect, the significance of \cref{thm:global-convergence} is more conceptual than technical, since once we identify the right quantity (algebraic connectivity) and utilize the self-normalizing property, the convergence result can be obtained using standard matrix analysis
and applying the theory of \citet{beck2013convergence} for block coordinate descent algorithms. Nevertheless, we feel the role of algebraic connectivity in the study of matrix balancing problems holds general significance and likely will lead to more results in related areas. See for example \citet{chang2024inferring} which highlights its importance for the statistical efficiency of a network traffic model based on matrix balancing.
\end{remark}
}

\revision{Although
\cref{thm:global-convergence} implies an $\mathcal{O}(\log(1/\varepsilon))$ iteration complexity, the complexity bound's dependence on problem parameters can be further improved. In particular, the constant $B$ which bounds $\|(u^{(t)},v^{(t)})\|_\infty$ can be hard to compute for some problems. We next establish an iteration complexity bound that does not depend exponentially on the implicit constant $B$ defined in \eqref{eq:B-bound}.
\begin{proposition}[\textbf{Iteration Complexity}]
\label{prop:iteration-complexity} Under \cref{ass:matrix-existence} and \cref{ass:matrix-uniqueness}, let $d^0_\ast,d^1_\ast$ be a pair of optimal scalings. Define $\|v\|_{-\infty}:=\min_i|v_i|$ and let $C=:\max \Big\{ \tfrac{\| d^0_{\ast} \|_{\infty}}{\|
d^0_{\ast} \|_{- \infty}}, \tfrac{1}{\| d^0_{\ast} \|_{- \infty} \|
d^1_{\ast} \|_{- \infty}}, \| d^0_{\ast} \|_{\infty} \| d^1_{\ast}
\|_{\infty} \Big\}$. Suppose Sinkhorn's algorithm initializes with $u^{(0)}=\mathbf{1}_m$. Then $\|(e^{u^{(t)}},e^{v^{(t)}})\|_\infty \leq C$ for all $t>0$.
Moreover, for any $\varepsilon \leq \frac{1}{2} \min \{ \| p \|_{- \infty}, \| q \|_{- \infty} \}$, after
      \begin{align}
      \label{eq:complexity-bound}
     \mathcal{O} \left (C^2\cdot \frac{\min\{\|p\|_\infty, \|q\|_\infty\}}{\lambda_{-2}(\mathcal{L})} \cdot (\log (1/\varepsilon) +\log\log C) \right )
      \end{align}
       iterations of Sinkhorn's algorithm, the optimality gap and the $\ell^1$ distance $\|r^{(t)}-p\|_1\leq \varepsilon$. 
\end{proposition} 
}

\revision{
\cref{prop:iteration-complexity} relies on two technical innovations. First, we bound $e^B \leq C$ where $C$ is explicitly constructed from any optimal solution pair and is invariant under rescalings. Second, we improve the dependence from $e^{4B}$ to $e^{2B}$, using the target marginals $p,q$ to quantify the smoothness of $g(u,v)$ instead. Our message here is that the convergence behavior of Sinkhorn's algorithm has two phases. Initially, we can apply a sub-linear complexity bound with $\mathcal{O}(1)$ iterations to  obtain Sinkhorn iterates sufficiently close to the optimal solution. Afterwards, the convergence can better be captured by a linear convergence with rate depending on the \emph{optimal} solution and target marginals $p,q$. The dependence of $C$ on problem dimension is problem specific. In the worst case, it can be exponential \citep{kalantari1993rate}. In Appendix \ref{sec:complexity-bound-plots}, we plot \eqref{eq:complexity-bound} as a function of problem dimension on randomly generated data, and find the dependence is quadratic. In contrast, sub-linear bounds, such as \citet{altschuler2017near,chakrabarty2021better}, have logarithmic dependence on problem dimension. It remains an interesting question to improve the dependence on problem dimension in \eqref{eq:complexity-bound}, and to study trade-offs with the dependence on $\varepsilon$.
\subsection{Strong vs. Weak Convergence of Sinkhorn's Algorithm}
We now discuss the two different convergence regimes of Sinkhorn's algorithm when $A\geq 0$. As mentioned in \cref{subsec:matrix-balancing,subsec:convergence-background}, when $A\geq 0$, the canonical matrix balancing problem with target marginals $p,q$ has a finite positive solution pair $D^0,D^1$ if and only if \cref{ass:matrix-existence} holds (which trivially holds when $A>0$). In this case, Sinkhorn's algorithm converges to $D^1AD^0$, which also solves the KL minimization problem \eqref{eq:relative-entropy-minimization}. We call this case the \emph{strong convergence} of Sinkhorn.}

\revision{However, even if \cref{ass:matrix-existence} fails and no positive finite scalings $D^0,D^1$ exist that solve the matrix balancing problem, the sequence of scaled matrices ${A}^{(t)}=D^{1(t)}AD^{0(t)}$ based on Sinkhorn's algorithm can still converge entry-wise to the solution of \eqref{eq:relative-entropy-minimization}, whenever it has a finite solution. This apparent discrepancy is explained by the fact that the solution of problem \eqref{eq:relative-entropy-minimization} requires a weaker condition than \cref{ass:matrix-existence} for the matrix balancing problem. It can be stated in the following equivalent forms.}
\begin{assumption}[\textbf{Weak Existence}]
    \label{ass:matrix-weak-existence}
    
    \textbf{(a)} There exists a non-negative matrix $A'\in \mathbb{R}_+^{n\times m}$ that inherits all zeros of $A$ and has row and column sums $p$ and $q$. Or, equivalently,
    
\textbf{(b)} For every pair of sets of indices $N \subsetneq [n]$ and $M \subsetneq [m]$ such that $A_{ij}=0$ for $i\notin N$ and $j\in M$, $\sum_{i\in N}p_i \geq \sum_{j\in M}q_j$.
\end{assumption}

The equivalence of the two conditions above follows from Theorem 4 in \citet{pukelsheim2009iterative}, which also shows that they are the minimal requirements for the convergence of Sinkhorn's algorithm. \cref{ass:matrix-weak-existence}(a) precisely guarantees 
that the constrained KL minimization problem \eqref{eq:relative-entropy-minimization} is feasible and bounded. It relaxes \cref{ass:matrix-existence}(a) by allowing additional zeros in the matrix $A'$. Similarly, \cref{ass:matrix-weak-existence}(b) relaxes \cref{ass:matrix-existence}(b) by allowing equality between $\sum_{i\in N}p_i$ and $\sum_{j\in M}q_j$ even when $M,N$ do not correspond to a block-diagonal structure. 

The distinction between \cref{ass:matrix-existence} and \cref{ass:matrix-weak-existence} is important for the matrix balancing problem and Sinkhorn's algorithm. \cref{ass:matrix-existence} guarantees the solutions of \eqref{eq:matrix-balancing} and \eqref{eq:relative-entropy-minimization} coincide, and have exactly the same zero pattern as $A$. 
If \cref{ass:matrix-weak-existence} holds but \cref{ass:matrix-existence} fails, then the solution $\hat A$ of \eqref{eq:relative-entropy-minimization}  has \emph{additional} zeros relative to $A$, and no direct (finite and positive) scaling $(D^{0},D^{1})$ exists such that $\hat A=D^1AD^0$. However, the sequence of scaled matrices $\hat A^{(t)}$  still converges to $\hat A$. We call this case the \emph{weak convergence} of Sinkhorn. In this case, the matrix balancing problem is said to have a \emph{limit} scaling, where some entries of $D^{0},D^{1}$ in Sinkhorn iterations approach 0 or $\infty$, resulting in additional zeros in $\hat A$. Below we give an example adapted from \citet{pukelsheim2009iterative}, where $p,q=(3,3)$ and the scaled matrices $\hat A^{(t)}$ converge but no direct scaling exists:
\begin{align}
\label{eq:counter-example}
A^{(t)}={D^{1}}^{(t)}A{D^{0}}^{(t)}={D^{1}}^{(t)}\begin{bmatrix}3&1\\
0 & 2
\end{bmatrix}{D^{0}}^{(t)} = \begin{bmatrix}1&0\\
0 & \frac{3t}{2}
\end{bmatrix}\begin{bmatrix}3&1\\
0 & 2
\end{bmatrix}\begin{bmatrix}1&0\\
0 & \frac{1}{t}
\end{bmatrix}	\rightarrow \begin{bmatrix}3&0\\
0 & 3
\end{bmatrix}.
\end{align}

Under \cref{ass:matrix-weak-existence}, \citet{leger2021gradient}proves the sub-linear convergence of Sinkhorn's algorithm, while it is known since at least \citet{soules1991rate,luo1992convergence} that the convergence is linear under \cref{ass:matrix-existence}. It is therefore important to clarify the convergence behaviors of Sinkhorn's algorithm in the two settings.
We next show that if \cref{ass:matrix-weak-existence} holds but \cref{ass:matrix-existence} fails, then there exists an entry of $A^{(t)}$ that converges at a lower bound rate $\Omega(1/t)$, i.e., sub-linear.
%Convergence can be sharpened to linear if and only if the weak existence condition (\cref{ass:matrix-weak-existence}) is replaced by the strong existence condition (\cref{ass:matrix-existence}). 
Together with existing and new results in this paper, \cref{thm:lower-bound} fully characterizes the following convergence behavior of Sinkhorn's algorithm: whenever a direct scaling exists for the matrix balancing problem, Sinkhorn's algorithm converges linearly. If only a limit scaling exists, then convergence deteriorates to sub-linear. \revision{This generalizes the observations made by \citet{sinkhorn1967concerning,achilles1993implications} for square matrices and uniform marginals.} 
\begin{proposition}[\revision{Linear vs. Sub-linear Convergence of Sinkhorn}]
\label{thm:lower-bound}
       For general non-negative matrices, Sinkhorn's algorithm converges linearly
iff $(A,p,q)$ satisfy \cref{ass:matrix-existence} and \cref{ass:matrix-uniqueness}. The convergence is sub-linear if and only if the weak existence condition \cref{ass:matrix-weak-existence} holds but \cref{ass:matrix-existence} fails.
\end{proposition}

The regime of sub-linear convergence also has an interpretation in the choice modeling framework. The weak existence condition \cref{ass:matrix-weak-existence}, when applied to $(A,p,q)$ constructed from a choice dataset, allows the case where some subset $S$ of items is always preferred over $S^C$, which implies, as observed already by the early work of \citet{ford1957solution}, that the log-likelihood function \eqref{eq:log-likelihood} is only maximized at the \emph{boundary} of the probability simplex, by shrinking  $s_j$ for $j\in S^C$ towards 0, i.e., $D^0_j \rightarrow 0$. Incidentally, \citet{bacharach1965estimating} also refers to the corresponding regime in matrix balancing as ``boundary solutions''. 

\subsection{Sharp Asymptotic Rate}
\label{subsec:sharp-rate}
\revision{Having established the global convergence and iteration complexity of Sinkhorn's algorithm when finite scalings exist, we now turn to the problem of characterizing the sharp, i.e., best possible asymptotic linear convergence rate as $t\rightarrow \infty$, for general non-negative $A$ and non-uniform marginals $(p,q)$. \citet{knight2008sinkhorn} computed this rate for uniform $(p,q)$ under an implicit metric.} Our analysis is distinct from theirs and relies on an intrinsic orthogonality structure of Sinkhorn's algorithm, which is also different from the auxiliary normalization in our global linear convergence analysis. Note that unlike the global rate, which depends on the initial problem data $A$ and $(u^{(0)},v^{(0)})$, the asymptotic rate now depends on the optimal solution $\hat A=D^1AD^0$, as expected.  
\begin{theorem}[\textbf{Sharp Asymptotic Rate}]
\label{thm:convergence}
Suppose $(A,p,q)$ satisfy \cref{ass:matrix-existence} and \cref{ass:matrix-uniqueness}. \revision{Let $\hat{A}$
be the unique scaled matrix with target marginals $p,q$ defined in \eqref{eq:scaled-matrix}. Then marginals $r^{(t)}=A^{(t)}\mathbf{1}$, where $A^{(t)}$ is the scaled matrix after $t$ iterations of Sinkhorn's algorithm, satisfy}
\begin{align}
\lim_{t\rightarrow \infty} \frac{\|r^{(t+1)}/\sqrt{p}-\sqrt{p}\|_{2} }{\|r^{(t)}/\sqrt{p}-\sqrt{p}\|_{2}} = \lambda_\infty,
\end{align}
where the asymptotic linear rate of convergence $\lambda_\infty$ is given by
\begin{align*}
\lambda_\infty & :=\lambda_{2}(\tilde{A}\tilde{A}^{T})=\lambda_{2}(\tilde{A}^{T}\tilde{A})\\
\tilde{A} & :=\mathcal{D}(1/\sqrt{p})\cdot\hat{A}\cdot\mathcal{D}(1/\sqrt{q}),
\end{align*}
and $\lambda_{2}(\cdot)$ denotes the second largest eigenvalue. 
\end{theorem}

In the special case of $m=n$ and $p=q=\mathbf{1}$, the asymptotic rate in \cref{thm:convergence} reduces to that in \citet{knight2008sinkhorn}. Note, however, that the convergence metric is different: we use the $\ell^2$ norm $\|r^{(t)}/\sqrt{p}-\sqrt{p}\|_2$ while \citet{knight2008sinkhorn} uses $\|r^{(t)}-p\|_\ast$ with an {implicit} norm $\|\cdot\|_\ast$ on $\mathbb{R}^n$. \minor{Moreover, one cannot directly extend results for square matrices such as theirs to non-square matrices by padding them with zeros, as doing so results in target marginals that are not strictly positive. See, however, \citet{knight2008sinkhorn} for a symmetrization proposal.} %Our analysis exploits the orthogonality structure of Sinkhorn's algorithm and more explicitly reveals the dependence of the convergence rate on the spectral structure of the data. 

The proof of \cref{thm:convergence} relies on a sequence of novel data-dependent mappings associated with Sinkhorn's algorithm. Intuitively, the dependence of the asymptotic linear rate on the second largest eigenvalue of $\tilde{A}^T\tilde{A}$ (and $\tilde{A}\tilde{A}^T$) is due to the fact that near the fixed point $\sqrt{p}$ of the mapping associated with Sinkhorn iterations, $\tilde{A}\tilde{A}^T$ (which is the Jacobian at $\sqrt{p}$) approximates the first order change in $r^{(t)}/\sqrt{p}$. Normally, the \emph{leading} eigenvalue quantifies this change. The unique leading eigenvalue of $\tilde{A}\tilde{A}^T$ is equal to 1 with eigenvector $\sqrt{p}$, which does not imply contraction. Fortunately, using the quantity  $r^{(t)}/\sqrt{p}$
allows us to exploit the following orthogonality structure:
\begin{align*}
(r^{(t)}/\sqrt{p}-\sqrt{p})^{T}\sqrt{p} & =\sum_{i}(r_{i}^{(t)}-p_{i})=0
\end{align*}
by virtue of Sinkhorn's algorithm preserving the quantities $r^{(t)T}\mathbf{1}_{n}$
for all $t$. Thus, the residual $r^{(t)}/\sqrt{p}-\sqrt{p}$ is always \emph{orthogonal} to
$\sqrt{p}$, which is both the leading eigenvector and the fixed point of the iteration. The convergence is then controlled by the \emph{second} largest eigenvalue of $\tilde{A}\tilde{A}^T$. This proof approach echoes that of the global linear convergence result in \cref{thm:global-convergence}, where we also exploit an orthogonality condition to obtain a meaningful bound. In \cref{thm:global-convergence} the bound depends on the second smallest eigenvalue of a Hessian matrix, while in \cref{thm:convergence} the bound depends on the second largest eigenvalue of a Jacobian matrix. 

{
Lastly, we note that the asymptotic rate $\lambda_\infty$ is itself a Fiedler eigenvalue, associated with the Laplacian that is the Schur complement of the \emph{scaled} graph Laplacian 
\begin{align*} \begin{bmatrix}\mathcal{D}(1/\sqrt{p}) & 0\\
0 & \mathcal{D}(1/\sqrt{q})
\end{bmatrix}\begin{bmatrix}\mathcal{D}(\hat A\mathbf{1}_{m}) & -\hat{A}\\
-\hat{A}^{T} & \mathcal{D}(\hat{A}^{T}\mathbf{1}_{n})
\end{bmatrix}\begin{bmatrix}\mathcal{D}(1/\sqrt{p}) & 0\\
0 & \mathcal{D}(1/\sqrt{q})
\end{bmatrix}.
\end{align*}
}

%% file: conclusion.tex
\section{Conclusion}
In this paper, we develop extensive connections between matrix balancing and choice modeling. We show that the maximum likelihood estimation of choice models based on the Luce axioms of choice is an instance of the canonical matrix balancing problem. Moreover, many algorithms in choice modeling can be viewed as special cases or analogs of Sinkhorn's algorithm for matrix balancing. These connections can potentially benefit multiple disciplines. For choice modeling, they open the door to tools and insights from well-studied topics in optimization and numerical linear algebra. For matrix balancing, the connections enable us to resolve some interesting open problems on the linear convergence of Sinkhorn's algorithm for non-negative matrices. We establish the first quantitative global linear convergence result for Sinkhorn's algorithm applied to general non-negative matrices. Our analysis reveals the importance of algebraic connectivity for matrix balancing. We also provide the first characterization of the exact asymptotic linear rate of convergence for general non-negative matrix and non-uniform target marginals. Lastly, we clarify the linear and sub-linear convergence behaviors of Sinkhorn's algorithm under the strong and weak existence assumptions for matrix balancing.
%Moreover, we propose regularization methods for Sinkhorn's algorithm inspired by works from choice modeling, in order to address existence and convergence issues for matrix balancing. 
Overall, we believe that the connections established in this paper are useful for researchers from different domains and can lead to further interesting results.

% Our results in this section are relevant in several respects. First, we clarify the gap between the $\mathcal O(1/t)$ and $\mathcal O(\lambda^t)$ convergence of Sinkhorn's algorithm: the slowdown happens if and only if Sinkhorn's algorithm converges but the canonical matrix balancing problem does not have a \emph{finite} scaling $(D^0,D^1)$. This slowdown has been observed in the literature but not systematically studied. Second, 
%   The importance of algebraic connectivity for Sinkhorn's algorithm becomes less surprising once we connect it to messages passing algorithms and the distributed optimization literature in \cref{sec:connections}, where it is well-known that the spectral gap of the \emph{gossip} matrix, which defines the decentralized communication network, governs the rates of convergence.

\ACKNOWLEDGMENT{This work is supported in part by NSF CAREER Award \#2143176 and European Research Council consolidator grant (ERC-CoG) No. 866274. We are very grateful for insightful comments and suggestions from Kurt Anstreicher, Dimitris Bertsimas, St\'ephane Bonhomme, Agostino Capponi, Serina Chang, Patrick Ding, Yanqin Fan, Robert Freund, Han Hong, Yuchen Hu, Guido Imbens, S{\"u}leyman Kerimov, Yongchan Kwon, Frederic Koehler, Flavien L\'eger, Gary Qian, Pavel Shibayev, Ruoxuan Xiong, and Yinyu Ye. We also thank the area editor, the associate editor, and the three anonymous referees for helping to improve the manuscript significantly.}

%% file: connection.tex
%\section{Applications to Existing Algorithms and Frameworks for Choice Modeling}
\section{Graph Laplacians and Algebraic Connectivity}
\label{subsec:graph-laplacian}
In this section, we introduce the quantities central to our global linear convergence analysis, especially the \emph{graph Laplacian} matrices associated with the graphs defined by a non-negative matrix $A$ and the Fielder eigenvalues.

Given a non-negative matrix $A\in \mathbb{R}_+^{n\times m}$, we define the associated (weighted) bipartite graph $G_b$ on $V\cup U$ by the adjacency matrix $A^b \in \mathbb{R}^{(m+n)\times (m+n)}$ defined as
\begin{align*}
    A^b := \begin{bmatrix}\mathbf{0} & {A}\\
{A}^{T} & \mathbf{0}
\end{bmatrix}.	
\end{align*}
The rows of $A$ correspond to vertices in $V$
with $|V|=n$, while the columns of $A$ correspond to vertices in
$U$ with $|U|=m$, and $V\cap U=\emptyset$. The matrix $A$ here is sometimes called the \emph{biadjacency matrix} of the bipartite graph.

In the context of choice modeling, $A$ also defines
an \emph{undirected} ``comparison'' graph $G_c$ on $m$ items. This is most easily understood when $A$ is binary and we can associate it with a choice dataset, with its rows encoding choice sets (see \cref{prop:mle-scaling}), but the definition below is more general. Define the adjacency matrix $A^c \in \mathbb{R}^{m\times m}$ by
\begin{align*}
{A^c}_{jj'} & =\begin{cases}
0 & j=j'\\
(A^{T}A)_{jj'} & j\neq j',
\end{cases}
\end{align*}
If $A$ is a binary participation matrix associated with a choice dataset, then there is a (weighted) edge in $G_{c}$ between items $j$ and $j'$
if and only the two appear in some choice set together, with the edge
weight equal to the number of times of their co-occurrence. This undirected comparison graph $G_c$ is not the same as the directed comparison graph in \cref{ass:strong-connected}, since it does not encode the \emph{choice} of each observation. However, it is also an important object in choice modeling. For example, the uniqueness condition in \cref{ass:weak-connected} for choice maximum likelihood estimation has a concise graph-theoretic interpretation as it is a requirement that $G_c$ be connected.

For a (generic) undirected graph $G$ with adjacency matrix $M$, the graph Laplacian matrix (or simply the Laplacian) is defined as $L(M):=\mathcal{D}(M\mathbf{1})-M$, where recall $\mathcal{D}$ is the diagonalization of a vector. The graph Laplacian $L(M)$ is always positive semidefinite as a result of the Gershgorin circle theorem, since $L(M)$ is diagonally dominant  with positive diagonal and non-positive off-diagonals. Moreover, the Laplacian always has $\mathbf{1}$ in its null space, and is spanned by it if and only if and only if $G$ is connected \citep{spielman2007spectral}.

For the graphs $G_b,G_c$, their Laplacians are  given respectively by 
\begin{align}
    \label{eq:Laplacians}
\mathcal{L}:&=L(A^b)=\begin{bmatrix}\mathcal{D}({A}\mathbf{1}_{m}) & -{A}\\
-{A}^{T} & \mathcal{D}({A}^{T}\mathbf{1}_{n})
\end{bmatrix} \\L:&=L(A^c)=\mathcal{D}(A^TA\mathbf{1}_m)-A^TA,
\end{align}
where we can
verify that for the comparison graph $G_c$, its Laplacian $L$ satisfies
\begin{align*}
L=\mathcal{D}(A^c\mathbf{1}_{m})-A^c & =\mathcal{D}(A^{T}A\mathbf{1}_{m})-A^{T}A.
\end{align*}
The graph Laplacian $\mathcal{L}$ based on $A^b$ and $L$ based on $A^c$ are closely connected through the identity
\begin{align*} (A^{b})^{2}&=\begin{bmatrix}AA^{T} & 0\\
0 & A^{T}A
\end{bmatrix},
\end{align*}
which implies that $L$ is the lower right block of the graph Laplacian $\mathcal{D}((A^{b})^{2}\mathbf{1}_{m+n})-(A^{b})^{2}$. Moreover, $L$ plays a central role in works on the statistical and computational efficiency in choice modeling \citep{shah2015estimation,seshadri2020learning,vojnovic2020convergence}.

An important concept in spectral graph theory is the \emph{algebraic connectivity} of a graph, quantified by the second smallest eigenvalue $\lambda_{-2}$ of the graph Laplacian matrix, also called the Fiedler eigenvalue \citep{fiedler1973algebraic,spielman2007spectral}.
% The Fiedler eigenvalue $\lambda_{-2}$ features prominently in Cheeger's inequality,
% \begin{align}
% \label{eq:Cheeger}
%   \frac{h^2(G)}{2\max_j(M\mathbf{1})_j}   \leq \lambda_{-2}(\mathcal{D}(M\mathbf{1})-M) \leq  2h(G),
% \end{align}
% where $h(G)\geq0$ is the Cheeger constant, positive if and only if $G$ is connected.
Intuitively, Fiedler eigenvalue quantifies how well-connected a graph is in terms of how many edges need to be removed for the graph to become disconnected. It is well-known that the multiplicity of the smallest eigenvalue of the graph Laplacian, which is 0, describes the number of connected components of a graph. The uniqueness condition for matrix balancing in \cref{ass:matrix-uniqueness} therefore guarantees that the Fiedler eigenvalue of $G_b$ is positive: $\lambda_{-2}(\mathcal{L}) >0$. This property is important for our results, since $\lambda_{-2}(\mathcal{L})$ quantifies the \emph{strong} convexity of the potential function and hence the linear convergence rate of Sinkhorn's algorithm.

\section{Further Connections to Choice Modeling and Optimization}
\label{sec:connections}
In this section, we demonstrate that our matrix balancing formulation \eqref{eq:equation-system} of the maximum likelihood problem \eqref{eq:log-likelihood} provides a unifying perspective on many existing works on choice modeling, and establishes interesting connections to distributed optimization as well. Throughout, Sinkhorn's algorithm will serve as the connecting thread. In particular, it reduces algebraically to the algorithms in \citet{zermelo1929berechnung,dykstra1956note,ford1957solution,hunter2004mm,maystre2017choicerank} in their respective choice model settings. This motivates us to provide an interpretation of Sinkhorn's algorithm as a ``minorization-maximization'' (MM) algorithm \citep{lange2000optimization}. Moreover, Sinkhorn's Algorithm is also related to the ASR algorithm of \citet{agarwal2018accelerated} for choice modeling, as they can both be viewed as message passing algorithms in distributed optimization \citep{balakrishnan2004polynomial}. Last but not least, we establish a connection between Sinkhorn's algorithm and the well-known BLP algorithm of \citet{berry1995automobile}, widely used in economics to estimate consumer preferences from data on market shares. 

\subsection{Pairwise Comparisons}
The same algorithmic idea in many works on pairwise comparisons appeared as early as \citet{zermelo1929berechnung}. For example,
\citet{dykstra1956note} gives the following update formula: 
\begin{align}
\label{eq:pairwise}
s_{j}^{(t+1)} & =W_{j}/\sum_{j\neq k}\frac{C_{jk}}{s_{j}^{(t)}+s_{k}^{(t)}},
% s_{j}^{(t+1)} & =W_{j}\cdot\left[\sum_{j\neq k}\frac{C_{jk}}{s_{j}^{(t)}+s_{k}^{(t)}}\right]^{-1},
\end{align}
 where again $W_{j}=|\{i\mid (j,S_i)\}|$ is the number of times
item $j$ is chosen (or ``wins''), and $C_{jk}$ is the number of comparisons between $j$
and $k$. \cref{ass:strong-connected} guarantees $C_{jk}>0$ for any $j,k$.  \citet{zermelo1929berechnung}
proved that under this assumption
$s^{(1)},s^{(2)},\dots$ converge to the unique maximum likelihood estimator,
and the sequence of log-likelihoods $\ell(s^{(1)}),\ell(s^{(2)}),\dots$
is monotone increasing. A cyclic version of \eqref{eq:pairwise} appeared in \citet{ford1957solution} with an independent proof of convergence. One can verify that by aggregating choice sets $S_i$ in \eqref{eq:scaling-iteration} over pairs of objects, it reduces to \eqref{eq:pairwise}. However, \eqref{eq:pairwise} as is written does not admit a matrix balancing formulation. A generalization of the algorithm of \citet{zermelo1929berechnung,ford1957solution,dykstra1956note} for pairwise comparison to ranking data was not achieved until the influential works of \citet{lange2000optimization} and \citet{hunter2004mm}. 
\subsection{MM Algorithm of \citet{hunter2004mm} for Ranking
Data} 
Motivated by the observation in \citet{lange2000optimization} that \eqref{eq:pairwise} is an instance of an minorization-maximization (MM) algorithm, the seminal work of \citet{hunter2004mm} proposed the general approach of solving ML estimation of choice models via MM algorithms, which relies on the inequality 
\begin{align*}
-\log x & \geq1-\log y-(x/y)
\end{align*}
to construct a lower bound (minorization) on the log-likelihood that has an explicit maximizer (maximization), and iterates between the two steps. \citet{hunter2004mm} develops such an algorithm for the Plackett--Luce
model for ranking data and proves its monotonicity and convergence.

Given $n$ partial rankings, where the $i$-th partial ranking on $l_{i}$ objects is indexed by $a(i,1)\rightarrow a(i,2)\rightarrow \cdots \rightarrow a(i,l_{i})$, the MM algorithm of \citet{hunter2004mm} takes the form 
\begin{align}
\label{eq:mm}
s_{k}^{(t+1)} & =\frac{w_{k}}{\sum_{i=1}^{n}\sum_{j=1}^{l_{i}-1}\delta_{ijk}[\sum_{j'=j}^{l_{i}}s_{a(i,j')}^{(t)}]^{-1}},
\end{align}
 where $\delta_{ijk}$ is the indicator that item $k$ ranks no better
than the $j$-th ranked item in the $i$-th ranking, and $w_{k}$ is the number of
rankings in which $k$ appears but is not ranked last. 
\begin{proposition}
\label{lem:mm}
Sinkhorn's algorithm applied to the ML estimation of the Plackett--Luce model is algebraically equivalent to \eqref{eq:mm}. 
\end{proposition}
Therefore, Sinkhorn's algorithm applied to the ML estimation of the Plackett--Luce model reduces algebraically to the MM algorithm of \citet{hunter2004mm}. However, \cref{alg:scaling} applies to more general choice models with minimal or no change, while the approach in \citet{hunter2004mm} requires deriving the minorization-majorization step for every new optimization objective. This was carried out, for example, in \citet{maystre2017choicerank} for a network choice model. We show in \cref{prop:choicerank} that their ChoiceRank algorithm is also a special case of (regularized) Sinkhorn's algorithm. From a computational perspective, even when algorithms are equivalent algebraically, their empirical performance can vary drastically depending
on the particular implementation. Another advantage of Sinkhorn's algorithm is that it computes \emph{all}
entries simultaneously through vector and matrix operations, while
the analytical formula in \eqref{eq:mm} is hard to parallelize. This distinction is likely behind
the discrepancy in \cref{sec:empirics} between our experiments and those in \citet{maystre2015fast},
who conclude that the MM version \eqref{eq:mm} is slower in terms of wall clock time than their Iterative Luce Spectral Ranking (I-LSR) algorithm for the Plackett--Luce model on $k$-way partial ranking data. 

\subsection{Markov Steady State Inversion and Network Choice}
\label{subsec:steady-state}
Our work is related to the works of \citet{kumar2015inverting,maystre2017choicerank} on Markov chains on graphs, where transition matrices are parameterized by node-dependent scores prescribed by Luce's choice axiom. More precisely, given a directed graph $G=(V,E)$ and $N^\out_j, N^\inn_j \subseteq V$ the neighbors with edges going out from and into $j\in V$, and a target stationary distribution $\pi$, the (unweighted) steady state inversion problem of \citet{kumar2015inverting} seeks scores $s_j$ such that the transition matrix $T_{j,k}=\frac{s_k}{\sum_{k'\in N^\out_j}s_{k'}}$ has the desired stationary distribution $\pi$. Their Theorem 13 shows that a bipartite version of this problem is equivalent to solving the ML estimation conditions \eqref{eq:optimality} of the choice model. Furthermore, one can verify that their bipartite inversion problem has the same form as \eqref{eq:bridge} in our paper, with the bipartite graph defined using $A$. Their existence condition (termed ``consistency'') is equivalent to \cref{ass:matrix-weak-existence}(a) \citep{menon1968matrix} for the matrix balancing problem. Despite these connections, the key difference in our work is the reformulation of \eqref{eq:optimality} as one involving diagonal scalings of rows and columns of $A$, which was absent in \citet{kumar2015inverting}. Consequently, they proposed a different algorithm instead of applying Sinkhorn's algorithm. 

Building on \citet{kumar2015inverting}, \citet{maystre2017choicerank} consider a similar Markov chain on $(V,E)$, where now for each edge $(j,k) \in E$ one observes a finite \emph{number} $c_{jk}$ of transitions along it, and take a maximum likelihood approach to estimate the scores $s_j$.  They show that, as one might expect, the steady state inversion problem of \citet{kumar2015inverting} is the asymptotic version of the ML estimation problem in their network choice model. 

An additional contribution of \citet{maystre2017choicerank} is the regularization of the inference problem via a Gamma prior on $s_j$'s, which eliminates the necessity of any assumptions on the choice dataset such as \cref{ass:strong-connected}. They then follow the proposal of \citet{hunter2004mm} and develop an MM algorithm for maximum likelihood estimation called ChoiceRank, the unregularized version of which can be written as follows:
\begin{align}
    \label{eq:choicrank}
    s_j^{(t+1)} = \frac{c_j^\inn}{\sum_{k\in N_j^\inn} \gamma_k^{(t)}}, \gamma_j^{(t)}=\frac{c_j^\out}{\sum_{k\in N_j^\out}s_k^{(t)}},
\end{align}
where $c_j^\inn=\sum_{k\in N_j^\inn}c_{kj}$ and $c_j^\out=\sum_{k\in N_j^\out}c_{jk}$ are the total number of observed transitions into and out of $j\in V$. 
\begin{proposition}
    \label{prop:choicerank}
   The network choice model of \citet{maystre2017choicerank} is a special case of the choice model \eqref{eq:model}, and  Sinkhorn's algorithm applied to this case reduces to an iteration algebraically equivalent to \eqref{eq:choicrank}.
\end{proposition}
We also explore the regularization approach of \citet{maystre2017choicerank} in \cref{subsec:regularization} and demonstrate that Sinkhorn's algorithm can easily accommodate this extension, resulting in a regularized version of Sinkhorn's algorithm for matrix balancing that \emph{always} converges. This is given in \cref{alg:regularized}. Once again, insights from choice modeling yield useful improvements in matrix balancing.

\subsection{Sinkhorn's Algorithm as a General MM Algorithm}
\label{sec:sinkhorn-MM}
That Sinkhorn's algorithm reduces to MM algorithms when applied to various choice models is not a coincidence. 
In this section, we establish the connection between choice modeling
and matrix balancing through an optimization perspective. This connection
provides an interesting interpretation of Sinkhorn's algorithm
as optimizing a dominating function, i.e., an MM algorithm. See \citet{lange2000optimization,lange2016mm} for a discussion of the general correspondence between block
coordinate descent algorithms and MM algorithms.

First, we discuss the connection between the log-likelihood function
\eqref{eq:log-likelihood} and the dual potential function \eqref{eq:log-barrier} when $(A,p,q)$ corresponds to a choice dataset. Consider maximizing the negative
dual potential function 
\begin{align*}
h(d^{0},d^{1}):=-g(d^{0},d^{1}) & =\sum_{j=1}^{m}q_{j}\log d_{j}^{0}+\sum_{i=1}^{n}p_{i}\log d_{i}^{1}-(d^{1})^{T}Ad^{0}.
\end{align*}
For each fixed $d^{0}$, the function is concave in $d^{1}$, and
maximization with respect to $d^{1}$ yields first order conditions
\begin{align*}
d^{1} & =p/(Ad^{0}).
\end{align*}
Substuting this back into $h(d^{0},d^{1})$, we obtain 
\begin{align*}
f(d^{0}):=h(d^{0},p/(Ad^{0})) & =\sum_{j=1}^{m}q_{j}\log d_{j}^{0}+\sum_{i=1}^{n}p_{i}\log(\frac{p_{i}}{(Ad^{0})_{i}})-\sum_{i}p_{i}\\
 & =\sum_{j=1}^{m}q_{j}\log d_{j}^{0}-\sum_{i=1}^{n}p_{i}\log(Ad^{0})_{i}+\sum_{i}p_{i}\log p_{i}-\sum_{i}p_{i}.
\end{align*}
If $A$ is a valid participation matrix for a choice dataset and $p,q$
are integers, we can identify $(A,p,q)$ with a choice dataset. Each
row of the participation matrix $A$ is the indicator vector of choice
set $S_{i}$, and $d_{j}^{0}$ is the quality score. In this
case $(Ad^{0})_{i}=\sum_{k\in S_{i}}d_{k}^{0}$, so that 
\begin{align*}
\sum_{j=1}^{m}q_{j}\log d_{j}^{0}-\sum_{i=1}^{n}p_{i}\log(Ad^{0})_{i} & =\sum_{j=1}^{m}q_{j}\log d_{j}^{0}-\sum_{i=1}^{n}p_{i}\sum_{k\in S_{i}}d_{k}^{0}=\ell(d^{0}).
\end{align*}

It then follows that
\begin{align*}
\min_{d^{0},d^{1}}g(d^{0},d^{1})\Leftrightarrow\max_{d^{0},d^{1}}h(d^{0},d^{1})\Leftrightarrow\max_{d^{0}}\max_{d^{1}}h(d^{0},d^{1})\Leftrightarrow\max_{d^{0}}f(d^{0})\Leftrightarrow\max_{d^{0}}\ell(d^{0}),
\end{align*}
so that minimizing the potential function $g$ is equivalent to maximizing
the log-likelihood function $\ell$. Moreover, the first order condition of maximizing $h$ with respect to $d^1$ is 
$d^0=q/(A^Td^1)$, which when $(A,p,q)$ is identified with a choice dataset reduces to 
\begin{align*}
    q_j =\sum_{i\mid j\in S_i}p_i \frac{d^0_j}{\sum_{k\in S_i}d^0_k},
\end{align*}
which is the optimality condition \eqref{eq:optimality} of the choice model.

Next, given $d^{0(t)}$ the estimate of $d^{0}$ after the
$t$-th iteration, define the function 
\begin{align*}
f(d^{0}\mid d^{0(t)}):=h(d^{0},p/(Ad^{0(t)})) & =\sum_{j=1}^{m}q_{j}\log d_{j}^{0}+\sum_{i=1}^{n}p_{i}\log\frac{p}{Ad^{0(t)}}-\frac{p_{i}}{(Ad^{0(t)})_{i}}(Ad^{0})_{i}.
\end{align*}
We can verify that 
\begin{align*}
f(d^{0(t)}\mid d^{0(t)}) & =f(d^{0(t)})\\
f(d^{0}\mid d^{0(t)}) & \leq f(d^{0}),
\end{align*}
 so that $f(d^{0}\mid d^{0(t)})$ is a valid minorizing function of
$f(d^{0})$ \citep{lange2000optimization} that guarantees the ascent property $f(d^{0(t+1)})\geq f(d^{0(t+1)}\mid d^{0(t)})=\max_{d^{0}}f(d^{0}\mid d^{0(t)})\geq f(d^{0(t)}\mid d^{0(t)})=f(d^{0(t)})$.
The update in the maximization step
\begin{align*}
d^{0(t+1)}=\arg\max_{d^{0}}f(d^{0}\mid d^{0(t)}) & =q/A^{T}\frac{p}{Ad^{0(t)}}
\end{align*}
 is precisely one full iteration of Sinkhorn's algorithm. Note that this interpretation of Sinkhorn's algorithm does not require $A$ to be binary and $p,q$ to be integers.

On the other hand, using the property $-\ln x\geq1-\ln y-(x/y)$,
we can \emph{directly} construct a minorizing function of $\ell$
by 
\begin{align*}
\ell(d^{0})=\sum_{j=1}^{m}q_{j}\log d_{j}^{0}-\sum_{i=1}^{n}p_{i}\log(Ad^{0})_{i} & \geq\sum_{j=1}^{m}q_{j}\log d_{j}^{0}+\sum_{i=1}^{n}p_{i}(-\frac{(Ad^{0})_{i}}{(Ad^{0(t)})_{i}}-\log(Ad^{0(t)})_{i}+1)=\ell(d^{0}\mid d^{0(t)}),
\end{align*}
 where $\ell(d^{0}\mid d^{0(t)})$ is a valid minorizing function
of $\ell$. Maximizing $\ell(d^{0}\mid d^{0(t)})$ with respect to
$d^{0}$, the update in the maximization step is
\begin{align*}
d_{j}^{0(t+1)} & =q_{j}/\sum_{i\mid j\in S_i}\frac{p_{i}}{(Ad^{0(t)})_{i}}
\end{align*}
 which again is one full iteration of Sinkhorn's algorithm applied to the Luce choice model. Moreover,
\begin{align*}
\ell(d^{0}\mid d^{0(t)})+\sum_{i}p_{i}\log p_{i}-\sum_{i}p_{i} & =\sum_{j=1}^{m}q_{j}\log d_{j}^{0}+\sum_{i=1}^{n}p_{i}(-\frac{(Ad^{0})_{i}}{(Ad^{0(t)})_{i}}-\log(Ad^{0(t)})_{i}+1)+\sum_{i}p_{i}\log p_{i}-\sum_{i}p_{i}\\
 & =\sum_{j=1}^{m}q_{j}\log d_{j}^{0}+\sum_{i=1}^{n}p_{i}(-\frac{(Ad^{0})_{i}}{(Ad^{0(t)})_{i}}+\log\frac{p_{i}}{(Ad^{0(t)})_{i}})\\
 & =f(d^{0}\mid d^{0(t)}).
\end{align*}
 Therefore, the minorizing function $\ell(d^{0}\mid d^{0(t)})$ constructed
using $-\ln x\geq1-\ln y-(x/y)$ for the log-likelihood and the minorizing
function $f(d^{0}\mid d^{0(t)})$ constructed for $\max_{d^{1}}h(d^{0},d^{1})$
are identical modulo a constant $\sum_{i}p_{i}\log p_{i}-\sum_{i}p_{i}$.
Sinkhorn's algorithm is in fact the MM algorithm corresponding to both minorizations. However, the perspective using $f(d^{0}\mid d^{0(t)})$ is more general since it applies to general $(A,p,q)$ as long as they satisfy Assumptions \cref{ass:matrix-existence} and \cref{ass:matrix-uniqueness}, whereas the MM algorithm based on $\ell(d^{0}\mid d^{0(t)})$ is designed for choice dataset, so requires $A$ to be binary.
{\subsection{Sinkhorn's Algorithm as a Partial Minimization Algorithm}}

Here we show that the Sinkhorn method can be also view a partial minimization
algorithm. Recall that
\[ g (d^0, d^1) = (d^1)^T A d^0 - \sum_{i = 1}^n p_i \log d_i^1 - \sum_{i =
   1}^m q_i \log d_j^0 \]
and we have $d^1_{\ast} = p / (A d^0_{\ast})$. Replacing $d^1$ by $p / (A
d^0)$, we get
\[ h (d^0) = \sum_{i = 1}^n p_i - \sum_{i = 1}^n p_i \log \frac{p_i}{a_i^T
   d^0} - \sum_{i = 1}^m q_i \log d_j^0 = C + \sum_{i = 1}^n p_i \log a_i^T
   d^0 - \sum_{i = 1}^m q_i \log d_j^0 . \]
where $a_i$ is the $i$-th row of $A$ and $C$ is some constant. Denote $h_1
(d^0) := \sum_{i = 1}^n p_i \log a_i^T d^0$ and $h_2 (d^0) := -
\sum_{i = 1}^m q_i \log d_j^0$. We can view each Sinkhorn iteration as
\[ d^{0 (t + 1)} = \arg \min_{d^0}  \{ h_1 (d^{0 (t)}) + \nabla h_1 (d^{0
   (t)})^T (d^{0 (t + 1)} - d^{0}) + h_2 (d^0) \}, \]
where $h_1 (d^0)$ is replaced by the first-order Taylor expansion and the
gradient of $h_1$ at $d^{0 (t)}$ is $\nabla h_1 (d^{0 (t)}) = p / A d^0 = d^{1
(t + 1)}$. Hence in terms of $d^0$ or $d^1$ alone, Sinkhorn algorithm can be
viewed as a partial minimization method: in our example $d^{1 (t + 1)}$ is
exactly the gradient of $h_1 (d^{0 (t)})$.

\subsection{Sinkhorn's Algorithm and Distributed Optimization} 
We now shift our focus to algorithms in distributed optimization, where Sinkhorn's algorithm can be interpreted as a message passing/belief propagation algorithm \citep{balakrishnan2004polynomial}. We start by observing a connection to the ASR algorithm for estimating Luce choice models \citep{agarwal2018accelerated}, which returns the same approximate ML estimators as the RC \citep{negahban2012iterative} and LSR \citep{maystre2015fast} algorithms, but has provably faster convergence. 

Consider the bipartite graph $G_b$ defined by $A$ in \cref{subsec:graph-laplacian}, which consists of choice set nodes $V$ on one hand and item nodes $U$ on the other, where there is an edge between $i\in V$ and $j\in U$ if and only if $j\in S_i$. \citet{agarwal2018accelerated} provide the following message passing interpretation of ASR on the bipartite graph: at every iteration, the item nodes send a ``message'' to their
neighboring choice set nodes consisting of each item node's current estimate of their own $s_j$; the choice set nodes then aggregate the messages they receive by summing up these estimates, and then sending back the sums to their neighboring item nodes. The item nodes use these sums to update estimates of their own $s_j$. \citet{agarwal2018accelerated} show that since the ASR algorithm is an instance of the message passing algorithm, it can be implemented in a distributed manner.% similar to the parallelization capability of Sinkhorn's algorithm. 

We now explain how Sinkhorn's algorithm is another instance of the message passing algorithm described above. Recall that $d^0$ is identified with the $s_j$'s in the Luce choice model, so that $d^{1} \leftarrow p/(Ad^{0})$ precisely corresponds to item nodes ``passing'' their current estimates to set nodes, which then sum up the received estimates and then take the weighted \emph{inverse} of this sum. Similarly, $d^0 \leftarrow q/(A^{T} d^{1})$ corresponds to choice set nodes passing their current estimates of $d^1$ back to item nodes, which then sum up the received messages and take the weighted inverse as their updated estimates of $s_j$. The main difference with ASR lies in how each item node $j$ updates its estimate of $s_j$ based on the messages it receives from neighboring set nodes. In Sinkhorn's algorithm, the update to $s_j$ is achieved by dividing $p$ by a weighted average of the \emph{inverse} of summed messages $1/\sum_{k\in S_{i}}s_{k}^{(t)}$:
\begin{align*}
  s_j^{(t+1)} \leftarrow  q_j/(A^{T} d^{1})_j=W_j/\sum_{i\in[n]\mid j\in S_{i}}\frac{R_i}{\sum_{k\in S_{i}}s_{k}^{(t)}},
\end{align*}
whereas in ASR, the update is an average of the summed messages $\sum_{k\in S_{i}}s_{k}^{(t)}$ without taking their inverses first: 
\begin{align*}
    s^{(t+1)}_j \leftarrow  \frac{1}{\sum_{i\in[n]\mid j\in S_{i}} R_i} \sum_{i\in[n]\mid j \in S_i}  W_{ji} \sum_{k\in S_{i}}s_{k}^{(t)},
\end{align*}
where $W_{ji}$ is the number of times item $j$ is selected from all observations having choice set $S_i$, with $\sum_i W_{ji} = W_j$.

From another perspective, the two algorithms arise from different \emph{moment} conditions. While Sinkhorn's algorithm is based on the optimality condition \eqref{eq:optimality}, ASR is based on the condition
\begin{align*}
    \sum_{i\in[n]\mid j\in S_{i}} R_i = \sum_{i\in[n]\mid j \in S_i}  W_{ji}/\frac{s_j}{\sum_{k\in S_{i}}s_{k}},
\end{align*}
which results in an approximate instead of exact MLE. In this regard, the algorithm of \citet{newman2023efficient} can be viewed as based on an alternative (but equivalent) form of the moment condition derived from ML estimation.

The message passing interpretation also provides further insights on the importance of algebraic connectivity to the convergence rate of Sinkhorn's algorithm. Graph theoretic conditions like \cref{ass:strong-connected} are related to network flow and belief propagation, and characterize how fast information can be distributed across the bipartite network with the target distributions $p,q$. It is well-known that convergence of distributed algorithms on networks depends critically on the network topology through the spectral gap of the associated averaging matrix. We can understand \cref{thm:convergence} on the asymptotic convergence rate of Sinkhorn's algorithm as a result of this flavor, although a precise equivalence is left for future works. 

\subsection{The Berry--Levinsohn--Pakes Algorithm}
Last but not least, our work is also closely related to the economics literature that studies consumer behavior based on discrete choices \citep{mcfadden1973conditional,mcfadden1978modelling,mcfadden1981econometric,berry1995automobile}. Here we discuss the particular connection with the work of Berry,
Levinsohn, and Pakes \cite{berry1995automobile}, often referred to as BLP. To estimate consumer preferences over automobiles across different markets (e.g., geographical), they propose a random utility (RUM) model indexed by individual $i$, product
$j$, and market $t$:
\begin{align*}
U_{ijt} & =\beta_{i}^{T}X_{jt}+\theta_{jt}+\epsilon_{ijt},
\end{align*}
 where $\theta_{jt}$ is an unobserved product characteristic, such as the overall popularity of certain types of cars in different regions, and $\epsilon_{ijt}$
are \emph{i.i.d.} double exponential random variables.
% \begin{align*}
% f(\epsilon) & =\exp(-\epsilon-\exp(-\epsilon))
% \end{align*}
The individual-specific coefficient $\beta_{i}$ is random with 
\begin{align*}
\beta_{i} & =Z_{i}^T\Gamma+\eta_{i}\\
\eta_{i}\mid Z_{i} & \sim\mathcal{N}(\beta,\Sigma),
\end{align*}
and the observations consist of \emph{market shares} $\hat{p}_{jt}$ of each
product $j$ in market $t$ and observable population characteristics $Z_{i}$ in each
market. Given a model with fixed $\beta,\Gamma,\Sigma$ and observations, the task is to estimate $\theta_{jt}$. 

For every value of $\theta_{jt}$, we can compute, or simulate if necessary, the \emph{expected} market shares $p_{jt}$, which is the likelihood of product $j$ being chosen in market $t$. For example, in the special case that $\beta,\Gamma,\Sigma \equiv 0$ and $\exp(\theta_{jt})\equiv s_{j}$ for all
$j,t$, i.e., (perceived) product characteristic does not vary across
markets, the expected market share reduces to the familiar formula
\begin{align*}
p_{jt} & =\frac{\exp(\theta_{jt})}{\sum_{k}\exp(\theta_{kt})}=\frac{s_{j}}{\sum_{k}s_{k}}.
\end{align*}
The generalized method of moments (GMM) approach of \citet{berry1995automobile} is to find $\theta_{jt}$ such that $p_{jt}=\hat{p}_{jt}$, i.e., the implied expected market share equals the observed share. Recall the similarity to the optimality condition \eqref{eq:optimality} of the Luce choice model. BLP propose the iteration 
\begin{align}
\label{eq:blp}
\theta_{jt}^{(m+1)} & =\theta_{jt}^{(m)}+\log\hat{p}_{jt}-\log p_{jt}(\theta^{(m)},\beta,\Gamma,\Sigma),
\end{align}
 and show that it is a \emph{contraction mapping}, whose fixed point is the desired estimates of $\theta_{jt}$. 
\begin{proposition}
    \label{prop:BLP}
    When $\beta,\Gamma,\Sigma \equiv 0$ and $\exp(\theta_{jt})\equiv s_{j}$, the GMM condition of BLP on market shares is equivalent to the optimality condition \eqref{eq:optimality} for a Luce choice model where all alternatives are available in every observation. Furthermore, the BLP algorithm is equivalent to Sinkhorn's algorithm in this model.
\end{proposition}
For a more detailed correspondence, see \citet{bonnet2022yogurts}. Importantly, \cref{prop:BLP} does not imply that the Luce choice model and Sinkhorn's algorithm is a strict special case of the BLP framework. The key difference is that BLP, and most discrete choice models in econometrics, implicitly assumes that the entire set of alternatives is always available in each observation. This assumption translates to a participation matrix $A$ in \eqref{eq:equation-system} that has 1's in all entries. In this setting, the MLE of $s_j$ is simply the empirical winning frequencies. %Note that \cref{ass:strong-connected} may still be violated in this case, highlighting again the gap between what is required for the ML estimation problem and by the matrix scaling problem.
On the other hand, while the Luce choice model allows different choice sets $S_i$ across observations, they do not include covariate information on the alternatives or decision makers, which is important in discrete choice modeling. One can reconcile this difference by relabeling alternatives with different covariates as \emph{distinct} items, and we leave investigations on further connections in this direction to future works.

%% file: regularization.tex
\section{Regularization of Luce Choice Models and Matrix Balancing Problems}
\label{subsec:regularization}
In practice, many choice and ranking datasets may not satisfy \cref{ass:strong-connected}, which is required for the maximum likelihood estimation to be well-posed. Equivalently, for the matrix balancing problem, when a triplet $(A,p,q)$ does not satisfy \cref{ass:matrix-existence} and \cref{ass:matrix-uniqueness}, no finite scalings exist and Sinkhorn's algorithm may diverge. In this section, we discuss some regularization techniques to address these problems. They are easy to implement and require minimal modifications to Sinkhorn's algorithm. Nevertheless, they can be very useful in practice to regularize ill-posed problems. Given the equivalence between the problem of computing the MLE of Luce choice models and the problem of matrix balancing, our proposed regularization methods apply to both. 
\subsection{Regularization via Gamma Prior}
 As discussed in \cref{sec:formulations},  for a choice dataset to have a well-defined maximum likelihood estimator, it needs to satisfy \cref{ass:strong-connected}, which requires the directed comparison graph to be strongly connected. Although this condition is easy to verify, the question remains what one can do in case it does \emph{not} hold. As one possibility, we may introduce a prior on the parameters $s_j$, which serves as a regularization of the log-likelihood that results in a unique maximizer. Many priors are possible. For example, \citet{maystre2017choicerank}, following \citet{caron2012efficient}, use independent Gamma priors on $s_j$. In view of the fact that the unregularized problem and algorithm in \citet{maystre2017choicerank} is a special case of the Luce choice model and Sinkhorn's algorithm, we can also incorporate the Gamma prior to the Luce choice model \eqref{eq:log-likelihood} to address identification problems. 

More precisely, suppose now that each $s_j$ in the Luce choice model are i.i.d. Gamma$(\alpha,\beta)\propto s_j^{\alpha-1}e^{-\beta s_j}$. This leads to the following regularized log-likelihood:
\begin{align}
\label{eq:regularized-log-likelihood}
\ell^R(s):=\sum_{i=1}^{n}\log s_{j_i}-\log\sum_{k\in S_{i}}s_{k} +(\alpha-1)\sum_{j=1}^{m}\log s_j - \beta \sum_{j=1}^{m} s_j.
\end{align}
The corresponding first order condition is given by
\begin{align}
\label{eq:optimality-regularized}
\frac{W_j+\alpha-1}{n} & = \frac{1}{n} \left(\sum_{i\in[n]\mid j\in S_{i}}R_i \frac{s_{j}}{\sum_{k\in S_{i}}s_{k}} +\beta s_j\right),
\end{align}
which leads to the following modified Sinkhorn's algorithm, which generalizes the ChoiceRank algorithm of \citet{maystre2017choicerank}: 
\begin{align}
\label{eq:regularized-sinkhorn}
 d^0 \leftarrow (q+\alpha -1)/ (A^Td^1+\beta),  \quad d^1 \leftarrow p/A d^0.
\end{align}
 The choice of $\beta$ determines the normalization of $s_j$. With $u_j=\log s_j$, we can show in a similar way as Theorem 2 of \citet{maystre2017choicerank} that \eqref{eq:regularized-log-likelihood} always has a unique maximizer whenever $\alpha>1$. Regarding the convergence, \citet{maystre2017choicerank} remarked that since their ChoiceRank algorithm can be viewed as an MM algorithm, it inherits the local linear convergence of MM algorithms \citep{lange2000optimization}, but ``a detailed investigation of convergence behavior is left for future works''. With the insights we develop in this paper, we can in fact provide an explanation for the validity of the Gamma priors from an optimization perspective. This perspective allows us to conclude directly that \eqref{eq:regularized-log-likelihood} always has a unique solution in the interior of the probability simplex, and that furthermore the iteration in \eqref{eq:regularized-sinkhorn} has global linear convergence. 
 Consider now the following regularized potential function
 \begin{align}
     g^R(d^0,d^1)	:= ((d^1)^{T}A+\beta (\mathbf{1}_m)^T)d^0-\sum_{i=1}^{n}p_{i}\log d^1_{i}-\sum_{j=1}^{m}(q_{j}+\alpha-1)\log d^0_{j}.
\end{align}
We can verify that by substituting the optimality condition of $d^1$ into $-g^R$, it reduces to the log posterior \eqref{eq:regularized-log-likelihood}. Moreover, the iteration \eqref{eq:regularized-sinkhorn} is precisely the alternating minimization algorithm for $g^R$. When $\alpha-1,\beta>0$, the reparameterized potential function
\begin{align*}
    \sum_{ij}(e^{-v_{i}}A_{ij}e^{u_{j}})+\beta\sum_{j}e^{u_{j}}+\sum_{i}p_{i}v_{i}-\sum_{j}(q_{j}+\alpha-1)u_{j}
\end{align*} 
is always coercive regardless of whether \cref{ass:matrix-existence} holds. Therefore, during the iterations \eqref{eq:regularized-sinkhorn}, $(u,v)$ stay \emph{bounded}. Moreover, the Hessian is 
\begin{align*} \begin{bmatrix}\sum_{j}e^{-v_{i}}A_{ij}e^{u_{j}} & -e^{-v_{i}}A_{ij}e^{u_{j}}\\
-e^{-v_{i}}A_{ij}e^{u_{j}} & \sum_{i}e^{-v_{i}}A_{ij}e^{u_{j}}+\beta e^{u_{j}}
\end{bmatrix}	\succ0,
\end{align*}
 which is now positive definite.  As a result, $g^{R}(u,v)$ is strongly convex and smooth, so that  \eqref{eq:regularized-sinkhorn} converges linearly. From the perspective of the matrix balancing problem, we have thus obtained a regularized version of Sinkhorn's algorithm, summarized in \cref{alg:regularized}, which is guaranteed to converge linearly to a finite solution $(D^1,D^0)$, even when \cref{ass:matrix-existence} does not hold for the triplet $(A,p,q)$. Moreover, the regularization also improves the convergence of Sinkhorn's algorithm even when it converges, as the regularized Hessian becomes more-behaved. This regularized algorithm could be very useful in practice to deal with real datasets that result in slow, divergent, or oscillating Sinkhorn iterations.
 \begin{algorithm}[tb]
\caption{Regularized Sinkhorn's Algorithm}
   \label{alg:regularized}
\begin{algorithmic}
   \STATE {\bfseries Input:}  $A, p, q,\alpha>1,\beta>0,\epsilon_{\text{tol}}$.
   \STATE {\bfseries initialize} $d^{0}\in\mathbb{R}_{++}^{m}$
   \REPEAT
   \STATE $d^{1} \leftarrow  p/( A d^0)$ 

   \STATE $d^{0}\leftarrow  (q+\alpha-1)/({A}^{T} d^{1}+\beta)$

   \STATE 
$\epsilon\leftarrow$  update of $(d^{0},d^1)$
\UNTIL{$\epsilon<\epsilon_{\text{tol}}$}
\end{algorithmic}
\end{algorithm}

\subsection{Regularization via Data Augmentation}
The connection between Bayesian methods and \emph{data augmentation} motivates us to also consider direct data augmentation methods. This is best illustrated in the choice modeling setting. Suppose for a choice dataset we construct participation matrix $A$, $p$ the counts of distinct choice sets, and $q$ the counts of each item being selected. We know that $(A,p,q)$ has a finite scaling solution if and only if  \cref{ass:strong-connected} holds, i.e., the directed comparison graph is strongly connected. We now propose the following modification of $(A,p,q)$ such that the resulting problem is always valid. 

First, if $A$ does not already contain a row equal to $\mathbf{1}_m^T$, i.e., containing all 1's, add this additional row to $A$. Call the resulting matrix $A'$. Then, expanding the dimension of $p$ if necessary, add $m\epsilon$ to the entry corresponding to $\mathbf{1}_m$, where we can assume for now that $\epsilon\geq 1$ is an integer.
This procedure effectively adds $m \epsilon$ ``observations'' that contain all $m$ items. For these additional observations, we let each item be selected exactly $\epsilon$ times. Luce's choice axiom guarantees that the exact choice of each artificial observation is irrelevant, and we just need to add $\epsilon \mathbf{1}_m$ to $q$. This represents augmenting each item with an additional $\epsilon$ ``wins'', resulting in the triplet $(A',p+(m\epsilon)\mathbf{e},q+\epsilon \mathbf{1}_m)$, where $\mathbf{e}$ is the one-hot indicator of the row $\mathbf{1}_m^T$ in $A'$. Now by construction, in any partition of $[m]$ into two non-empty subsets, any item from one subset is selected at least $\epsilon$ times over any item from the other subset. Therefore, \cref{ass:strong-connected} holds, and the maximum likelihood estimation problem, and equivalently the matrix balancing problem with $(A',p+(m\epsilon)\mathbf{e},q+\epsilon \mathbf{1}_m)$, is well-defined. This regularization method applies more generally to any non-negative $A$, even if it is not a participation matrix, i.e., binary. Although in the above construction based on choice dataset, $\epsilon$ is taken to be an integer, for the regularized matrix balancing problem with $(A',p+(m\epsilon)\mathbf{e},q+\epsilon \mathbf{1}_m)$, we can let $\epsilon \rightarrow 0$.

%% file: empirical.tex
\section{Applications of Matrix Balancing}
\label{app:related-works}
This section contains a brief survey on the applications of matrix balancing in a diverse range of disciplines.

\textbf{Traffic and Transportation Networks.} These applications are some of the earliest and most popular uses of matrix balancing. \citet{kruithof1937telefoonverkeersrekening} considered the problem of estimating new telephone traffic patterns among telephone exchanges given existing traffic volumes and marginal densities of departing and terminating traffic for each exchange when their subscribers are updated. A closely related problem in transportation networks is to use observed \emph{total} traffic flows out of each origin and into each destination to estimate \emph{detailed} traffics between origin-destination pairs \citep{carey1981method,nguyen1984estimating,sheffi1985urban,chang2021mobility}. The key idea is to find a traffic assignment satisfying the total flow constraints that is ``close'' to some known reference traffic pattern. The resulting (relative) entropy minimization principle, detailed in \eqref{eq:relative-entropy-minimization}, is an important optimization perspective on Sinkhorn's algorithm. Leveraging the connections between the maximum likelihood objective of Luce choice models and the dual objective of the relative entropy minimization problem, 
\citet{chang2024inferring} propose a statistical model for network traffic. The maximum likelihood estimate of the resulting model is precisely the optimal solution to a matrix balancing algorithm. They also quantify the statistical efficiency of the MLE in terms of the algebraic connectivity of the bipartite graph defined in \eqref{eq:bipartite-adjacency}. 

\textbf{Demography.} A problem similar to that in networks arises in demography. Given out-of-date inter-regional migration statistics and up-to-date net migrations from and into each region, the task is to estimate migration flows that are consistent with the marginal statistics \citep{plane1982information}.

\textbf{Economics.}
General equilibrium models in economics employ \emph{social accounting matrices}, which record the flow of funds between important (aggregate) agents in an economy at different points in time \citep{stone1962multiple,pyatt1985social}. Often accurate data is available on the total expenditure and receipts for each agent, but due to survey error or latency, detailed flows are not always consistent with these marginal statistics. Thus they need to be ``adjusted'' to satisfy consistency requirements.  Other important applications of the matrix balancing problem in economics include the estimation of gravity equations in inter-regional and international trade \citep{uribe1966information,wilson1969use,anderson2003gravity,silva2006log} and coefficient matrices in input-output models \citep{leontief1965structure,stone1971computable,bacharach1970biproportional}. In recent years, optimal transport \citep{villani2009optimal} has found great success in economics \citep{carlier2016vector,galichon2018optimal,galichon2021unreasonable,galichon2021matching}. As matrix balancing and Sinkhorn's algorithm are closely connected to optimal transport (\cref{sec:linear-convergence}), they are likely to have more applications in economics. %Indeed, we show in \cref{sec:connections} that the well-known algorithm of \citet{berry1995automobile} for estimating random utility choice models in economics is related to Sinkhorn's algorithm. 

\textbf{Statistics.}
A contingency table encodes frequencies of subgroups of populations, where the rows and columns correspond to values of two categorical variables, such as gender and age. Similar to social accounting matrices, a common problem is to adjust out-of-date or inaccurate cell values of a table given accurate marginal frequencies. The problem is first studied by  \citet{deming1940least}, who proposed the classic iterative algorithm. \citet{ireland1968contingency} formalized its underlying entropy optimization principle, and \citet{fienberg1970iterative} analyzed its convergence. 

\textbf{Optimization and Machine Learning.} Matrix balancing plays a different but equally important role in optimization. Given a linear system $Ax=b$ with non-singular $A$, it is well-known that the convergence of first order solution methods depends on the \emph{condition number} of $A$, and an important problem is to find diagonal \emph{preconditioners} $D^1,D^0$ such that $D^1AD^0$ has smaller condition number. Although it is possible to find optimal diagonal preconditioners via semidefinite programming \citep{boyd1994linear,qu2022optimal}, matrix balancing methods remain very attractive heuristics due to their low computational costs, and continue to be an important component of modern workhorse optimization solvers \citep{ruiz2001scaling,bradley2010algorithms,knight2013fast,stellato2020osqp,gao2022hdsdp}.

In recent years, optimal transport distances have become an important tool in machine learning and optimization for measuring the similarity between probability distributions \citep{arjovsky2017wasserstein,peyre2019computational,blanchet2019robust,mohajerin2018data,kuhn2019wasserstein}. Besides appealing theoretical properties, efficient methods to approximate them in practice have also contributed to their wide adoption. This is achieved through an entropic regularization of the OT problem, which is precisely equivalent to the matrix balancing problem and solved via Sinkhorn's algorithm \citep{cuturi2013sinkhorn,altschuler2017near,dvurechensky2018computational}.  %Although various algorithms have since been proposed that enjoy better theoretical runtime complexities \citep{dvurechensky2018computational,blanchet2018towards,ge2019interior,jambulapati2019direct,chambolle2022accelerated,lin2022efficiency,li2023fast}, Sinkhorn-based approximations remain popular in practice despite sub-optimal theoretical bounds. The linear convergence results in this paper reveal that Sinkhorn's algorithm actually enjoys the same state-of-the-art $\tilde{\mathcal{O}}(n^2/\varepsilon)$ runtime bound achieved by other more complicated \emph{first-order} algorithms.  

\textbf{Political Representation.} The apportionment of representation seats based on election results has found unexpected solution in matrix balancing. A standard example consists of the matrix recording the votes each party received from different regions. The marginal constraints are that each party's total number of seats be proportional to the number of votes they receive, and similar for each region. A distinct feature of this problem is that the final apportionment matrix must have integer values, and variants of the standard algorithm that incorporate \emph{rounding} have been proposed \citep{balinski2006matrices,pukelsheim2006current,maier2010divisor}. More than just mathematical gadgets, they have found real-world implementations in Swiss cities such as Zurich \citep{pukelsheim2009iterative}.

\textbf{Markov Chains.}
Last but not least, Markov chains and related topics offer another rich set of applications for matrix balancing. \citet{schro1931uber} considered a continuous version of the following problem. Given a ``prior'' transition matrix $A$ of a Markov chain and \emph{observed} distributions $p^0,p^1$ before and after the transition, find the most probable transition matrix (or path) $\hat A$ that satisfies $\hat A p^0=p^1$. This is a variant of the matrix balancing problem and has been studied and generalized in a long line of works \citep{fortet1940resolution,beurling1960automorphism,ruschendorf1995convergence,gurvits2004classical,georgiou2015positive,friedland2017schrodinger}. Applications in marketing estimate customers' transition probabilities between different brands using market share data
\citep{theil1966quadratic}. Coming full circle back to  choice modeling, matrix balancing has also been used to rank nodes of a network. \citet{knight2008sinkhorn} explains how the inverses of left and right scalings of the adjacency matrix with uniform target marginals (stationary distributions) can be naturally interpreted as measures of their ability to attract and emit traffic. This approach is also related to the works of \citet{lamond1981bregman,kleinberg1999authoritative,tomlin2003new}.

\section{Numerical Experiments}
\label{sec:empirics}
In this section, we study the empirical performance of Sinkhorn's algorithm relative to a few popular alternatives. In addition, we present simulation results that illustrate the magnitude of our complexity bounds on randomly generated problem instances.\footnote{The replication package is available at \url{https://github.com/zhaonanq/choice-Sinkhorn}} 
\subsection{Performance on Real Datasets}
We compare the empirical performance of Sinkhorn's algorithm with the iterative LSR (I-LSR) algorithm of \citet{maystre2015fast} on real choice datasets. Because the implementation of I-LSR by \citet{maystre2015fast}
only accommodates pairwise comparison data and partial ranking data, but does not easily generalize to multi-way choice data, we focus on data with pairwise comparisons. 

We use the
natural parameters $\log s_{j}$ (logits) instead of $s_{j}$ when computing and evaluating
the updates, as the probability of $j$ winning
over $k$ is proportional to the ratio $s_{j}/s_{k}$, so that $s_j$ are usually logarithmically spaced. To make sure that estimates
are normalized, we impose the normalization that $\sum_{j}s_{j}=m$,
the number of objects, at the end of each iteration, although due to
the logarithm scale of the convergence criterion, the choice of normalization
does not seem to significantly affect the performances of the algorithms.
%generated based on the proposal in \citet{agarwal2018accelerated}, and partial ranking data, generated with code from \citet{maystre2015fast}.

We evaluate the algorithms on five real-world datasets consisting
of partial ranking or pairwise comparison data. The NASCAR dataset
consists of ranking results of the 2002 season NASCAR races. The SUSHI
datasets consist of rankings of sushi items. The Youtube dataset consists
of pairwise comparisons between videos and which one was considered
more entertaining by users. The GIFGIF dataset similarly consists
of pairwise comparisons of GIFs that are rated based on which one
is closer to describing a specific sentiment, such as happiness and
anger. We downsampled the Youtube and GIFGIF datasets due to memory
constraints. 
\begin{table*}[t]
\begin{centering}
\hspace*{-3cm}%
\begin{tabular}{c c c c c c c c c}
\hline 
\multirow{2}{*}{dataset} & \multirow{2}{*}{data type} & \multirow{2}{*}{items} & \multirow{2}{*}{observations} & \multirow{2}{*}{$k$} & \multicolumn{2}{c}{Sinkhorn} & \multicolumn{2}{c}{I-LSR}\tabularnewline
\cline{6-9} \cline{7-9} \cline{8-9} \cline{9-9} 
 &  &  &  &  & iterations & time & iterations & time\tabularnewline
\hline 
\hline 
NASCAR & $k$-way ranking & 83 & 36 & 43 & 20 & \textbf{0.029} & 13 & 0.960\tabularnewline
\hline 
SUSHI-10 & $k$-way ranking & 10 & 5000 & 10 & 16 & \textbf{0.025} & 8 & 1.708\tabularnewline
\hline 
SUSHI-100 & $k$-way ranking & 100 & 5000 & 10 & 21 & \textbf{0.253} & 9 & 2.142\tabularnewline
\hline 
Youtube & pairwise comparison & 2156 & 28134 & 2 & 89 & \textbf{7.984} & 33 & 15.026\tabularnewline
\hline 
GIFGIF & pairwise comparison & 2503 & 6876 & 2 & 1656 & \textbf{30.97} & 315 & 138.63\tabularnewline
\hline 
\end{tabular}\hspace*{-3cm}
\par\end{centering}
\caption{Performance of iterative ML inference algorithms on five real datasets. Youtube and GIFGIF data were subsampled. 
Convergence is declared when the maximum entry-wise change of an update
is less than $10^{-8}$. At convergence, the ML estimates returned
by the algorithms have entry-wise difference of at most $10^{-10}$.}
\label{tab:empirical}
\end{table*}
In \cref{tab:empirical} we report the running time of the three algorithms on different
datasets. Convergence is declared when the maximum entry-wise change
of an update to the natural parameters $\log s_{i}$ is less than
$10^{-8}$. At convergence, the MLEs returned by the algorithms have entry-wise difference of at most $10^{-10}$. %Since
% we know that the sequence of MM algorithm must converge to a stationary
% point, and that the matrix scaling algorithm is equivalent to the
% MM algorithm for pairwise comparisons and partial rankings, we know
% that the final estimates must be the MLE. 
We see that Sinkhorn's algorithm consistently outperforms
the I-LSR algorithm in terms of convergence speed. It also has the additional advantage of being parallelized with
elementary matrix-vector operations, whereas the iterative I-LSR algorithm needs to repeatedly
compute the steady-state of a continuous-time Markov chain, which
is prone to problems of ill-conditioning. This also explains why Sinkhorn's algorithm may take more iterations but has better wall clock time, since each iteration is much less costly. On the other hand, we note
that for large datasets, particularly those with a large number of
observations or alternatives, the dimension of $A$ used
may become too large for the memory of a single machine. If this is still a problem after removing duplicate rows and columns according to \cref{sec:equivalence}, we can use distributed implementations of Sinkhorn's algorithm, which in view of its connections to message passing algorithms, is a standard procedure.
\begin{figure}[hbt!]
\centering
\includegraphics[scale=0.3]{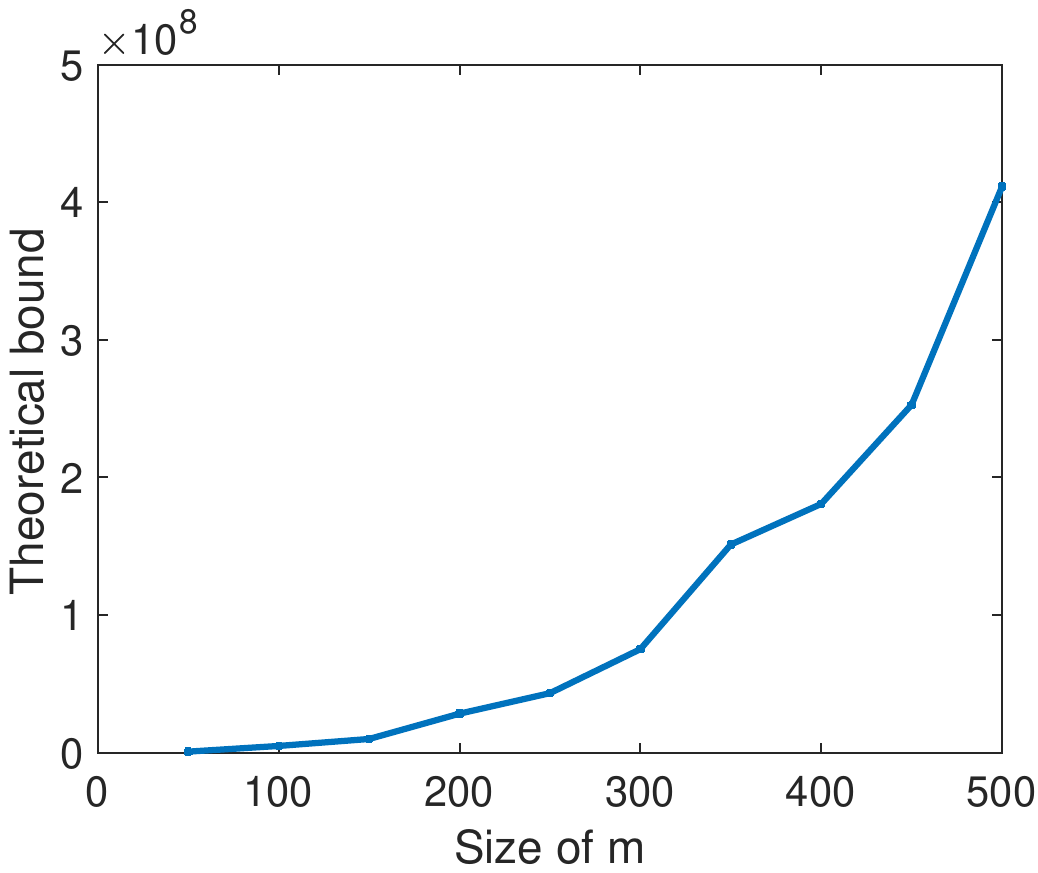}
\includegraphics[scale=0.3]{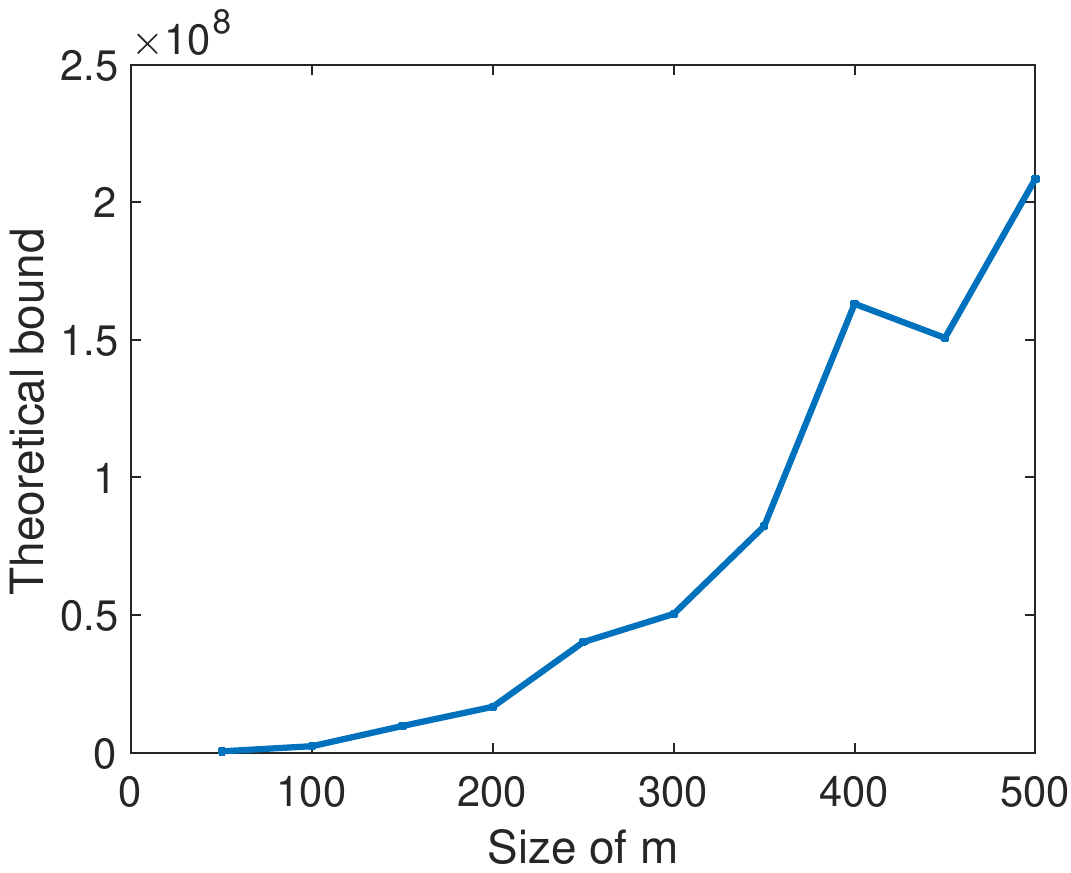}
  \caption{Left: folded Gaussian distribution. Right: uniform distribution}
  \label{fig:complexity}
\end{figure}
\revision{\subsection{Empirical Evaluation of Complexity Bounds}
\label{sec:complexity-bound-plots}
In this section, in order to study the dependence on problem dimensions of our results, we provide simulation results for the rate of linear convergence in our complexity bound in \cref{prop:iteration-complexity}, namely 
\[ \xi(A, p, q) = \max \Big\{ \tfrac{\| d^0_{\ast} \|_{\infty}}{\|
d^0_{\ast} \|_{- \infty}}, \tfrac{1}{\| d^0_{\ast} \|_{- \infty} \|
d^1_{\ast} \|_{- \infty}}, \| d^0_{\ast} \|_{\infty} \| d^1_{\ast}
\|_{\infty} \Big\}^2\cdot\tfrac{\min\{\max_j q_j, \max_i p_i\}}{\lambda_{-2}(\mathcal{L})}  . \]
We generate random matrices $A \in \mathbb{R}^{n \times m}$ element-wise from two distributions:
folded Gaussian $| \mathcal{N} (0, 1) |$ and uniform distribution $\mathcal{U}
[0, 1]$. The sparsity of $A$ is set to be around 80\%, and we only include cases when the Sinkhorn iterations converge. The target marginals $p$ and $q$ are generated from
uniform distribution $\mathcal{U} [0, 1]$ and scaled to ensure their sum are
the same.
We run simulation 100 times on different sizes of matrices $(n, m)$ where $n =
50 k + 50, k = 0, \ldots, 5$ and $m = 2 n$. We report the median of $\xi (A, p,
q)$ as $n$ increases in \cref{fig:complexity}. As we see, the rate approximately depends quadratically on the problem dimension. It remains an interesting question whether the rate in our bound can be improved to have logarithmic dependence on $n$.
}

%% file: proofs.tex
\section{Proofs}
\label{app:proofs}
\subsection{Proof of \cref{prop:mle-scaling}}
By construction, any normalized $s_j$ solving the optimality conditions in \eqref{eq:optimality} satisfy the matrix equations in \eqref{eq:equation-system}. 
It remains to show that a solution to \eqref{eq:equation-system} uniquely determines a solution to \eqref{eq:optimality}. Suppose two positive diagonal matrices $D^{0}$ and $D^{1}$ satisfy \cref{eq:scaled-matrix,eq:marginal}, i.e., $D^{1}AD^{0}\mathbf{1}_m =p$. Since the $i$-th row of $A$ is the indicator of the choice set $S_i$, we must have
\begin{align*}
D_{j}^{0} & =d_{j}\\
D_{i}^{1} & =\frac{p_i}{\sum_{k\in S_{i}}d_{k}},
\end{align*}
 for some positive $d_j$'s. The condition in \eqref{eq:bridge} that $\hat{A}^{T}\mathbf{1}_m = q$ then implies 
\begin{align*}
   q_j &= (D^0A^TD^1\mathbf{1}_m)_j\\
   &= \sum_{i\in[n]\mid j\in S_{i}}p_i\frac{d_{j}}{\sum_{k\in S_{i}}d_{k}},
\end{align*}
so that $d_j$'s, the diagonal entries of $D^0$, satisfy the optimality condition \eqref{eq:optimality} of the maximum likelihood estimation problem of Luce choice models. $\hfill \square$ 
 % However, if we impose the constraint that $\sum_{j}D_{j}^{0}=\sum_{j}s_{j}=1$,
% $D^{0}$ and $D^{1}$ are uniquely defined to be the above. 

 \subsection{Proof of \cref{thm:necessary-and-sufficient}}
First, we can verify that \cref{ass:weak-connected} is equivalent to the uniqueness condition \cref{ass:matrix-uniqueness} of the matrix balancing problem, namely the participation matrix $A$ is not permutation equivalent to a block-diagonal matrix. 

Now we prove that for the triplet $(A,p,q)$ constructed from the choice dataset, \cref{ass:strong-connected} on the choice dataset is equivalent to \cref{ass:matrix-existence} combined with \cref{ass:matrix-uniqueness} when $p, q$ are strictly positive. Consider an arbitrary pair of sets of indices $N \subsetneq [n]$ and $M \subsetneq [m]$ such that $A_{ij}=0$ for $i\notin N$ and $j\in M$. In the choice problem this condition implies that items in $M$ only appear in choice sets index by $N$. Then \cref{ass:strong-connected} implies that there is at least one item $k\notin M$ that is chosen over some item $j \in M$, which means items in $M$ are not always selected in observations with choice sets indexed by $N$, i.e., $\sum_{i\in N}p_i > \sum_{j\in M}q_j$. Moreover, \cref{ass:strong-connected} implies \cref{ass:weak-connected}, which is equivalent to \cref{ass:matrix-uniqueness}, so we have shown that \cref{ass:strong-connected} implies \cref{ass:matrix-existence} and \cref{ass:matrix-uniqueness}.

The converse direction is slightly less obvious. Suppose the $(A, p^0, p^1)$ constructed from a choice dataset satisfies Assumptions \ref{ass:matrix-existence} and \ref{ass:matrix-uniqueness}. In the choice dataset, consider an arbitrary partition of $[m]$ into $M$ and $M^C$. There are two cases to consider. First, suppose items in $M$ do not appear in all choice sets, i.e., there exists $N \subsetneq [n]$ such that $A_{ij}=0$ for $i\notin N$ and $j\in M$. Then \cref{ass:matrix-existence} of the matrix balancing problem implies $\sum_{i\in N}p^0_i > \sum_{j\in M}p^1_j$, i.e., some item $k \in M^C$ is selected over some item $j\in M$, which is what \cref{ass:strong-connected} requires. Second, suppose that every choice set contains least one item in $M$. Then the fact that $q>0$ implies that $q_k>0$ for any $k\in M^C$, i.e., at least some item $k\in M^C$ is selected over some item $j\in M$, which is again the condition in \cref{ass:strong-connected}. We have thus shown that there is some $k\in M^C$ that is chosen over some $j\in M$, and similarly vice versa, as required by \cref{ass:strong-connected}.
%To reiterate, the condition $\sum_{i\in N}p_i \geq \sum_{j\in M}q_j$ ensures a maximizer exists, although not necessarily in the interior of the simplex. When equality holds only if $A$ is in block diagonal form, maximizers are in the interior. Otherwise, maximizers are on the boundary, i.e., $s_j=0$ for some $j$.
$\hfill \square$

 \subsection{Proof of \cref{thm:global-convergence}}
 \begin{algorithm}[tb]
\caption{Normalized Sinkhorn's Algorithm}
   \label{alg:scaling-normalized}
\begin{algorithmic}
   \STATE {\bfseries Input:}  $A, p, q,\epsilon_{\text{tol}}$.
   \STATE {\bfseries initialize} $d^{0}\in\mathbb{R}_{++}^{m}$
   \REPEAT
   \STATE $d^{1} \leftarrow  p/( A d^0)$ 
   \STATE 
  normalization  $(d^0,d^1) \leftarrow (d^0/c,c d^1),c>0$
   \STATE $d^{0}\leftarrow  q/({A}^{T} d^{1})$
   \STATE 
  normalization  $(d^0,d^1) \leftarrow (d^0/c,c d^1),c>0$ 
   \STATE 
$\epsilon\leftarrow$  maximal update in $(d^{0},d^1)$ or distance between $D^1Ad^0$ and $p$
\UNTIL{$\epsilon<\epsilon_{\text{tol}}$}
\end{algorithmic}
\end{algorithm}
Suppose first that Sinkhorn's algorithm is normalized at each iteration as in \cref{alg:scaling-normalized} with $c$ defined in \eqref{eq:normalization}. We will first prove the linear convergence result for this normalized version, then argue that it applies to the standard Sinkhorn's algorithm \cref{alg:scaling}.

 \textbf{Outline of Proof.} \textbf{I.} We show that $g(u,v)$ is coercive on $\mathbf{1}_{m+n}^\perp$ under \cref{ass:matrix-existence}, which then guarantees that normalized iterates $(u^{(t)},v^{(t)})$ stay bounded. \textbf{II.} The next step is to study the Hessian of $g(u,v)$ on $\mathbf{1}_{m+n}^\perp$. The key observation is that the Hessian $\nabla^2 g(0,0)$ is the Laplacian matrix $\mathcal{L}$. The boundedness of normalized iterates $(u^{(t)},v^{(t)})$ then allows us to bound $\nabla^2g(u^{(t)},v^{(t)})$. \textbf{III.} Given the (subspace) strong convexity and smoothness of $g(u,v)$, the linear convergence rate of Sinkhorn's algorithm can then be obtained by applying the result of \citet{beck2013convergence} on the linear convergence of alternating minimization methods, and invoking the invariance of $g(u,v)$ under constant translations.

\textbf{I.} Recall the reparameterized potential function, for which Sinkhorn's
algorithm is the alternating minimization algorithm:
\begin{align*}
\min_{u\in\mathbb{R}^{m},v\in\mathbb{R}^{n}}g(u,v) & :=\sum_{ij}A_{ij}e^{-v_{i}+u_{j}}+\sum_{i=1}^{n}p_{i}v_{i}-\sum_{j=1}^{m}q_{j}u_{j}.
\end{align*}
 First, we show that $g(u,v)$ is coercive  on $\mathbf{1}_{m+n}^\perp$, i.e., $g(u,v)\rightarrow +\infty$ if $(u,v)\in \mathbf{1}_{m+n}^\perp$ and $\|(u,v)\|\rightarrow \infty$, whenever \cref{ass:matrix-existence} and \cref{ass:matrix-uniqueness} are satisfied. There are several cases.
\begin{enumerate}
     \item If $u_j \rightarrow +\infty$ for some $j$ but all $v_i$ stay bounded from above, then since $A$ does not contain zero columns or rows, the term $\sum_{ij}A_{ij}e^{-v_{i}+u_{j}}$ dominates, so $g(u,v)\rightarrow +\infty$.

     \item  Similarly, if some $v_i \rightarrow -\infty$ but all $u_j$ stay bounded from below, the term $\sum_{ij}A_{ij}e^{-v_{i}+u_{j}}$ dominates, so $g(u,v)\rightarrow +\infty$.

     \item If $u_j \rightarrow -\infty$ for some $j$ but all $v_i$ stay bounded from below, the term $-\sum_{j=1}^{m}q_{j}u_{j}$ dominates, so $g(u,v)\rightarrow +\infty$.

     \item Similarly, if $v_i \rightarrow+\infty$ for some $i$ but all $u_j$ stay bounded from above, the term $\sum_{i=1}^{n}p_{i}v_{i}$ dominates, so $g(u,v)\rightarrow +\infty$.

     \item Suppose now $v_i, u_j\rightarrow +\infty$ for some $i,j$. Then the subsets 
\begin{align*}
    I:=\{i: v_i \rightarrow +\infty\},\quad J:=\{j: u_j \rightarrow +\infty\}.
\end{align*}
    are both non-empty. Now either the exponential terms or the linear terms could dominate. If $u_j-v_i \rightarrow +\infty$ for some $i\in I, j\in J$, then one of the exponential terms dominates and $g(u,v)\rightarrow +\infty$. Otherwise, $-\min_{i\in I} v_i+\max_{j\in J} u_j$ stay bounded above, and the sum $\sum_{i\in I, j\in J}A_{ij}e^{-v_{i}+u_{j}}$ over $I,J$ stays bounded. There are now two sub-cases. 
     \begin{itemize}
         \item If there exists $i \notin I,j\in J$ such that $A_{ij} >0$, then  $A_{ij}e^{-v_{i}+u_{j}} \rightarrow +\infty$, since $v_i$ for $i\notin I$ is bounded from above.
         \item If $A_{i,j}=0$ for all $j\in J,i\notin I$. Suppose first $I=[n]$, i.e., all $v_i \rightarrow +\infty$. Now since we require that $(u,v) \in 
 \mathbf{1}_{m+n}^\perp$, there must exist some $j$ such that $u_j \rightarrow -\infty$, i.e., $J\subsetneq [m]$. Then we have $\sum_{i\in I} p_i=\sum_{i=1}^n p_i=\sum_{j=1}^m q_j > \sum_{j\in J} q_j$. 
 Thus, the linear terms dominate:  
  \begin{align*}
\sum_{i=1}^{n}p_{i}v_{i}-\sum_{j=1}^{m}q_{j}u_{j} \geq  \sum_{i=1}^{n}p_{i}v_{i}-\sum_{j\in J}q_{j}u_{j} \rightarrow +\infty.
     \end{align*}
         Now suppose $I \subsetneq [n]$. If $J=[m]$, i.e., $u_j \rightarrow+\infty$ for all $j$, then the requirement that $(u,v) \in 
 \mathbf{1}_{m+n}^\perp$ again guarantees that for some $i\notin I$, $v_i \rightarrow -\infty$. Since $A$ does not contain zero rows, the sum $\sum_{j=1}^m A_{ij}e^{-v_{i}+u_{j}} \rightarrow +\infty$.
         Lastly, if both  $I\subsetneq [n]$ $J\subsetneq [m]$,  then \cref{ass:matrix-existence} applies to $I,J$, and guarantees that $\sum_{i\in I} p_i > \sum_{j\in J} q_j$, so that
          \begin{align*}
\sum_{i\in I}p_{i}v_{i}-\sum_{j\in J}q_{j}u_{j} \rightarrow +\infty.
     \end{align*}
     If, in addition, for all $i\notin I$, $v_i$ is bounded below, then the linear terms dominate
     \begin{align*}
\sum_{i=1}^{n}p_{i}v_{i}-\sum_{j=1}^{m}q_{j}u_{j} \rightarrow +\infty.
     \end{align*}
     And if for some $i\notin I$, $v_i \rightarrow -\infty$, then since $A$ does not contain zero rows, the sum $\sum_{j=1}^m A_{ij}e^{-v_{i}+u_{j}} \rightarrow +\infty$.  
     \end{itemize}
     \item The case $v_i, u_j\rightarrow -\infty$ for some $i,j$ is symmetric to the previous case and we omit the detailed reasoning.   
 \end{enumerate}

 Note that coercivity does not hold if we do not restrict to $\mathbf{1}_{m+n}^\perp$, since $g(u,v)=g(u+c,v+c)$ for any constant $c\in\mathbb{R}$. This is one reason why \emph{normalization} is essential in our analysis.
 
 Coercivity on $\mathbf{1}_{m+n}^\perp$ guarantees that there exist finite optimal solutions $(u,v)$. More importantly, it implies that the sub-level sets
 \begin{align*}
   S_g^\perp(\alpha):  \{(u,v)\in \mathbf{1}_{m+n}^\perp:g(u,v)\leq \alpha \}
 \end{align*}
 are \emph{bounded}. In particular, this property holds for $\alpha_0=g(u^{(0)},v^{(0)})$,\footnote{In practice, we can for example take $(u^{(0)},v^{(0)})=0$, so that $\alpha_0=\sum_{i,j}A_{ij}$.} so that there exists $B<\infty$ such that $\|(u,v)\|_{\infty}\leq B$ whenever $(u,v)\in \mathbf{1}_{m+n}^\perp$ and  $g(u,v)\leq g(u^{(0)},v^{(0)})$. Since Sinkhorn's algorithm is the alternating minimization
of $g(u,v)$,
\begin{align*}
    g(u^{(t+1)},v^{(t+1)})\leq g(u^{(t)},v^{(t)})\leq g(u^{(0)},v^{(0)})
\end{align*}
for all $t$. It then follows that $(u^{(t)},v^{(t)}) \in S_g^\perp(\alpha_0)$, so that $\|(u^{(t)},v^{(t)})\|_{\infty}\leq B$ for
all $t>0$. In summary, we have shown that normalized Sinkhorn iterations stay bounded by $B$, which will be important for lower bounding the Hessian of $g$. 

\textbf{II.} Next, we show that $g$ is strongly convex when restricted to $ \mathbf{1}_{m+n}^\perp$. The gradient of $g(u,v)$ is given by 
\begin{align*}
\partial_{u_{j}}g(u,v) & =\sum_{i}A_{ij}e^{-v_{i}+u_{j}}-q_{j}\\
\partial_{v_{i}}g(u,v) & =-\sum_{j}A_{ij}e^{-v_{i}+u_{j}}+p_{i},
\end{align*}
and the Hessian is given by 
\begin{align*}
\nabla^{2}g(u,v)=\begin{bmatrix}\mathcal{D}(\sum_{j}A_{ij}e^{u_{j}-v_{i}}) & -\hat{A}\\
-\hat{A}^{T} & \mathcal{D}(\sum_{i}A_{ij}e^{u_{j}-v_{i}})
\end{bmatrix} & =\begin{bmatrix}\mathcal{D}(\hat{A}\mathbf{1}_{m}) & -\hat{A}\\
-\hat{A}^{T} & \mathcal{D}(\hat{A}^{T}\mathbf{1}_{n})
\end{bmatrix}
\end{align*}
 where 
\begin{align*}
\hat{A}_{ij}(u,v) & =A_{ij}e^{-v_{i}+u_{j}}.
\end{align*}
Note that $\mathbf{1}_{m+n}$ is in the null space of the Hessian
at any $(u,v)$, since
\begin{align*}
\begin{bmatrix}\mathcal{D}(\hat{A}\mathbf{1}_{m}) & -\hat{A}\\
-\hat{A}^{T} & \mathcal{D}(\hat{A}^{T}\mathbf{1}_{n})
\end{bmatrix}\begin{bmatrix}\mathbf{1}_{n}\\
\mathbf{1}_{m}
\end{bmatrix} & =\begin{bmatrix}\hat{A}\mathbf{1}_{m}-\hat{A}\mathbf{1}_{m}\\
-\hat{A}^{T}\mathbf{1}_{n}+\hat{A}^{T}\mathbf{1}_{n}
\end{bmatrix}=\begin{bmatrix}\mathbf{0}\\
\mathbf{0}
\end{bmatrix}.
\end{align*}
As a result, if we restrict to the subspace $(u,v)\perp\begin{bmatrix}\mathbf{1}_{m}\\
\mathbf{1}_{n}
\end{bmatrix}$, which is achieved with the appropriate normalization in Sinkhorn's
algorithm, the Hessian of the potential function $g(u,v)$ is lower
bounded by $\lambda_{-2}(\nabla^{2}g(u,v))\cdot I_{m+n}$. This is the second reason why normalization is important for the linear convergence analysis of Sinkhorn's algorithm.

The next key observation is that the Hessian $\nabla^{2}g(u,v)$ is precisely the \emph{Laplacian}
matrix of the bipartite graph defined by $\hat{A}$: 
\begin{align*}
\mathcal{L}(\hat{A}):=\begin{bmatrix}\mathcal{D}(A\mathbf{1}_{m}) & -\hat{A}\\
-\hat{A}^{T} & \mathcal{D}(\hat{A}^{T}\mathbf{1}_{n})
\end{bmatrix} & =\mathcal{D}(\begin{bmatrix}0 & \hat{A}\\
\hat{A}^{T} & 0
\end{bmatrix}\begin{bmatrix}\mathbf{1}_{n}\\
\mathbf{1}_{m}
\end{bmatrix})-\begin{bmatrix}0 & \hat{A}\\
\hat{A}^{T} & 0
\end{bmatrix},
\end{align*}
where we recognize 
\begin{align*}
\begin{bmatrix}0 & \hat{A}\\
\hat{A}^{T} & 0
\end{bmatrix}
\end{align*}
as the adjacency matrix of the bipartite graph defined by $\hat{A}$, and $\mathcal{L}$ defined in \eqref{eq:Laplacians} is precisely $\nabla^2g(0,0)$.

Our next step is then to connect $\mathcal{L}(\hat{A})$ to the Laplacian $\mathcal{L}$ at the origin, so that we can lower bound $\lambda_{-2}(\nabla^{2}g(u,v))$. 
% We use the
% well-known fact that the bipartite graph Laplacian $\nabla^{2}g(u,v)$
% is \emph{similar} to the following ``signless'' Laplacian, so they share the
% same spectrum:
% \begin{align*}
% \mathcal{L}'(\hat{A}):= & \mathcal{D}(\begin{bmatrix}0 & \hat{A}\\
% \hat{A}^{T} & 0
% \end{bmatrix}\begin{bmatrix}\mathbf{1}_{n}\\
% \mathbf{1}_{m}
% \end{bmatrix})+\begin{bmatrix}0 & \hat{A}\\
% \hat{A}^{T} & 0
% \end{bmatrix},
% \end{align*}
%  so that 
% \begin{align*}
% \lambda_{-2}(\nabla^{2}g(u,v)) & =\lambda_{-2}(\mathcal{L}'(\hat{A})).
% \end{align*}
We claim that for all $(u,v) \in \mathbf{1}_{m+n}^\perp$ with $\|(u,v)\|_\infty\leq B$, 
\begin{align*}
\mathcal{L}(\hat{A})=\mathcal{D}(\begin{bmatrix}0 & \hat{A}\\
\hat{A}^{T} & 0
\end{bmatrix}\begin{bmatrix}\mathbf{1}_{n}\\
\mathbf{1}_{m}
\end{bmatrix})-\begin{bmatrix}0 & \hat{A}\\
\hat{A}^{T} & 0
\end{bmatrix} & \succeq e^{-2B}(\mathcal{D}(\begin{bmatrix}0 & A\\
A^{T} & 0
\end{bmatrix}\begin{bmatrix}\mathbf{1}_{n}\\
\mathbf{1}_{m}
\end{bmatrix})-\begin{bmatrix}0 & A\\
A^{T} & 0
\end{bmatrix})=e^{-2B}\mathcal{L}(A).
\end{align*}
 Consider first the off-diagonal blocks. We have 
\begin{align*}
\hat{A}_{ij} & =A_{ij}e^{-v_{i}+u_{j}}\geq e^{-2B}A_{ij}.
\end{align*}
 Similarly, for the diagonal blocks, we have 
\begin{align*}
\mathcal{D}(\hat{A}\mathbf{1}_{m})_{i} & =\sum_{j}A_{ij}e^{u_{j}-v_{i}}\geq e^{-2B}\sum_{j}A_{ij}\\
\mathcal{D}(\hat{A}^{T}\mathbf{1}_{n})_{j} & =\sum_{i}A_{ij}e^{u_{j}-v_{i}}\geq e^{-2B}\sum_{i}A_{ij}.
\end{align*}
 The above inequalities imply that the diagonal entries of the following difference 
\begin{align*}
\mathcal{L}(\hat{A})-e^{-2B}\mathcal{L}(A)
\end{align*}
are non-negative and its off-diagonal terms are non-positive. Moreover, since both $\mathcal{L}(\hat{A})$ and
$\mathcal{L}(A)$ have $\mathbf{1}_{n+m}$ in their
null spaces, so does the difference above. Gershgorin circle theorem
then guarantees that the eigenvalues of $\mathcal{L}(\hat{A})-e^{-2B}\mathcal{L}(A)$
are all non-negative. 
% Finally, since the signless Laplacian $\mathcal{L}'(A)$
% shares the same spectrum as the Laplacian 
% \begin{align*}
% \mathcal{L} & =\mathcal{D}(\begin{bmatrix}0 & A\\
% A^{T} & 0
% \end{bmatrix}\begin{bmatrix}\mathbf{1}_{n}\\
% \mathbf{1}_{m}
% \end{bmatrix})-\begin{bmatrix}0 & A\\
% A^{T} & 0
% \end{bmatrix},
% \end{align*}
We can conclude that 
\begin{align*}
\lambda_{-2}(\nabla^{2}g(u,v)) & \geq e^{-2B}\lambda_{-2}(\mathcal{L}).
\end{align*}
It then follows that $g(u,v)$ is $e^{-2B}\lambda_{-2}(\mathcal{L})$-strongly
convex on the subspace $(u,v)\perp\mathbf{1}_{n+m},\|(u,v)\|_{\infty}\leq B$. Recall that coercivity of $g(u,v)$ on $\mathbf{1}_{m+n}^\perp$ precisely guarantees $\|(u^{(t)},v^{(t)})\|_{\infty}\leq B$ for all $(u^{(t)},v^{(t)})$ during the iterations of normalized Sinkhorn's algorithm.

Next, we compute the smoothness constants $L_{0},L_{1}$ of $g(u,v)$
when restricted to one of the two blocks of variables. Recall that the gradient of $g$ is given by 
\begin{align*}
\partial_{u_{j}}g(u,v) & =\sum_{i}A_{ij}e^{-v_{i}+u_{j}}-q_{j}\\
\partial_{v_{i}}g(u,v) & =-\sum_{j}A_{ij}e^{-v_{i}+u_{j}}+p_{i},
\end{align*}
so that for any $i,j$,
\begin{align*}
|\partial_{u_{j}}^{2}g(u,v)| & \leq e^{2B}\sum_{i}A_{ij}\\
|\partial_{v_{i}}^{2}g(u,v)| & \leq e^{2B}\sum_{j}A_{ij}.
\end{align*}
 It then follows that the Lipschitz constants of the two blocks are
given by 
\begin{align*}
L_{0} & =e^{2B}\max_{j}\sum_{i}A_{ij}=e^{2B}l_{0}\\
L_{1} & =e^{2B}\max_{i}\sum_{j}A_{ij}=e^{2B}l_{1}.
\end{align*}
\textbf{III.} %Recall that the function $g(u,v)$, its gradient, and its Hessian are all invariant under the transformation $(u,v)\rightarrow (u-\log c, v-\log c)$ for any $c>0$. It follows that the strong convexity constant of $g(u,v)$ restricted to the subspace $\mathbf{1}_{m+n}$ in fact applies to $g(u,v)$ on $\mathbb{R}^{m+n}$. 
Given the strong convexity of $g(u,v)$ on the bounded subspace of $\mathbf{1}_{m+n}$ containing the normalized Sinkhorn iterates, we can then apply Theorem 5.2 of \citet{beck2013convergence} to conclude that for any $t\geq 0$, normalized Sinkhorn iterates $(u^{(t+1)},v^{(t+1)})$ resulting from \cref{alg:scaling-normalized} with $c$ defined in \eqref{eq:normalization} satisfies
\begin{align*}
g(u^{(t+1)},v^{(t+1)})-g^{\ast} & \leq(1-e^{-4B}\frac{\lambda_{-2}(\mathcal{L})}{\min\{l_{0},l_{1}\}})\cdot\left(g(u^{(t)},v^{(t)})-g^{\ast}\right).
\end{align*}
Since the objective $g(u,v)$ is \emph{invariant} under normalization, one can verify that the sequence of objective values $g(u^{(t)},v^{(t)})$ remains unchanged regardless of the specific normalization constant $c$ used. The linear convergence result therefore carries over to the original Sinkhorn's algorithm without normalization, i.e., $c=1$, as well as to any other normalized versions of Sinkhorn's algorithm, even when $c$ changes with iteration number $t$.

Lastly, note that we do not need to compute the smoothness constant $L$ of
the entire function $g(u,v)$, since convergence results on alternating
minimization algorithms only require the \emph{minimum} of smoothness
constants of the function restricted to each block, which is upper
bounded by $L$. Nevertheless, for $g(u,v)$, we compute its $L$
explicitly for completeness. Note that by a similar reasoning used
to lower bound $\mathcal{L}(\hat{A})$, we have
\begin{align*}
\mathcal{L}(\hat{A}) & \preceq e^{2B}\mathcal{L}(A).
\end{align*}
Gershgorin circle theorem then bounds the maximal eigenvalue
\begin{align*}
    \lambda_{1}(\mathcal{L}(A))=\lambda_{1}(\mathcal{L})\leq2\max\{\max_{i}\sum_{j}A_{ij},\max_{j}\sum_{i}A_{ij}\}=L,
\end{align*}
which is always strictly larger than $l=\min\{\max_{i}\sum_{j}A_{ij},\max_{j}\sum_{i}A_{ij}\}$.
  $\hfill \square$

  \input{skbnd}

  \subsection{Proof of \cref{thm:lower-bound}}
  We will show that there exists some row $i$ of the scaled matrix $A^{(t)}$ such that $|\sum_{j=1}^m A^{(t)}_{ij}-p_i|$ vanishes at a rate bounded below by $\Omega(1/t)$. Because \cref{ass:matrix-existence} is satisfied but \cref{ass:matrix-weak-existence} fails, no finite solution to the
matrix balancing problem exists, but Sinkhorn's algorithm converges
to $\hat{A}$. Without loss of generality, we may assume column scalings
do not vanish, since we can always multiply and divide the scalings
without altering the scaled matrix. Lemma 1 of \citet{pukelsheim2014biproportional} states that
there exists $I\subseteq[n]$ and $J\subseteq[m]$ such that, after
necessary permutation, $A$ can be written in block form as 
\begin{align*}
A=\begin{bmatrix}A_{IJ} & 0\\
A_{I^{C}J} & A_{I^{C}J^{C}}
\end{bmatrix},
\end{align*}
 and 
\begin{align*}
\sum_{j\in J^{C}}q_{j} & =\sum_{i\in I^{C}}p_{i}\\
\sum_{j\in J}q_{j} & =\sum_{i\in I}p_{i}.
\end{align*}
Moreover, $I$ is the set of rows with row scalings not vanishing,
and $J$ is the set of columns with column scalings not diverging:
\begin{align*}
I=\{i:d_{i}^{1(t)}\not\rightarrow0\},\quad J=\{j:d_{j}^{0(t)}\not\rightarrow\infty\}.
\end{align*}
Thus there exists a constant $c$ such that $\sum_{k\in J}d_{k}^{0(t)}\leq c$
for (a subsequence of) $t=1,2,\dots$. 

For any non-negative matrix inheriting the zeros of $A$ and having
marginals $p,q$, the lower left block indexed by $I^{C}J$ must be
identically 0. This block is referred to as a fading block. Thus the
limiting matrix $\hat{A}$ of Sinkhorn's algorithm, which solves the
primal KL divergence minimization problem, must have the form 
\begin{align*}
\hat{A}=\begin{bmatrix}\hat{A}_{IJ} & 0\\
0 & \hat{A}_{I^{C}J^{C}}
\end{bmatrix},
\end{align*}
and all entries in the lower left block in $A^{(t)}$ converge to
$0$ as the number of iterations increases.

Without loss of generality, we may assume that $\min_{ij:A_{ij}>0}A_{ij}\geq1$,
since we can rescale $A$ by a fixed constant without altering the
convergence rate. Define $\overline{A}:=\max_{ij}A_{ij}$. We claim
that $\sum_{j\in J^{C}}d_{j}^{0(t)}\leq c\overline{A}t$ for all $t$.
We can initialize with $d^{0(1)}$ that satisfies $\sum_{j\in J^{C}}d_{j}^{0(1)}\leq\overline{A}c$,
and prove the general case by induction. Suppose $\sum_{j\in J^{C}}d_{j}^{0(t)}\leq c\overline{A}t$
and we want to show $\sum_{j\in J^{C}}d_{j}^{0(t+1)}\leq c\overline{A}(t+1)$.
We have 
\begin{align*}
\sum_{j\in J^{C}}d_{j}^{0(t+1)}=\sum_{j\in J^{C}}q_{j}/\sum_{i}A_{ij}\frac{p_{i}}{\sum_{k}A_{ik}d_{k}^{0(t)}} & =\sum_{j\in J^{C}}q_{j}/\sum_{i\in I^{C}}A_{ij}\frac{p_{i}}{\sum_{k}A_{ik}d_{k}^{0(t)}},
\end{align*}
 where the last equality follows from the block structure of $A$.
Note that since $A$ cannot have zero columns, for any $j\in J^{C}$,
there is at least one $i\in I^{C}$ such that $A_{ij}>0$. 

Next, for any $i\in I^{C}$, we have 
\begin{align*}
\sum_{k}A_{ik}d_{k}^{0(t)} & =\sum_{k\in J}A_{ik}d_{k}^{0(t)}+\sum_{k\in J^{C}}A_{ik}d_{k}^{0(t)}\\
 & \leq\overline{A}\cdot\sum_{k\in J}d_{k}^{0(t)}+\overline{A}\cdot\sum_{k\in J^{C}}d_{k}^{0(t)}\\
 & \leq c\overline{A}+c\overline{A}t=c\overline{A}(t+1),
\end{align*}
 where the last inequality follows from the induction assumption and
the condition that the column scalings in $J$ stay bounded. Using
the inequality above and $\min_{ij:A_{ij}>0}A_{ij}\geq1$, we have
\begin{align*}
\sum_{j\in J^{C}}d_{j}^{0(t+1)} & =\sum_{j\in J^{C}}q_{j}/\sum_{i\in I^{C}}A_{ij}\frac{p_{i}}{\sum_{k}A_{ik}d_{k}^{0(t)}}\\
 & \leq\sum_{j\in J^{C}}q_{j}/(\sum_{i\in I^{C}}\frac{p_{i}}{c\overline{A}(t+1)})\\
 & =c\overline{A}(t+1)\cdot(\sum_{j\in J^{C}}q_{j}/\sum_{i\in I^{C}}p_{i})\\
 & =c\overline{A}(t+1).
\end{align*}
 where the last equality precisely follows from the condition $\sum_{j\in J^{C}}q_{j}=\sum_{i\in I^{C}}p_{i}$
in \cref{ass:matrix-weak-existence}. We have therefore proved that $\sum_{j\in J^{C}}d_{j}^{0(t)}\leq c\overline{A}t$
for all $t$. The pigeonhole principle implies that there exists a
$j\in J^{C}$ and a subsequence $t_{1},t_{2},\cdots\rightarrow\infty$
such that $d_{j}^{0(t_{l})}\leq\frac{1}{|J^{C}|}c\overline{A}t_{l}$
for all $l=1,2,\dots$. Since the $j$-th column of $\hat{A}$ cannot
all vanish, there exists $i\in I^{C}$ such that $d_{i}^{1(t_{l})}A_{ij}d_{j}^{0(t_{l})}\rightarrow\hat{A}_{ij}>0$,
so that, selecting a further subsequence if necessary, we have 
\begin{align*}
d_{i}^{1(t_{l})} & \geq\hat{A}_{ij}/(\frac{1}{|J^{C}|}c\overline{A}t_{l}).
\end{align*}
 In other words, the row scaling $d_{i}^{1(t)}$ for row $i\in I^{C}$
vanishes at a rate bounded below by $\Omega(1/t)$. Recall that by assumption column scalings
are non-vanishing, which implies $\sum_{j\in J}^m A^{(t)}_{ij}\rightarrow 0$ at a rate bounded below by $\Omega(1/t)$. The same reasoning in fact guarantees that $\sum_{j\in J}^m A^{(t)}_{ij}= \Omega(1/t)$ for all $i\in I^C$, i.e., the fading block vanishes at rate $\Omega(1/t)$, so that $D_{\text{KL}}(r^{(t)}\| p)=\Omega(1/t)$.
$\hfill \square$
 
\subsection{Proof of \cref{thm:convergence}}
 \textbf{Outline of Proof.} \textbf{I.} We define a novel sequence of data-dependent mappings induced
by Sinkhorn's iterations, $f^{(t)}(x):\mathbb{R}_{++}^{n}\rightarrow\mathbb{R}_{++}^{n}$,
that map $r^{(t)}/\sqrt{p}$ to $r^{(t+1)}/\sqrt{p}$, starting from
column normalized $A^{(t)}$. For all $t$, $\sqrt{p}$ is always
the fixed point of $f^{(t)}$, and the construction guarantees that
the residuals $r^{(t)}/\sqrt{p}-\sqrt{p}$ are always orthogonal to
the fixed point $\sqrt{p}$. \textbf{II.} The Jacobian $J^{(t)}$ at $\sqrt{p}$ is given by
\begin{align*}
J^{(t)}= & \mathcal{D}(\sqrt{p}/r^{(t)})A^{(t)}\mathcal{D}(1/q)A^{(t)T}\mathcal{D}(1/\sqrt{p})
\end{align*}
 which has unique maximal eigenvector $\sqrt{p}$ with eigenvalue
1. Moreover, $J^{(t)}\rightarrow J=\mathcal{D}(1/\sqrt{p})\hat{A}\mathcal{D}(1/q)\hat{A}^{T}\mathcal{D}(1/\sqrt{p})$
as $t\rightarrow\infty$, and $\tilde{A}=\mathcal{D}(1/\sqrt{p})\cdot\hat{A}\cdot\mathcal{D}(1/\sqrt{q})$ also has unique maximal
eigenvector $\sqrt{p}$ with eigenvalue 1. \textbf{III.} Thus we have a sequence of mappings $f^{(t)}$ with
fixed points $\sqrt{p}$, whose Jacobians $J^{(t)}$ converge to $J$
with maximal eigenvector $\sqrt{p}$. The orthogonality property $(r^{(t)}/\sqrt{p}-\sqrt{p})^{T}\sqrt{p}=0$
and uniform boundedness of second derivatives of $f^{(t)}$ near $\sqrt{p}$
then allows us to conclude that the asymptotic linear convergence rate is exactly
equal to the subdominant eigenvalue of $J$.

\textbf{I.} Starting with column normalized matrix $A^{(t)}$ at the
$t$-th iteration, we first write down the general formula for the
row sums $r^{(t+1)}=A^{(t+1)}\mathbf{1}_{m}$ after one iteration
of Sinkhorn's algorithm:

\begin{align*}
r_{i}^{(t+1)} & =\\
\frac{A_{i1}^{(t)}\cdot p_{i}/r_{i}^{(t)}\cdot q_{1}}{A_{11}^{(t)}\frac{p_{1}}{r_{1}^{(t)}}+A_{21}^{(t)}\frac{p_{2}}{r_{2}^{(t)}}+\cdots+A_{n1}^{(t)}\frac{p_{n}}{r_{n}^{(t)}}}+\frac{A_{i2}^{(t)}\cdot p_{i}/r_{i}^{(t)}\cdot q_{2}}{A_{12}^{(t)}\frac{p_{1}}{r_{1}^{(t)}}+A_{22}^{(t)}\frac{p_{2}}{r_{2}^{(t)}}+\cdots+A_{n2}^{(t)}\frac{p_{n}}{r_{n}^{(t)}}}+\cdots+\frac{A_{im}^{(t)}\cdot p_{i}/r_{i}^{(t)}\cdot q_{m}}{A_{1m}^{(t)}\frac{p_{1}}{r_{1}^{(t)}}+A_{2m}^{(t)}\frac{p_{2}}{r_{2}^{(t)}}+\cdots+A_{nm}^{(t)}\frac{p_{n}}{r_{n}^{(t)}}}
\end{align*}
Although the update formula looks complicated, it can be interpreted
in the following manner. Given positive $r_{1}^{(0)}/p_{1},r_{2}^{(0)}/p_{2},\dots,r_{n}^{(0)}/p_{n}$,
we take the\emph{ convex} combination of their \emph{inverses}:
\begin{align*}
c_{j}^{(t)}/q_{j} & :=A_{1j}^{(t)}/q_{j}\cdot\frac{p_{1}}{r_{1}^{(t)}}+A_{2j}^{(t)}/q_{j}\cdot\frac{p_{2}}{r_{2}^{(t)}}+\cdots+A_{nj}^{(t)}/q_{j}\cdot\frac{p_{n}}{r_{n}^{(t)}},
\end{align*}
 where $A^{(t)}$ is assumed to be column normalized, hence have column
sums equal to $q$. After we have formed $c_{1}^{(t)}/q_{1},\dots,c_{m}^{(t)}/q_{m}$,
we simply repeat the process by taking the \emph{convex} combination
of their inverses: 
\begin{align*}
r_{i}^{(1)}/p_{i}:= & \frac{A_{i1}^{(t)}}{r_{i}^{(t)}}\cdot\frac{q_{1}}{c_{1}^{(t)}}+\frac{A_{i2}^{(t)}}{r_{i}^{(t)}}\cdot\frac{q_{2}}{c_{2}^{(t)}}+\cdots+\frac{A_{im}^{(t)}}{r_{i}^{(t)}}\cdot\frac{q_{m}}{c_{m}^{(t)}},
\end{align*}
 where by definition of $r_{i}^{(t)}$ we always have $\sum_{j}A_{ij}^{(t)}=r_{i}^{(t)}$.
We will work with this update formula to obtain a sequence of mappings. 

Recall that $/$ denotes entry-wise division whenever the quantities
are vectors, and similarly for $\sqrt{\cdot}$. Instead of $r^{(t)}$,
we will use $r^{(t)}/\sqrt{p}$ as the natural quantity to measure
the progress of convergence and show that $\|r^{(t)}/\sqrt{p}-\sqrt{p}\|_{2}\rightarrow0$
linearly with rate $\lambda$. The reason for using $r^{(t)}/\sqrt{p}$
is because the residual $r^{(t)}/\sqrt{p}-\sqrt{p}$ satisfies 
\begin{align*}
(r^{(t)}/\sqrt{p}-\sqrt{p})^{T}\sqrt{p} & =\sum_{i}(r_{i}^{(t)}-p_{i})=0
\end{align*}
by virtue of Sinkhorn's algorithm preserving the quantities $r^{(t)T}\mathbf{1}_{n}$
for all $t$, so that the residual is always \emph{orthogonal} to
$\sqrt{p}$. This orthogonality property is crucial in identifying the rate of convergence,
as we will show that $\sqrt{p}$ is also the unique maximal eigenvector
of the limiting Jacobian with eigenvalue 1. 

Now rewrite the update formula as 
\begin{align*}
(r^{(t+1)}/\sqrt{p})_{i} & =\\
\frac{A_{i1}^{(t)}\cdot\sqrt{p}_{i}/r_{i}^{(t)}\cdot q_{1}}{A_{11}^{(t)}p_{1}/r_{1}^{(t)}+A_{21}^{(t)}p_{2}/r_{2}^{(t)}+\cdots+A_{n1}^{(t)}p_{n}/r_{n}^{(t)}}+\cdots+\frac{A_{im}^{(t)}\cdot\sqrt{p}_{i}/r_{i}^{(t)}\cdot q_{m}}{A_{1m}^{(t)}p_{1}/r_{1}^{(t)}+A_{2m}^{(t)}p_{2}/r_{2}^{(t)}+\cdots+A_{nm}^{(t)}p_{n}/r_{n}^{(t)}} & =\\
\frac{A_{i1}^{(t)}\cdot\sqrt{p}_{i}/r_{i}^{(t)}\cdot q_{1}}{A_{11}^{(t)}\sqrt{p_{1}}/(r_{1}^{(t)}/\sqrt{p_{1}})+A_{21}^{(t)}\sqrt{p_{2}}/(r_{2}^{(t)}/\sqrt{p_{2}})+\cdots+A_{n1}^{(t)}\sqrt{p_{n}}/(r_{n}^{(t)}/\sqrt{p_{n}})}+\cdots & +\\
\frac{A_{im}^{(t)}\cdot\sqrt{p}_{i}/r_{i}^{(t)}\cdot q_{m}}{A_{1m}^{(t)}\sqrt{p_{1}}/(r_{1}^{(t)}/\sqrt{p_{1}})+A_{2m}^{(t)}\sqrt{p_{2}}/(r_{2}^{(t)}/\sqrt{p_{2}})+\cdots+A_{nm}^{(t)}\sqrt{p_{n}}/(r_{n}^{(t)}/\sqrt{p_{n}})}.
\end{align*}
And define the mapping $f^{(t)}(x):\mathbb{R}_{++}^{n}\rightarrow\mathbb{R}_{++}^{n}$,
given $A^{(t)}$ and $r^{(t)}$, as 
\begin{align*}
f_{i}^{(t)}(x) & :=\frac{A_{i1}^{(t)}\cdot\sqrt{p}_{i}/r_{i}^{(t)}\cdot q_{1}}{A_{11}^{(t)}\sqrt{p_{1}}/x_{1}+A_{21}^{(t)}\sqrt{p_{2}}/x_{2}+\cdots+A_{n1}^{(t)}\sqrt{p_{n}}/x_{n}}+\cdots\\
 & +\frac{A_{im}^{(t)}\cdot\sqrt{p}_{i}/r_{i}^{(t)}\cdot q_{m}}{A_{1m}^{(t)}\sqrt{p_{1}}/x_{1}+A_{2m}^{(t)}\sqrt{p_{2}}/x_{2}+\cdots+A_{nm}^{(t)}\sqrt{p_{n}}/x_{n}}.
\end{align*}
Written more compactly, 
\begin{align*}
f^{(t)}(x) & =\mathcal{D}(\sqrt{p}/r^{(t)})A^{(t)}\mathcal{D}(q)\mathcal{D}(A^{(t)T}(\sqrt{p}/x))^{-1}\mathbf{1}_{m}.
\end{align*}

Our first observation is that $f^{(t)}(x)$ has fixed point $\sqrt{p}$
for all $t$, since 
\begin{align*}
f^{(t)}(\sqrt{p}) & =\mathcal{D}(\sqrt{p}/r^{(t)})A^{(t)}\mathcal{D}(q)\mathcal{D}(A^{(t)T}(\sqrt{p}/\sqrt{p}))^{-1}\mathbf{1}_{m}\\
 & =\mathcal{D}(\sqrt{p}/r^{(t)})A^{(t)}\mathcal{D}(q)\mathcal{D}(A^{(t)T}\mathbf{1}_{n})^{-1}\mathbf{1}_{m}\\
 & =\mathcal{D}(\sqrt{p}/r^{(t)})A^{(t)}\mathcal{D}(q)\cdot\mathbf{1}_{m}/q\\
 & =\mathcal{D}(\sqrt{p}/r^{(t)})A^{(t)}\mathbf{1}_{m}=\mathcal{D}(\sqrt{p}/r^{(t)})r^{(t)}=\sqrt{p}
\end{align*}
 We also define the ``limiting'' mapping by replacing $r^{(t)}$ with
$p$ and $A^{(t)}$ with $\hat{A}$ in $f^{(t)}$:
\begin{align*}
f(x) & :=\mathcal{D}(1/\sqrt{p})\hat{A}\mathcal{D}(q)\mathcal{D}(\hat{A}^{T}(\sqrt{p}/x))^{-1}\mathbf{1}_{m}.
\end{align*}
 It is straightforward to verify that $f(\sqrt{p})=\sqrt{p}$ as well.
Sinkhorn's iteration can then be represented as the following sequence:
\begin{align*}
r^{(0)}/\sqrt{p}\rightarrow f^{(0)}(r^{(0)}/\sqrt{p}) & =r^{(1)}/\sqrt{p}\\
r^{(1)}/\sqrt{p}\rightarrow f^{(1)}(r^{(1)}/\sqrt{p}) & =r^{(2)}/\sqrt{p}\\
\cdots
\end{align*}
 We reiterate that $f^{(t)}(x)$ is data-dependent, as it uses $r^{(t)}$
and $A^{(t)}$. Although $f$ is the pointwise limit of $f^{(t)}$,
we do not need any uniform convergence, thanks to the orthogonality
condition and the fact that $\sqrt{p}$ is the fixed point and unique
maximal eigenvector for \emph{every }Jacobian at $\sqrt{p}$ by construction. We only need the uniform boundedness of $r^{(t)},A^{(t)}$, which is guaranteed by coercivity of the potential function. 

\textbf{II. }The partial derivatives of $f^{(t)}$ are given by 
\begin{align*}
\partial_{j}f_{i}^{(t)}(x) & =\frac{A_{i1}^{(t)}\cdot\sqrt{p}_{i}/r_{i}^{(t)}\cdot q_{1}}{(A_{11}^{(t)}\sqrt{p_{1}}/x_{1}+A_{21}^{(t)}\sqrt{p_{2}}/x_{2}+\cdots+A_{n1}^{(t)}\sqrt{p_{n}}/x_{n})^{2}}\cdot\frac{A_{j1}^{(t)}\sqrt{p_{j}}}{x_{j}^{2}}+\cdots\\
 & +\frac{A_{im}^{(t)}\cdot\sqrt{p}_{i}/r_{i}^{(t)}\cdot q_{m}}{(A_{1m}^{(t)}\sqrt{p_{1}}/x_{1}+A_{2m}^{(t)}\sqrt{p_{2}}/x_{2}+\cdots+A_{nm}^{(t)}\sqrt{p_{n}}/x_{n})^{2}}\cdot\frac{A_{jm}^{(t)}\sqrt{p_{j}}}{x_{j}^{2}}
\end{align*}
and we can verify that the Jacobian $J^{(t)}(x)$ of $f^{(t)}$ can
be written in matrix-vector notation as 
\begin{align*}
J^{(t)}(x) & =\mathcal{D}(\sqrt{p}/r^{(t)})A^{(t)}\mathcal{D}(q)\cdot\mathcal{D}(A^{(t)T}(\sqrt{p}/x))^{-2}\cdot A^{(t)T}\mathcal{D}(1/x^{2})\mathcal{D}(\sqrt{p})
\end{align*}
When evaluated at the fixed point $x=\sqrt{p}$ of $f^{(t)}(x)$,
we obtain 
\begin{align*}
J^{(t)}:=J^{(t)}(\sqrt{p})=\mathcal{D}(\sqrt{p}/r^{(t)})\cdot A^{(t)}\cdot\mathcal{D}(1/q)A^{(t)T}\mathcal{D}(1/\sqrt{p})
\end{align*}
Our second observation is that $J^{(t)}$ has $\sqrt{p}$ as the unique
maximal eigenvector with eigenvalue equal to 1. Letting $t\rightarrow\infty$,
we have
\begin{align*}
J^{(t)}\rightarrow J:=\mathcal{D}(1/\sqrt{p})\hat{A}\mathcal{D}(1/q)\hat{A}^{T}\mathcal{D}(1/\sqrt{p})
\end{align*}
The unique maximal eigenvector with eigenvalue 1 of $J$ is also $\sqrt{p}$.
Moreover, $J$ is the Jacobian matrix of the limiting map $f$. 

\textbf{III. }We have so far obtained a sequence of mappings $f^{(t)}$
all with fixed point $\sqrt{p}$. The limiting map $f$ has Jacobian
equal to $J=\mathcal{D}(1/\sqrt{p})\hat{A}\mathcal{D}(1/q)\hat{A}^{T}\mathcal{D}(1/\sqrt{p})$,
which has $\sqrt{p}$ as unique maximal eigenvector. Moreover, the
residuals $r^{(t)}/\sqrt{p}-\sqrt{p}$ are always orthogonal to $\sqrt{p}$,
so that convergence is eventually governed by the subdominant eigenvalue
of $J$. We formalize this argument next.\textbf{ }

For any fixed $\epsilon$, since $J^{(t)}\rightarrow J=\mathcal{D}(1/\sqrt{p})\hat{A}\mathcal{D}(1/q)\hat{A}^{T}\mathcal{D}(1/\sqrt{p})$,
there exists $T_{0}$ such that for all $t\geq T_{0}$, $\lambda_{2}(J^{(t)})\leq\lambda_{2}(J)+\epsilon/2$.
Moreover, since $r^{(t)}\rightarrow p$, we can find $T_{1}$ such
that for all $t\geq T_{1}$, $\|r^{(t)}/\sqrt{p}-\sqrt{p}\|_{2}\leq\delta$,
with $\delta$ to be determined later. We use the approximation formula
\begin{align*}
r^{(t+1)}/\sqrt{p} & =f^{(t)}(r^{(t)}/\sqrt{p})\\
 & =f^{(t)}(\sqrt{p}+r^{(t)}/\sqrt{p}-\sqrt{p})\\
 & =\sqrt{p}+J^{(t)}\cdot(r^{(t)}/\sqrt{p}-\sqrt{p})+R_{2}^{(t)}(r^{(t)}/\sqrt{p}-\sqrt{p}).
\end{align*}
Note that $f^{(t)}(x)=\mathcal{D}(\sqrt{p}/r^{(t)})A^{(t)}\mathcal{D}(q)\mathcal{D}(A^{(t)T}(\sqrt{p}/x))^{-1}\mathbf{1}_{m}$
are continuously differentiable, so Taylor's Theorem implies there
exists $T_{2}$ such that for all $t\geq T_{2}$,
\begin{align*}
\|R_{2}^{(t)}(r^{(t)}/\sqrt{p}-\sqrt{p})\|_{2} & \leq C\cdot\|r^{(t)}/\sqrt{p}-\sqrt{p}\|_{2}^{2}
\end{align*}
for some constant $C\geq1$. The bound $C$ can be taken independent
of $t$ because second derivatives of $f^{(t)}$ have no singularities
near $\sqrt{p}$, and are parameterized by $r^{(t)},A^{(t)}$ which
are both uniformly bounded in $t$ by some $C$ near $\sqrt{p}$, using the coercivity of the potential function \eqref{eq:transformed-potential}, which we showed in the proof of \cref{thm:global-convergence}.

Now the key is that $(r^{(t)}/\sqrt{p}-\sqrt{p})\perp\sqrt{p}$, which
is the unique maximal eigenvalue of $J^{(t)}$, so that for all $t\geq0$,
\begin{align*}
\|J^{(t)}\cdot(r^{(t)}/\sqrt{p}-\sqrt{p})\|_{2} & \leq\|J^{(t)}\|_{*}\cdot\|r^{(t)}/\sqrt{p}-\sqrt{p}\|_{2}
\end{align*}
 where $\|J^{(t)}\|_{*}$ is the operator norm of $J^{(t)}$ restricted
to the subspace orthogonal to $\sqrt{p}$. Since we are using the
2-norm and $J^{(t)}$ has unique maximal eigenvalue 1, $\|J^{(t)}\|_{*}$
is precisely $\lambda_{2}(J^{(t)})<1$. Now with $\delta\leq\epsilon/2C$,
for all $t\geq\max\{T_{0},T_{1},T_{2}\}$, we can then bound
\begin{align*}
\|r^{(t+1)}/\sqrt{p}-\sqrt{p}\|_{2} & \leq\|J^{(t)}\cdot(r^{(t)}/\sqrt{p}-\sqrt{p})\|_{2}+\|R_{2}^{(t)}(r^{(t)}/\sqrt{p}-\sqrt{p})\|_{2}\\
 & \leq\lambda_{2}(J^{(t)})\cdot\|r^{(t)}/\sqrt{p}-\sqrt{p}\|_{2}+C\cdot\|r^{(t)}/\sqrt{p}-\sqrt{p}\|_{2}^{2}\\
 & \leq(\lambda_{2}(J)+\epsilon/2)\cdot\|r^{(t)}/\sqrt{p}-\sqrt{p}\|_{2}+C\delta\cdot\|r^{(t)}/\sqrt{p}-\sqrt{p}\|_{2}\\
 & \leq(\lambda_{2}(J)+\epsilon)\cdot\|r^{(t)}/\sqrt{p}-\sqrt{p}\|_{2}.
\end{align*}
 which completes the proof. 

The proof again illustrates the importance of the strong existence and uniqueness conditions in Assumptions \ref{ass:matrix-existence} and \ref{ass:matrix-uniqueness}. When only the weak existence condition in \cref{ass:matrix-weak-existence} is satisfied,
the maximal eigenvalue 1 of $J$ has multiplicity greater than 1. To prove this, we can again use Lemma 1 in \citet{pukelsheim2014biproportional}, which shows that \cref{ass:matrix-weak-existence} holds but \cref{ass:matrix-existence} fails if and only if the limit $\hat A$ of Sinkhorn's algorithm is equivalent to a block diagonal matrix under permutation. Let $N\subseteq [n],M\subseteq[m]$ be such that $\hat A_{ij}=0$ for all $i\in N^C,j\in M$ or $i\in N,j\in M^C$. Then the vector $\tilde{\sqrt{p}}$ with 
\begin{align*}
    (\tilde{\sqrt{p}})_{i}&:=\begin{cases}
\sqrt{p_{i}} & i\in N\\
-\sqrt{p_{i}} & i\in N^{C}
\end{cases}
\end{align*}
is another eigenvector of $J=\mathcal{D}(1/\sqrt{p})\hat{A}\mathcal{D}(1/q)\hat{A}^{T}\mathcal{D}(1/\sqrt{p})$ not in the span of $\sqrt{p}$ with eigenvalue 1. As a result, the orthogonality condition is not sufficient to guarantee contraction and hence asymptotic linear
convergence. $\hfill \square$

 \subsection{Proof of \cref{lem:mm}}
We will show that the iteration in \eqref{eq:scaling-iteration}, which is equivalent to Sinkhorn's algorithm when it is applied to the Luce choice model, is algebraically equivalent to \eqref{eq:mm}. Recall the iteration \eqref{eq:mm} of the MM algorithm in \citet{hunter2004mm}:
\begin{align*}
s_{k}^{(t+1)} & =\frac{w_{k}}{\sum_{i=1}^{n}\sum_{j=1}^{l_{i}-1}\delta_{ijk}[\sum_{j'=j}^{l_{i}}s_{a(i,j')}^{(t)}]^{-1}},
\end{align*}
where $\delta_{ijk}$ is the indicator that item $k\in[m]$ appears in the $i$-th partial ranking $a(i,1)\rightarrow a(i,2)\rightarrow \cdots \rightarrow a(i,l_{i})$, and is not ranked higher than the $j$-th ranked item in that partial ranking. Suppose item $k$ is ranked $\kappa(i,k)$-th in the $i$-th partial ranking,  where $1\leq \kappa(i,k)\leq l_i-1$, i.e., it is not ranked last. Note that $s_{a(i,\kappa(i,k))}=s_k$. We have
\begin{align}
\label{eq:mm-proof-1}
\sum_{j=1}^{l_{i}-1}\delta_{ijk}[\sum_{j'=j}^{l_{i}}s_{a(i,j')}^{(t)}]^{-1} = \frac{1}{s_{a(i,1)}+\cdots+s_{a(i,l_i)}}+\frac{1}{s_{a(i,2)}+\cdots+s_{a(i,l_i)}}+\cdots+\frac{1}{s_{a(i,\kappa(i,k))}+\cdots+s_{a(i,l_i)}},
\end{align}
and if item $k$ is ranked last in the $i$-th partial ranking or does not appear in the $i$-th partial ranking, $\sum_{j=1}^{l_{i}-1}\delta_{ijk}[\sum_{j'=j}^{l_{i}}s_{a(i,j')}^{(t)}]^{-1}\equiv0$. Now if we associate each term in \eqref{eq:mm-proof-1} with a ``choice'' set consisting of the items that appear in the denominator, with the highest ranked item being the ``selected'' item, item $k$ is selected exactly once, in the choice set consisting of $a(i,\kappa(i,k)), \dots,a(i,a(i,l_i))$, corresponding to the last term in \eqref{eq:mm-proof-1}. Thus, each partial ranking of $l_i$ items gives rise to $l_i-1$ observations of choices. Summing over all partial rankings where item $k$ appears and is not last, we see that ${\sum_{i=1}^{n}\sum_{j=1}^{l_{i}-1}\delta_{ijk}[\sum_{j'=j}^{l_{i}}s_{a(i,j')}^{(t)}]^{-1}}$ is exactly equal to the denominator in \eqref{eq:scaling-iteration}:
\begin{align*}
s_{k}^{(t+1)} & =\frac{W_{k}}{\sum_{i\in[n]\mid k\in S_{i}}\frac{R_i}{\sum_{k'\in S_{i}}s_{k'}^{(t)}}},
\end{align*}
and $w_k$, the number of partial rankings in which item $k$ appears and is not ranked last, is exactly the number of times it is ``selected'', i.e., $W_k$.$\hfill \square$

\subsection{Proof of \cref{prop:choicerank}}
Recall the random walk on
$G$  in the network choice model in \citet{maystre2017choicerank}. A user at node $j$ decides to move to node $k\in N_{j}^{\text{out}}$ with probability proportional to $s_{k}$.
To interpret the network choice model as an MNL model, we further define a \emph{choice set} parameterized by each node $j\in V$, consisting of
all items in $N_{j}^{\text{out}}$. Thus there are $m$ \emph{unique} types of choice sets. The number $c_{j}^{\text{out}}$ of transitions
out of node $j$ is exactly the number of observations whose choice sets are those in $N_{j}^{\text{out}}$.  
Therefore, the choice model in \citet{maystre2017choicerank} is a special case of the general Luce choice model with $m$ objects as well as $m$ unique types of choice sets, each appearing in $c_{j}^{\text{out}}$ observations, for a total of $m\cdot c_{j}^{\text{out}}$ choice observations. We index these observations by $i$ as in the Luce choice model. The selected item in an observation at node $j$ corresponds to the node 
into which the transition occurs. The set of nodes $N_{j}^{\text{in}}$
with edges pointing into node $j$ can then be interpreted as the set of all unique choice sets in which node $j$ appears. The number $c_{j}^{\text{in}}$ of observed transitions into a node
$j$ therefore corresponds to the number of observations in which
node $j$ is selected, i.e., it is equal to $W_{j}$ in \eqref{eq:scaling-iteration}. As observed in \citet{maystre2017choicerank}, the exact winner of each observation does not matter, as the sufficient statistics of their model are the aggregate transition counts at each node, i.e., $c_{j}^{\text{out}}, c_{j}^{\text{in}}$. This feature echoes our observation that in the Luce choice model, only the winning frequency $p^1$ matters for solving the ML estimation problem. In fact, $(c_{j}^{\text{out}}, c_{j}^{\text{in}}) = f(A,p^1)$ for some deterministic mapping $f$. 

Now we show that the unregularized version of the ChoiceRank algorithm in \citet{maystre2017choicerank}, given by
\begin{align}
\label{eq:choicerank-proof}
    s_j^{(t+1)} = \frac{c_j^\inn}{\sum_{k\in N_j^\inn} \gamma_k^{(t)}}, \gamma_j^{(t)}=\frac{c_j^\out}{\sum_{k\in N_j^\out}s_k^{(t)}},
\end{align}
is equivalent to \eqref{eq:scaling-iteration} when it is applied to the equivalent Luce choice model. First, note that ${\sum_{k\in N_{j}^{\text{out}}}s_{k}^{(t)}}$
corresponds to the sum of current estimates of $s_k$ over all items in the choice set indexed by node $j$. Furthermore, the summation in $\sum_{k\in N_{j}^{\text{in}}}\gamma_{k}^{(t)}$
is over all unique types of choice sets indexed by node $j$ in which item $k$ appears. As a result, 
\begin{align*}
    \sum_{k\in N_j^\inn} \gamma_k^{(t)} & = \sum_{k\in N_j^\inn} \frac{c_k^\out}{\sum_{k'\in N_k^\out}s_{k'}^{(t)}}\\
    &={\sum_{i\in[n]\mid j\in S_{i}}\frac{1}{\sum_{k'\in S_{i}}s_{k'}^{(t)}}}
\end{align*}
and since $c_{j}^{\text{in}}=W_j$, we have proved that the iteration in \eqref{eq:choicerank-proof} is equal to \begin{align*}
s_{j}^{(t+1)} & =\frac{W_{j}}{\sum_{i\in[n]\mid j\in S_{i}}\frac{R_i}{\sum_{k'\in S_{i}}s_{k'}^{(t)}}},
\end{align*}
which is equal to \eqref{eq:scaling-iteration}. $\hfill \square$

%\subsection{Proof of \cref{thm:sinkhorn-mm}}

\subsection{Proof of \cref{prop:BLP}}
In the simplified setting with $\exp(\theta_{jt})\equiv s_{j}$, the BLP algorithm in \eqref{eq:blp} reduces to \begin{align*}
s_{j}^{(t+1)} & =s_{j}^{(t)}\cdot\hat{p}_{j}\cdot\frac{1}{p_{j}(s^{(t)},\beta,\Gamma,\Sigma)}
\end{align*}
where $p_{j}(\theta^{(m)},\beta,\Gamma,\Sigma)$ denotes the expected market share of the $j$-th product as a function of model parameters. 

Setting $\beta,\Gamma,\Sigma \equiv 0$, the expected market shares are given by 
\begin{align*}
p_{j} = \frac{s_{j}}{\sum_{k}s_{k}}= \frac{1}{n}\sum_{i=1}^n \frac{s_{j}}{\sum_{k}s_{k}}
\end{align*}
Now, consider a Luce choice model of $n$ observations, where the choice set in each observation contains all alternatives. The ``market share'' of product $j$ above is exactly equal to the right hand side in \eqref{eq:optimality} for this model: 
\begin{align*}
\sum_{i\in[n]\mid j\in S_{i}}\frac{1}{n}\frac{s_{j}}{\sum_{k\in S_{i}}s_{k}},
\end{align*}
 while the ``observed market share'' of product $j$ in the BLP model is the fraction of observations where product $j$ is selected, i.e.,  $\hat{p}_{j}=\frac{1}{n}|\{i\in[n]\mid (j,S_i)\}|$, which is equal to the left hand side of \eqref{eq:optimality-original}. The condition $p_{j}=\hat{p}_{j}$ from BLP is exactly 
\begin{align*}
\frac{1}{n}|\{i\in[n]\mid (j,S_i)\}| & =\sum_{i\in[n]\mid j\in S_{i}}\frac{1}{n}\frac{s_{j}}{\sum_{k\in S_{i}}s_{k}}.
\end{align*}

The BLP iteration takes the form 
\begin{align*}
 s_j^{(t+1)}& =s_{j}^{(t)}\cdot\frac{1}{n}|\{i\in[n]\mid (j,S_i)\}|\cdot\left[\sum_{i\in[n]\mid j\in S_{i}}\frac{1}{n}\frac{s_{j}^{(t)}}{\sum_{k\in S_{i}}s_{k}^{(t)}}\right]^{-1}\\
 & =W_{i}\cdot\left[\sum_{i\in[n]\mid j\in S_{i}}\frac{1}{\sum_{k\in S_{i}}s_{k}^{(t)}}\right]^{-1},
\end{align*}
which we recognize to be the algebraic expression \eqref{eq:scaling-iteration} of Sinkhorn's algorithm.  $\hfill \square$
% This argument can be extended
% further to accommodate $s_{j}=\exp(\beta^{T}x_{j})$, in which case
% we use the above iteration to identify $s_{j}$, and then $X\beta=\log s$
% uniquely identifies $\beta$ through 
% \begin{align*}
% \beta & =(X^{T}X)^{-1}X^{T}\log s
% \end{align*}
%  where $X\in\mathbb{R}^{m\times p}$ contains rows $x_{j}\in\mathbb{R}^{p}$,
% if $X^{T}X$ has full rank. 

%% file: skbnd.tex
\newcommand{\assign}{:=}
\newcommand{\tmcolor}[2]{{\color{#1}{#2}}}
\newcommand{\tmop}[1]{\ensuremath{\operatorname{#1}}}
\newcommand{\tmtextbf}[1]{\text{{\bfseries{#1}}}}
\newcommand{\1}{\mathbf{1}}
\renewcommand{\0}{\textbf{0}}

\revision{
\subsection{Proof of Proposition \ref{prop:iteration-complexity}}
We first introduce two supporting lemmas for Proposition \ref{prop:iteration-complexity}.
\begin{lemma}[Theorem 1 of \cite{chakrabarty2021better}]
 Suppose a matrix balancing problem with $(A,p,q)$ has finite positive scaling solutions. For any $\delta>0$, after
  \[ T =\mathcal{O} ( \tfrac{\| p \|_1^2 \log ( \frac{\Delta
     \rho}{\nu} )}{\delta^2} ) \]
  iterations, scaled Sinkhorn matrices $\hat{A}$ satisfy $\| \hat{A} \1_m
  - p \|_1 \leq \delta$ and $\| \hat{A}^T \1_n - q \|_1
  \leq \delta$, where $\Delta$ is the maximum number of nonzeros in any
  column of $A$, $\rho = \max \{ \| p \|_{\infty}, \| q \|_{\infty} \}$ and
  $\nu = \frac{\min_{i j, A_{i j} > 0} A_{i j}}{\max_{i j} A_{i j}} $.
\end{lemma}
\begin{lemma} \label{lem:skbnd}
  Let $d^0_s = \1_n$ be the starting point of the Sinkhorn iteration and $(d^0_\ast, d^1_\ast)$ be an arbitrary pair of optimal scalings. Then for
  all $d^0$ and $d^1$ generated by the Sinkhorn iteration, we have
  \[ \tfrac{1}{\| d^0_{\ast} \|_{\infty}} d^0_{\ast} \leq d^0 \leq
     \tfrac{1}{\| d^0_{\ast} \|_{- \infty}} d^0_{\ast} \]
  \[ \| d^0_{\ast} \|_{- \infty} d^1_{\ast} \leq d^1 \leq \| d^0_{\ast}
     \|_{\infty} d^1_{\ast} . \]
\end{lemma}
\paragraph{Proof of Lemma \ref{lem:skbnd}}
We prove the result by showing that with  $\gamma = \tfrac{1}{\| d^0_{\ast} \|_{\infty}},\eta=\tfrac{1}{\| d^0_{\ast} \|_{- \infty}}$ the following inequality holds for all $d^0$ during the Sinkhorn iterations:
\[ \gamma d^0_{\ast} \leq d^0 \leq \eta d^0_{\ast}. \]
For $d^0_s = \1_n$ the relation $\frac{1}{\| d^0_{\ast} \|_{\infty}} \cdot
\1_n \leq (D^0_{\ast})^{- 1} d^0_s \leq \frac{1}{\| d^0_{\ast} \|_{-
\infty}} \cdot \1_n$ holds by definition. We prove the inequality by induction. Denote $(d^0, d^1)$ and
$(d^0_+, d^1_+)$ as two consecutive Sinkhorn iterations, and assume
$\gamma d^0_{\ast} \leq d^0 \leq \eta d^0_{\ast}$. Then we can write
\begin{align}
  d^0_+ ={} & q / A^{\top} d^1_+ \nonumber\\
  ={} & q / A^{\top} (p / A d^0) \nonumber\\
  ={} & D^0_{\ast} A^{\top} d^1_{\ast} / (A^{\top} ((D^1_{\ast} A
  d^0_{\ast}) / (A d^0))) \nonumber\\
  ={} & D^0_{\ast} \mathcal{D} [A^{\top} ((D^1_{\ast} A d^0_{\ast}) / (A
  d^0))]^{- 1} A^{\top} d^1_{\ast}  \nonumber\\
  ={} & D^0_{\ast} \mathcal{D} [A^{\top} D^1_{\ast}  (\mathcal{D} [A d^0]^{-
  1} (A d^0_{\ast}))]^{- 1} A^{\top} D^1_{\ast} \1_n, \nonumber
\end{align}
where we use the relation $(D a) / b = D\mathcal{D} (b)^{- 1} a$. Since
$\gamma d^0_{\ast} \leq d^0 \leq \eta d^0_{\ast}$ and $A \geq \0$, we get
\[ \gamma A d^0_{\ast} \leq A d^0 \leq \eta A d^0_{\ast} \]
and therefore, $\eta^{- 1} \1_n \leq \mathcal{D} (A d^0)^{- 1} A d^0_{\ast}
\leq \gamma^{- 1} \1_n$ and define $g \assign \mathcal{D} [A d^0]^{- 1} (A
d^0_{\ast})$, we get
\[ d^0_+ = D^0_{\ast} \mathcal{D} [A^{\top} D^1_{\ast} g]^{- 1} A^{\top}
   D^1_{\ast} \1_n, \]
where $\eta^{- 1} \1_n \leq g \leq \gamma^{- 1} \1_n$. Since $A^{\top}
D^1_{\ast} \geq \0$, we can again deduce that
\[ \eta^{- 1} A^{\top} D^1_{\ast}  \1_n \leq A^{\top} D^1_{\ast} g \leq
   \gamma^{- 1} A^{\top} D^1_{\ast} \1_n \]
with $\gamma \1_n \leq \mathcal{D} [A^{\top} D^1_{\ast} g]^{- 1} A^{\top}
D^1_{\ast} \1_n \leq \eta \1_n$. Finally, dividing both sides by $D_\ast^0$ we have $\gamma d^0_{\ast} \leq d^0_+ \leq
\eta d^0_{\ast}$ and this proves the inequality.
}
\revision{
Next we consider $d^1$. Since $d^1 = p / A d^0 = D^1_{\ast} (A d^0_{\ast} /
A d^0)$ and $\gamma d^0_{\ast} \leq d^0 \leq \eta d^0_{\ast}$, we have
\[ \tfrac{1}{\eta} d^1_{\ast} \leq d^1 \leq \tfrac{1}{\gamma} d^1_{\ast}, \]
and this completes the proof.	
\paragraph{Proof of Proposition \ref{prop:iteration-complexity}} Since normalization does not change $\nabla g (u, v)$ or $\nabla^2 g (u, v)$, it
is sufficient to consider the unnormalized version of Sinkhorn. First we show
that the Lipschitz constant of the gradients of the two sub-blocks can be improved over the constants $e^{2B}l_0,e^{2B}l_1$, with $l_0,l_1$ defined in \eqref{eq:smoothness}. Taking $\delta
= \frac{1}{2} \min \{ \| p \|_{- \infty}, \| q \|_{- \infty} \}$ in Lemma 1, we have,
after $\mathcal{O} ( \frac{\| p \|_1^2 \log ( \frac{\Delta
\rho}{\nu} )}{\min \{\| p \|_{- \infty}, \| q \|_{- \infty}\}} )$
iterations, that
\[ \| \hat{A} \1_m - p \|_{\infty} \leq \| \hat{A} \1_m - p
   \|_1 \leq \tfrac{1}{2} \| p \|_{- \infty} \]
and $\hat{A} \1_m \leq \frac{3}{2} p$. Similarly, we obtain $\hat{A}^T \1_n
\leq \frac{3}{2} q$. Recall that $\mathcal{D} ( \hat{A} \1_m )$ and
$\mathcal{D} ( \hat{A}^T \1_n )$ are respectively the Hessians of $g
(u, v)$ with respect to $u$ and $v$, so we get refined Lipschitz constants
\[ l'_0 = \tfrac{3}{2} \| p \|_{\infty}, \quad l'_1 = \tfrac{3}{2} \| q
   \|_{\infty}, \]
which do not depend on $B$. Compared to the bounds in \cref{thm:global-convergence} on the Hessians of sub-blocks of $g(u,v)$, which are given by $e^{2B}l_0,e^{2B}l_1$, the new smoothness constants get rid of the dependence on $e^{2B}$, which allows us to improve the dependence of the convergence rate on $e^B$ in the proof of \cref{thm:global-convergence} from $e^{4B}$ to $e^{2B}$. Moreover, 
since \cref{thm:global-convergence} requires bounds on the smoothness constants that holds for all $t>0$, the global convergence rate depends on $l_0=\|A\mathbf{1}_m\|_\infty,l_1=\|A^T\mathbf{1}_n\|_\infty$; on the other hand, for the complexity bound, we care about the behavior of the algorithm for all large $t>t_0$, which is better captured by $\| p \|_{\infty},\| q \|_{\infty}=\|\hat A\mathbf{1}_m\|_\infty,\|{\hat A}^T\mathbf{1}_n\|_\infty$, defined in terms of the optimal scaled matrix. Now we invoke Lemma \ref{lem:skbnd}, which states
\[ \| d^0_{\ast} \|_{\infty}^{- 1} d^0_{\ast} \leq d^0 \leq \| d^0_{\ast}
   \|_{- \infty}^{- 1} d^0_{\ast}, \qquad \| d^0_{\ast} \|_{- \infty}
   d^1_{\ast} \leq d^1 \leq \| d^0_{\ast} \|_{\infty} d^1_{\ast} . \]
With $u = \log d^0$ and $v = - \log d^1$,
\begin{align}
  \| u \|_{\infty} ={} & \| \log d^0 \| \leq \log ( \tfrac{\| d^0_{\ast}
  \|_{\infty}}{\| d^0_{\ast} \|_{- \infty}} ) \nonumber\\
  \| v \|_{\infty} ={} & \| \log d^1 \| \leq \max \{ \log (
  \tfrac{1}{\| d^0_{\ast} \|_{- \infty} \| d^1_{\ast} \|_{- \infty}}
  ), \log (\| d^0_{\ast} \|_{\infty} \| d^1_{\ast} \|_{\infty})
  \} \nonumber
\end{align}
which implies $e^B \leq C =:\max \{ \tfrac{\| d^0_{\ast} \|_{\infty}}{\|
d^0_{\ast} \|_{- \infty}}, \tfrac{1}{\| d^0_{\ast} \|_{- \infty} \|
d^1_{\ast} \|_{- \infty}}, \| d^0_{\ast} \|_{\infty} \| d^1_{\ast}
\|_{\infty} \}$. We can then invoke the Hessian lower bound
\[ \lambda_{- 2} (\nabla^2 g (u, v)) \geq e^{- 2 B} \lambda_{- 2}
   (\mathcal{L}) \]
in the proof of \cref{thm:global-convergence} to conclude that 
\[ \lambda_{- 2} (\nabla^2 g (u, v)) \geq \frac{1}{C^2} \cdot \lambda_{- 2}
   (\mathcal{L}) \]
and again, invoking Theorem 5.2 of \cite{beck2013convergence} shows that the linear convergence rate is at least
\begin{align*}
    \frac{1}{C^2}\cdot \frac{\min\{\max_j q_j, \max_i p_i\}}{\lambda_{- 2}
   (\mathcal{L})}.
\end{align*}
}
\revision{
Now we show how the linear convergence result on $g(u,v)$ can be
converted to a bound on the $\ell^{1}$ distance $\|r^{(t)}-p\|_{1}$.
Our derivation relies on the property of the optimality gap in \eqref{eq:optimality-gap-property} or the
the following connection between $g(u,v)$ and the KL divergence
$D_{\text{KL}}(p\|r^{(t)})$ (see for example \citet{altschuler2017near}, Lemma 2): 
\begin{align*}
g(u^{(t)},v^{(t)})-g(u^{(t+1)},v^{(t+1)}) & =D_{KL}(p\|r^{(t)})+D_{KL}(q\|c^{(t)}).
\end{align*}
Note that the KL divergence $D_{KL}(p\|r^{(t)})$ above has a different
order from that in \citet{leger2021gradient}, which is $D_{\text{KL}}(r^{(t)}\|p)$, but
this difference does not matter for our analysis, since $D_{KL}(p\|r^{(t)})$
is just an intermediate quantity that is then converted to $\ell^{1}$
distance via Pinsker's inequality: 
\begin{align*}
\|p-r^{(t)}\|_{1} & \leq\sqrt{2D_{KL}(p\|r^{(t)})}.
\end{align*}
 Applying this inequality, we have 
\begin{align*}
\|p-r^{(t)}\|_{1}^{2}	&\leq2(g(u^{(t)},v^{(t)})-g(u^{(t+1)},v^{(t+1)}))
	\\ &\leq2(g(u^{(t)},v^{(t)})-g(u^{\ast},v^{\ast}))
	\\ &\leq2(1-\frac{1}{C^2}\frac{\lambda_{-2}(\mathcal{L})}{\min\{\max_j q_j, \max_i p_i\}})^{t}(g(u^{(0)},v^{(0)})-g(u^{\ast},v^{\ast}))
	\\ &\leq(8\log (C)\sum_{i}p_{i})\cdot \exp{(-\frac{1}{C^2}\frac{\lambda_{-2}(\mathcal{L})}{\min\{\max_j q_j, \max_i p_i\}}\cdot t)}. 
\end{align*}
Suppose after $T$ steps the last line is bounded above by $\varepsilon^2$, then solving for $T$, we have the iteration complexity bound that $\|r^{(t)}-p\|_1\leq \varepsilon$ after
        \begin{align*}
    T = C^2\cdot\frac{\min\{\max_j q_j, \max_i p_i\}}{\lambda_{-2}(\mathcal{L})} \cdot ( \log (1/\varepsilon) + \log(8\log (C)\sum_{i}p_{i}))
      \end{align*}
      iterations of Sinkhorn's algorithm.
$\hfill \square$
}

%% file: main.bbl
\begin{thebibliography}{167}
\providecommand{\natexlab}[1]{#1}
\providecommand{\url}[1]{\texttt{#1}}
\providecommand{\urlprefix}{URL }

\bibitem[{Abowd et~al.(1999)Abowd, Kramarz, \protect\BIBand{} Margolis}]{abowd1999computing}
Abowd JM, Kramarz F, Margolis DN (1999) Computing person and firm effects using linked longitudinal employer-employee data. \emph{Econometrica} 67(6):1411--1453.

\bibitem[{Achilles(1993)}]{achilles1993implications}
Achilles E (1993) Implications of convergence rates in sinkhorn balancing. \emph{Linear Algebra and its Applications} 187:109--112.

\bibitem[{Agarwal et~al.(2010)Agarwal, Negahban, \protect\BIBand{} Wainwright}]{agarwal2010fast}
Agarwal A, Negahban S, Wainwright MJ (2010) Fast global convergence rates of gradient methods for high-dimensional statistical recovery. \emph{Advances in Neural Information Processing Systems} 23.

\bibitem[{Agarwal et~al.(2018)Agarwal, Patil, \protect\BIBand{} Agarwal}]{agarwal2018accelerated}
Agarwal A, Patil P, Agarwal S (2018) Accelerated spectral ranking. \emph{International Conference on Machine Learning}, 70--79.

\bibitem[{Allen-Zhu et~al.(2017)Allen-Zhu, Li, Oliveira, \protect\BIBand{} Wigderson}]{allen2017much}
Allen-Zhu Z, Li Y, Oliveira R, Wigderson A (2017) Much faster algorithms for matrix scaling. \emph{2017 IEEE 58th Annual Symposium on Foundations of Computer Science (FOCS)}, 890--901 (IEEE).

\bibitem[{Altschuler et~al.(2017)Altschuler, Niles-Weed, \protect\BIBand{} Rigollet}]{altschuler2017near}
Altschuler J, Niles-Weed J, Rigollet P (2017) Near-linear time approximation algorithms for optimal transport via sinkhorn iteration. \emph{Advances in neural information processing systems} 30.

\bibitem[{Anderson \protect\BIBand{} Van~Wincoop(2003)}]{anderson2003gravity}
Anderson JE, Van~Wincoop E (2003) Gravity with gravitas: A solution to the border puzzle. \emph{American Economic Review} 93(1):170--192.

\bibitem[{Arjovsky et~al.(2017)Arjovsky, Chintala, \protect\BIBand{} Bottou}]{arjovsky2017wasserstein}
Arjovsky M, Chintala S, Bottou L (2017) Wasserstein generative adversarial networks. \emph{International Conference on Machine Learning}, 214--223 (PMLR).

\bibitem[{Bacharach(1965)}]{bacharach1965estimating}
Bacharach M (1965) Estimating nonnegative matrices from marginal data. \emph{International Economic Review} 6(3):294--310.

\bibitem[{Bacharach(1970)}]{bacharach1970biproportional}
Bacharach M (1970) \emph{Biproportional matrices and input-output change}, volume~16 (CUP Archive).

\bibitem[{Balakrishnan et~al.(2004)Balakrishnan, Hwang, \protect\BIBand{} Tomlin}]{balakrishnan2004polynomial}
Balakrishnan H, Hwang I, Tomlin CJ (2004) Polynomial approximation algorithms for belief matrix maintenance in identity management. \emph{2004 43rd IEEE Conference on Decision and Control (CDC)(IEEE Cat. No. 04CH37601)}, volume~5, 4874--4879 (IEEE).

\bibitem[{Balinski \protect\BIBand{} Pukelsheim(2006)}]{balinski2006matrices}
Balinski M, Pukelsheim F (2006) Matrices and politics. \emph{Festschrift for Tarmo Pukkila on His 60th Birthday} (Department of Mathematics, Statistics and Philosophy, University of Tampere).

\bibitem[{Batsell \protect\BIBand{} Polking(1985)}]{batsell1985new}
Batsell RR, Polking JC (1985) A new class of market share models. \emph{Marketing Science} 4(3):177--198.

\bibitem[{Bauer(1963)}]{bauer1963optimally}
Bauer FL (1963) Optimally scaled matrices. \emph{Numerische Mathematik} 5(1):73--87.

\bibitem[{Beck \protect\BIBand{} Tetruashvili(2013)}]{beck2013convergence}
Beck A, Tetruashvili L (2013) On the convergence of block coordinate descent type methods. \emph{SIAM journal on Optimization} 23(4):2037--2060.

\bibitem[{Berkson(1944)}]{berkson1944application}
Berkson J (1944) Application of the logistic function to bio-assay. \emph{Journal of the American Statistical Association} 39(227):357--365.

\bibitem[{Berry et~al.(1995)Berry, Levinsohn, \protect\BIBand{} Pakes}]{berry1995automobile}
Berry S, Levinsohn J, Pakes A (1995) Automobile prices in market equilibrium. \emph{Econometrica: Journal of the Econometric Society} 841--890.

\bibitem[{Bertsekas(2016)}]{bertsekas2016nonlinear}
Bertsekas D (2016) \emph{Nonlinear Programming}, volume~4 (Athena Scientific).

\bibitem[{Bertsekas(1997)}]{bertsekas1997nonlinear}
Bertsekas DP (1997) Nonlinear programming. \emph{Journal of the Operational Research Society} 48(3):334--334.

\bibitem[{Beurling(1960)}]{beurling1960automorphism}
Beurling A (1960) An automorphism of product measures. \emph{Annals of Mathematics} 189--200.

\bibitem[{Birch(1963)}]{birch1963maximum}
Birch M (1963) Maximum likelihood in three-way contingency tables. \emph{Journal of the Royal Statistical Society Series B: Statistical Methodology} 25(1):220--233.

\bibitem[{Bishop \protect\BIBand{} Nasrabadi(2006)}]{bishop2006pattern}
Bishop CM, Nasrabadi NM (2006) \emph{Pattern recognition and machine learning}, volume~4 (Springer).

\bibitem[{Blanchet et~al.(2022)Blanchet, Chen, \protect\BIBand{} Zhou}]{blanchet2022distributionally}
Blanchet J, Chen L, Zhou XY (2022) Distributionally robust mean-variance portfolio selection with wasserstein distances. \emph{Management Science} 68(9):6382--6410.

\bibitem[{Blanchet et~al.(2016)Blanchet, Gallego, \protect\BIBand{} Goyal}]{blanchet2016markov}
Blanchet J, Gallego G, Goyal V (2016) A markov chain approximation to choice modeling. \emph{Operations Research} 64(4):886--905.

\bibitem[{Blanchet et~al.(2019)Blanchet, Kang, \protect\BIBand{} Murthy}]{blanchet2019robust}
Blanchet J, Kang Y, Murthy K (2019) Robust wasserstein profile inference and applications to machine learning. \emph{Journal of Applied Probability} 56(3):830--857.

\bibitem[{Bong \protect\BIBand{} Rinaldo(2022)}]{bong2022generalized}
Bong H, Rinaldo A (2022) Generalized results for the existence and consistency of the mle in the bradley-terry-luce model. \emph{International Conference on Machine Learning}, 2160--2177 (PMLR).

\bibitem[{Bonnet et~al.(2022)Bonnet, Galichon, Hsieh, O’hara, \protect\BIBand{} Shum}]{bonnet2022yogurts}
Bonnet O, Galichon A, Hsieh YW, O’hara K, Shum M (2022) Yogurts choose consumers? estimation of random-utility models via two-sided matching. \emph{The Review of Economic Studies} 89(6):3085--3114.

\bibitem[{Boyd et~al.(1994)Boyd, El~Ghaoui, Feron, \protect\BIBand{} Balakrishnan}]{boyd1994linear}
Boyd S, El~Ghaoui L, Feron E, Balakrishnan V (1994) \emph{Linear matrix inequalities in system and control theory} (SIAM).

\bibitem[{Boyd et~al.(2006)Boyd, Ghosh, Prabhakar, \protect\BIBand{} Shah}]{boyd2006randomized}
Boyd S, Ghosh A, Prabhakar B, Shah D (2006) Randomized gossip algorithms. \emph{IEEE Transactions on Information Theory} 52(6):2508--2530.

\bibitem[{Bradley(2010)}]{bradley2010algorithms}
Bradley AM (2010) \emph{Algorithms for the equilibration of matrices and their application to limited-memory Quasi-Newton methods}. Ph.D. thesis, Stanford University Stanford University, CA.

\bibitem[{Bradley \protect\BIBand{} Terry(1952)}]{bradley1952rank}
Bradley RA, Terry ME (1952) Rank analysis of incomplete block designs: I. the method of paired comparisons. \emph{Biometrika} 39(3/4):324--345.

\bibitem[{Bregman(1967{\natexlab{a}})}]{bregman1967proof}
Bregman LM (1967{\natexlab{a}}) Proof of the convergence of sheleikhovskii's method for a problem with transportation constraints. \emph{USSR Computational Mathematics and Mathematical Physics} 7(1):191--204.

\bibitem[{Bregman(1967{\natexlab{b}})}]{bregman1967relaxation}
Bregman LM (1967{\natexlab{b}}) The relaxation method of finding the common point of convex sets and its application to the solution of problems in convex programming. \emph{USSR Computational Mathematics and Mathematical Physics} 7(3):200--217.

\bibitem[{Brualdi(1968)}]{brualdi1968convex}
Brualdi RA (1968) Convex sets of non-negative matrices. \emph{Canadian Journal of Mathematics} 20:144--157.

\bibitem[{Bushell(1973)}]{bushell1973hilbert}
Bushell PJ (1973) Hilbert's metric and positive contraction mappings in a banach space. \emph{Archive for Rational Mechanics and Analysis} 52:330--338.

\bibitem[{Carey et~al.(1981)Carey, Hendrickson, \protect\BIBand{} Siddharthan}]{carey1981method}
Carey M, Hendrickson C, Siddharthan K (1981) A method for direct estimation of origin/destination trip matrices. \emph{Transportation Science} 15(1):32--49.

\bibitem[{Carlier(2022)}]{carlier2022linear}
Carlier G (2022) On the linear convergence of the multimarginal sinkhorn algorithm. \emph{SIAM Journal on Optimization} 32(2):786--794.

\bibitem[{Carlier et~al.(2016)Carlier, Chernozhukov, \protect\BIBand{} Galichon}]{carlier2016vector}
Carlier G, Chernozhukov V, Galichon A (2016) Vector quantile regression: An optimal transport approach. \emph{Annals of Statistics} 44(3):1165--1192.

\bibitem[{Carlier et~al.(2023)Carlier, Dupuy, Galichon, \protect\BIBand{} Sun}]{carlier2023sista}
Carlier G, Dupuy A, Galichon A, Sun Y (2023) Sista: learning optimal transport costs under sparsity constraints. \emph{Communications on Pure and Applied Mathematics} 76(9):1659--1677.

\bibitem[{Caron \protect\BIBand{} Doucet(2012)}]{caron2012efficient}
Caron F, Doucet A (2012) Efficient bayesian inference for generalized bradley--terry models. \emph{Journal of Computational and Graphical Statistics} 21(1):174--196.

\bibitem[{Chakrabarty \protect\BIBand{} Khanna(2021)}]{chakrabarty2021better}
Chakrabarty D, Khanna S (2021) Better and simpler error analysis of the sinkhorn--knopp algorithm for matrix scaling. \emph{Mathematical Programming} 188(1):395--407.

\bibitem[{Chang et~al.(2024)Chang, Koehler, Qu, Leskovec, \protect\BIBand{} Ugander}]{chang2024inferring}
Chang S, Koehler F, Qu Z, Leskovec J, Ugander J (2024) Inferring dynamic networks from marginals with iterative proportional fitting. \emph{Forty-first International Conference on Machine Learning}.

\bibitem[{Chang et~al.(2021)Chang, Pierson, Koh, Gerardin, Redbird, Grusky, \protect\BIBand{} Leskovec}]{chang2021mobility}
Chang S, Pierson E, Koh PW, Gerardin J, Redbird B, Grusky D, Leskovec J (2021) Mobility network models of covid-19 explain inequities and inform reopening. \emph{Nature} 589(7840):82--87.

\bibitem[{Chen et~al.(2022)Chen, Kyng, Liu, Peng, Gutenberg, \protect\BIBand{} Sachdeva}]{chen2022maximum}
Chen L, Kyng R, Liu YP, Peng R, Gutenberg MP, Sachdeva S (2022) Maximum flow and minimum-cost flow in almost-linear time. \emph{2022 IEEE 63rd Annual Symposium on Foundations of Computer Science (FOCS)}, 612--623 (IEEE).

\bibitem[{Cohen et~al.(2017)Cohen, Madry, Tsipras, \protect\BIBand{} Vladu}]{cohen2017matrix}
Cohen MB, Madry A, Tsipras D, Vladu A (2017) Matrix scaling and balancing via box constrained newton's method and interior point methods. \emph{2017 IEEE 58th Annual Symposium on Foundations of Computer Science (FOCS)}, 902--913 (IEEE).

\bibitem[{Cottle et~al.(1986)Cottle, Duvall, \protect\BIBand{} Zikan}]{cottle1986lagrangean}
Cottle RW, Duvall SG, Zikan K (1986) A lagrangean relaxation algorithm for the constrained matrix problem. \emph{Naval Research Logistics Quarterly} 33(1):55--76.

\bibitem[{Critchlow et~al.(1991)Critchlow, Fligner, \protect\BIBand{} Verducci}]{critchlow1991probability}
Critchlow DE, Fligner MA, Verducci JS (1991) Probability models on rankings. \emph{Journal of mathematical psychology} 35(3):294--318.

\bibitem[{Cuturi(2013)}]{cuturi2013sinkhorn}
Cuturi M (2013) Sinkhorn distances: Lightspeed computation of optimal transport. \emph{Advances in neural information processing systems} 26.

\bibitem[{De~Paula(2017)}]{de2017econometrics}
De~Paula A (2017) Econometrics of network models. \emph{Advances in economics and econometrics: Theory and applications, eleventh world congress}, 268--323 (Cambridge University Press Cambridge).

\bibitem[{Deming \protect\BIBand{} Stephan(1940)}]{deming1940least}
Deming WE, Stephan FF (1940) On a least squares adjustment of a sampled frequency table when the expected marginal totals are known. \emph{The Annals of Mathematical Statistics} 11(4):427--444.

\bibitem[{Djokovi{\'c}(1970)}]{djokovic1970note}
Djokovi{\'c} D (1970) Note on nonnegative matrices. \emph{Proceedings of the American Mathematical Society} 25(1):80--82.

\bibitem[{Dvurechensky et~al.(2018)Dvurechensky, Gasnikov, \protect\BIBand{} Kroshnin}]{dvurechensky2018computational}
Dvurechensky P, Gasnikov A, Kroshnin A (2018) Computational optimal transport: Complexity by accelerated gradient descent is better than by sinkhorn’s algorithm. \emph{International Conference on Machine Learning}, 1367--1376 (PMLR).

\bibitem[{Dwork et~al.(2001)Dwork, Kumar, Naor, \protect\BIBand{} Sivakumar}]{dwork2001rank}
Dwork C, Kumar R, Naor M, Sivakumar D (2001) Rank aggregation methods for the web. \emph{Proceedings of the 10th International Conference on World Wide Web}, 613--622 (ACM).

\bibitem[{Dykstra(1956)}]{dykstra1956note}
Dykstra O (1956) A note on the rank analysis of incomplete block designs--applications beyond the scope of existing tables. \emph{Biometrics} 12(3):301--306.

\bibitem[{Elo(1978)}]{elo1978rating}
Elo AE (1978) \emph{The rating of chessplayers, past and present} (Arco Pub.).

\bibitem[{Esfahani \protect\BIBand{} Kuhn(2018)}]{mohajerin2018data}
Esfahani PM, Kuhn D (2018) Data-driven distributionally robust optimization using the wasserstein metric: performance guarantees and tractable reformulations. \emph{Mathematical Programming} 171(1-2):115--166.

\bibitem[{Fiedler(1973)}]{fiedler1973algebraic}
Fiedler M (1973) Algebraic connectivity of graphs. \emph{Czechoslovak mathematical journal} 23(2):298--305.

\bibitem[{Fienberg(1970)}]{fienberg1970iterative}
Fienberg SE (1970) An iterative procedure for estimation in contingency tables. \emph{The Annals of Mathematical Statistics} 41(3):907--917.

\bibitem[{Fofana et~al.(2002)Fofana, Lemelin, \protect\BIBand{} Cockburn}]{fofana2002balancing}
Fofana I, Lemelin A, Cockburn J (2002) Balancing a social accounting matrix. \emph{CREFA-Universit{\'e} Laval} .

\bibitem[{Ford(1957)}]{ford1957solution}
Ford LR (1957) Solution of a ranking problem from binary comparisons. \emph{The American Mathematical Monthly} 64(8P2):28--33.

\bibitem[{Ford \protect\BIBand{} Fulkerson(1956)}]{ford1956maximal}
Ford LR, Fulkerson DR (1956) Maximal flow through a network. \emph{Canadian Journal of Mathematics} 8:399--404.

\bibitem[{Ford \protect\BIBand{} Fulkerson(1957)}]{ford1957simple}
Ford LR, Fulkerson DR (1957) A simple algorithm for finding maximal network flows and an application to the hitchcock problem. \emph{Canadian Journal of Mathematics} 9:210--218.

\bibitem[{Fortet(1940)}]{fortet1940resolution}
Fortet R (1940) R{\'e}solution d’un systeme d’{\'e}quations de m. schr{\"o}dinger. \emph{J. Math. Pure Appl. IX} 1:83--105.

\bibitem[{Franklin \protect\BIBand{} Lorenz(1989)}]{franklin1989scaling}
Franklin J, Lorenz J (1989) On the scaling of multidimensional matrices. \emph{Linear Algebra and its Applications} 114:717--735.

\bibitem[{Friedland(2017)}]{friedland2017schrodinger}
Friedland S (2017) On schr{\"o}dinger's bridge problem. \emph{Sbornik: Mathematics} 208(11):1705.

\bibitem[{Gale et~al.(1957)}]{gale1957theorem}
Gale D, et~al. (1957) A theorem on flows in networks. \emph{Pacific J. Math} 7(2):1073--1082.

\bibitem[{Galichon(2018)}]{galichon2018optimal}
Galichon A (2018) \emph{Optimal transport methods in economics} (Princeton University Press).

\bibitem[{Galichon(2021)}]{galichon2021unreasonable}
Galichon A (2021) The unreasonable effectiveness of optimal transport in economics. \emph{arXiv preprint arXiv:2107.04700} .

\bibitem[{Galichon \protect\BIBand{} Salani{\'e}(2021)}]{galichon2021matching}
Galichon A, Salani{\'e} B (2021) Matching with trade-offs: Revealed preferences over competing characteristics. \emph{arXiv preprint arXiv:2102.12811} .

\bibitem[{Gao \protect\BIBand{} Kleywegt(2023)}]{gao2023distributionally}
Gao R, Kleywegt A (2023) Distributionally robust stochastic optimization with wasserstein distance. \emph{Mathematics of Operations Research} 48(2):603--655.

\bibitem[{Gao et~al.(2022)Gao, Ge, \protect\BIBand{} Ye}]{gao2022hdsdp}
Gao W, Ge D, Ye Y (2022) Hdsdp: Software for semidefinite programming. \emph{arXiv preprint arXiv:2207.13862} .

\bibitem[{Georgiou \protect\BIBand{} Pavon(2015)}]{georgiou2015positive}
Georgiou TT, Pavon M (2015) Positive contraction mappings for classical and quantum schr{\"o}dinger systems. \emph{Journal of Mathematical Physics} 56(3):033301.

\bibitem[{Ghosal \protect\BIBand{} Nutz(2022)}]{ghosal2022convergence}
Ghosal P, Nutz M (2022) On the convergence rate of sinkhorn's algorithm. \emph{arXiv preprint arXiv:2212.06000} .

\bibitem[{Good(1963)}]{good1963maximum}
Good IJ (1963) Maximum entropy for hypothesis formulation, especially for multidimensional contingency tables. \emph{The Annals of Mathematical Statistics} 34(3):911--934.

\bibitem[{Gurvits(2004)}]{gurvits2004classical}
Gurvits L (2004) Classical complexity and quantum entanglement. \emph{Journal of Computer and System Sciences} 69(3):448--484.

\bibitem[{Hajek et~al.(2014)Hajek, Oh, \protect\BIBand{} Xu}]{hajek2014minimax}
Hajek B, Oh S, Xu J (2014) Minimax-optimal inference from partial rankings. \emph{Advances in Neural Information Processing Systems}, 1475--1483.

\bibitem[{Hall(1935)}]{hall1935representatives}
Hall P (1935) On representatives of subsets. \emph{Journal of the London Mathematical Society} 1(1):26--30.

\bibitem[{Hausman \protect\BIBand{} Ruud(1987)}]{hausman1987specifying}
Hausman JA, Ruud PA (1987) Specifying and testing econometric models for rank-ordered data. \emph{Journal of Econometrics} 34(1-2):83--104.

\bibitem[{Hazan et~al.(2016)}]{hazan2016introduction}
Hazan E, et~al. (2016) Introduction to online convex optimization. \emph{Foundations and Trends{\textregistered} in Optimization} 2(3-4):157--325.

\bibitem[{Hendrickx et~al.(2020)Hendrickx, Olshevsky, \protect\BIBand{} Saligrama}]{hendrickx2020minimax}
Hendrickx J, Olshevsky A, Saligrama V (2020) Minimax rate for learning from pairwise comparisons in the btl model. \emph{International Conference on Machine Learning}, 4193--4202 (PMLR).

\bibitem[{Hunter(2004)}]{hunter2004mm}
Hunter DR (2004) Mm algorithms for generalized bradley-terry models. \emph{The annals of statistics} 32(1):384--406.

\bibitem[{Idel(2016)}]{idel2016review}
Idel M (2016) A review of matrix scaling and sinkhorn's normal form for matrices and positive maps. \emph{arXiv preprint arXiv:1609.06349} .

\bibitem[{Ireland \protect\BIBand{} Kullback(1968)}]{ireland1968contingency}
Ireland CT, Kullback S (1968) Contingency tables with given marginals. \emph{Biometrika} 55(1):179--188.

\bibitem[{Jochmans \protect\BIBand{} Weidner(2019)}]{jochmans2019fixed}
Jochmans K, Weidner M (2019) Fixed-effect regressions on network data. \emph{Econometrica} 87(5):1543--1560.

\bibitem[{Kalantari \protect\BIBand{} Khachiyan(1993)}]{kalantari1993rate}
Kalantari B, Khachiyan L (1993) On the rate of convergence of deterministic and randomized ras matrix scaling algorithms. \emph{Operations Research Letters} 14(5):237--244.

\bibitem[{Kalantari \protect\BIBand{} Khachiyan(1996)}]{kalantari1996complexity}
Kalantari B, Khachiyan L (1996) On the complexity of nonnegative-matrix scaling. \emph{Linear Algebra and its Applications} 240:87--103.

\bibitem[{Kalantari et~al.(2008)Kalantari, Lari, Ricca, \protect\BIBand{} Simeone}]{kalantari2008complexity}
Kalantari B, Lari I, Ricca F, Simeone B (2008) On the complexity of general matrix scaling and entropy minimization via the ras algorithm. \emph{Mathematical Programming} 112(2):371--401.

\bibitem[{Keener(1993)}]{keener1993perron}
Keener JP (1993) The perron--frobenius theorem and the ranking of football teams. \emph{SIAM review} 35(1):80--93.

\bibitem[{Kendall \protect\BIBand{} Smith(1940)}]{kendall1940method}
Kendall MG, Smith BB (1940) On the method of paired comparisons. \emph{Biometrika} 31(3/4):324--345.

\bibitem[{Kleinberg(1999)}]{kleinberg1999authoritative}
Kleinberg JM (1999) Authoritative sources in a hyperlinked environment. \emph{Journal of the ACM (JACM)} 46(5):604--632.

\bibitem[{Knight(2008)}]{knight2008sinkhorn}
Knight PA (2008) The sinkhorn--knopp algorithm: convergence and applications. \emph{SIAM Journal on Matrix Analysis and Applications} 30(1):261--275.

\bibitem[{Knight \protect\BIBand{} Ruiz(2013)}]{knight2013fast}
Knight PA, Ruiz D (2013) A fast algorithm for matrix balancing. \emph{IMA Journal of Numerical Analysis} 33(3):1029--1047.

\bibitem[{Kruithof(1937)}]{kruithof1937telefoonverkeersrekening}
Kruithof J (1937) Telefoonverkeersrekening. \emph{De Ingenieur} 52:15--25.

\bibitem[{Kuhn et~al.(2019)Kuhn, Esfahani, Nguyen, \protect\BIBand{} Shafieezadeh-Abadeh}]{kuhn2019wasserstein}
Kuhn D, Esfahani PM, Nguyen VA, Shafieezadeh-Abadeh S (2019) Wasserstein distributionally robust optimization: Theory and applications in machine learning. \emph{Operations research \& management science in the age of analytics}, 130--166 (Informs).

\bibitem[{Kullback(1997)}]{kullback1997information}
Kullback S (1997) \emph{Information theory and statistics} (Courier Corporation).

\bibitem[{Kumar et~al.(2015)Kumar, Tomkins, Vassilvitskii, \protect\BIBand{} Vee}]{kumar2015inverting}
Kumar R, Tomkins A, Vassilvitskii S, Vee E (2015) Inverting a steady-state. \emph{Proceedings of the Eighth ACM International Conference on Web Search and Data Mining}, 359--368 (ACM).

\bibitem[{Lamond \protect\BIBand{} Stewart(1981)}]{lamond1981bregman}
Lamond B, Stewart NF (1981) Bregman's balancing method. \emph{Transportation Research Part B: Methodological} 15(4):239--248.

\bibitem[{Landau(1895)}]{landau1895relativen}
Landau E (1895) Zur relativen wertbemessung der turnierresultate. \emph{Deutsches Wochenschach} 11:366--369.

\bibitem[{Lange(2016)}]{lange2016mm}
Lange K (2016) \emph{MM optimization algorithms} (SIAM).

\bibitem[{Lange et~al.(2000)Lange, Hunter, \protect\BIBand{} Yang}]{lange2000optimization}
Lange K, Hunter DR, Yang I (2000) Optimization transfer using surrogate objective functions. \emph{Journal of computational and graphical statistics} 9(1):1--20.

\bibitem[{L{\'e}ger(2021)}]{leger2021gradient}
L{\'e}ger F (2021) A gradient descent perspective on sinkhorn. \emph{Applied Mathematics \& Optimization} 84(2):1843--1855.

\bibitem[{Leontief(1965)}]{leontief1965structure}
Leontief WW (1965) The structure of the us economy. \emph{Scientific American} 212(4):25--35.

\bibitem[{Linial et~al.(1998)Linial, Samorodnitsky, \protect\BIBand{} Wigderson}]{linial1998deterministic}
Linial N, Samorodnitsky A, Wigderson A (1998) A deterministic strongly polynomial algorithm for matrix scaling and approximate permanents. \emph{Proceedings of the thirtieth annual ACM symposium on Theory of computing}, 644--652.

\bibitem[{Luce(1959)}]{luce2012individual}
Luce RD (1959) \emph{Individual choice behavior: A theoretical analysis} (Wiley).

\bibitem[{Luo \protect\BIBand{} Tseng(1992)}]{luo1992convergence}
Luo ZQ, Tseng P (1992) On the convergence of the coordinate descent method for convex differentiable minimization. \emph{Journal of Optimization Theory and Applications} 72(1):7--35.

\bibitem[{Maier et~al.(2010)Maier, Zachariassen, \protect\BIBand{} Zachariasen}]{maier2010divisor}
Maier S, Zachariassen P, Zachariasen M (2010) Divisor-based biproportional apportionment in electoral systems: A real-life benchmark study. \emph{Management Science} 56(2):373--387.

\bibitem[{Marino \protect\BIBand{} Gerolin(2020)}]{marino2020optimal}
Marino SD, Gerolin A (2020) An optimal transport approach for the schr{\"o}dinger bridge problem and convergence of sinkhorn algorithm. \emph{Journal of Scientific Computing} 85(2):27.

\bibitem[{Maystre \protect\BIBand{} Grossglauser(2015)}]{maystre2015fast}
Maystre L, Grossglauser M (2015) Fast and accurate inference of plackett--luce models. \emph{Advances in neural information processing systems}, 172--180.

\bibitem[{Maystre \protect\BIBand{} Grossglauser(2017)}]{maystre2017choicerank}
Maystre L, Grossglauser M (2017) Choicerank: identifying preferences from node traffic in networks. \emph{International Conference on Machine Learning}, 2354--2362 (PMLR).

\bibitem[{McFadden(1978)}]{mcfadden1978modelling}
McFadden D (1978) Modelling the choice of residential location. \emph{Spatial interaction Theory and Planning Models} .

\bibitem[{McFadden(1981)}]{mcfadden1981econometric}
McFadden D (1981) Econometric models of probabilistic choice. \emph{Structural analysis of discrete data with econometric applications} 198272.

\bibitem[{McFadden \protect\BIBand{} Train(2000)}]{mcfadden2000mixed}
McFadden D, Train K (2000) Mixed mnl models for discrete response. \emph{Journal of Applied Econometrics} 15(5):447--470.

\bibitem[{McFadden et~al.(1973)}]{mcfadden1973conditional}
McFadden D, et~al. (1973) Conditional logit analysis of qualitative choice behavior .

\bibitem[{Menon(1968)}]{menon1968matrix}
Menon M (1968) Matrix links, an extremization problem, and the reduction of a non-negative matrix to one with prescribed row and column sums. \emph{Canadian Journal of Mathematics} 20:225--232.

\bibitem[{Negahban et~al.(2012)Negahban, Oh, \protect\BIBand{} Shah}]{negahban2012iterative}
Negahban S, Oh S, Shah D (2012) Iterative ranking from pair-wise comparisons. \emph{Advances in neural information processing systems}, 2474--2482.

\bibitem[{Negahban et~al.(2016)Negahban, Oh, \protect\BIBand{} Shah}]{negahban2016rank}
Negahban S, Oh S, Shah D (2016) Rank centrality: Ranking from pairwise comparisons. \emph{Operations Research} 65(1):266--287.

\bibitem[{Nemirovski \protect\BIBand{} Rothblum(1999)}]{nemirovski1999complexity}
Nemirovski A, Rothblum U (1999) On complexity of matrix scaling. \emph{Linear Algebra and its Applications} 302:435--460.

\bibitem[{Newman(2023)}]{newman2023efficient}
Newman M (2023) Efficient computation of rankings from pairwise comparisons. \emph{Journal of Machine Learning Research} 24(238):1--25.

\bibitem[{Nguyen(1984)}]{nguyen1984estimating}
Nguyen S (1984) Estimating origin destination matrices from observed flows. \emph{Publication of: Elsevier Science Publishers BV} .

\bibitem[{Noothigattu et~al.(2020)Noothigattu, Peters, \protect\BIBand{} Procaccia}]{noothigattu2020axioms}
Noothigattu R, Peters D, Procaccia AD (2020) Axioms for learning from pairwise comparisons. \emph{Advances in Neural Information Processing Systems} 33:17745--17754.

\bibitem[{Orabona(2019)}]{orabona2019modern}
Orabona F (2019) A modern introduction to online learning. \emph{arXiv preprint arXiv:1912.13213} .

\bibitem[{Page et~al.(1999)Page, Brin, Motwani, \protect\BIBand{} Winograd}]{page1999pagerank}
Page L, Brin S, Motwani R, Winograd T (1999) The pagerank citation ranking: Bringing order to the web. Technical report, Stanford InfoLab.

\bibitem[{Peyr{\'e} et~al.(2019)Peyr{\'e}, Cuturi et~al.}]{peyre2019computational}
Peyr{\'e} G, Cuturi M, et~al. (2019) Computational optimal transport: With applications to data science. \emph{Foundations and Trends{\textregistered} in Machine Learning} 11(5-6):355--607.

\bibitem[{Plackett(1975)}]{plackett1975analysis}
Plackett RL (1975) The analysis of permutations. \emph{Journal of the Royal Statistical Society: Series C (Applied Statistics)} 24(2):193--202.

\bibitem[{Plane(1982)}]{plane1982information}
Plane DA (1982) An information theoretic approach to the estimation of migration flows. \emph{Journal of Regional Science} 22(4):441--456.

\bibitem[{Pukelsheim(2006)}]{pukelsheim2006current}
Pukelsheim F (2006) Current issues of apportionment methods. \emph{Mathematics and democracy: recent advances in voting systems and collective choice}, 167--176 (Springer).

\bibitem[{Pukelsheim(2014)}]{pukelsheim2014biproportional}
Pukelsheim F (2014) Biproportional scaling of matrices and the iterative proportional fitting procedure. \emph{Annals of Operations Research} 215:269--283.

\bibitem[{Pukelsheim \protect\BIBand{} Simeone(2009)}]{pukelsheim2009iterative}
Pukelsheim F, Simeone B (2009) On the iterative proportional fitting procedure: Structure of accumulation points and l1-error analysis .

\bibitem[{Pyatt \protect\BIBand{} Round(1985)}]{pyatt1985social}
Pyatt G, Round JI (1985) \emph{Social accounting matrices: A basis for planning} (The World Bank).

\bibitem[{Qu et~al.(2024)Qu, Gao, Hinder, Ye, \protect\BIBand{} Zhou}]{qu2022optimal}
Qu Z, Gao W, Hinder O, Ye Y, Zhou Z (2024) Optimal diagonal preconditioning. \emph{Operations Research} .

\bibitem[{Ragain \protect\BIBand{} Ugander(2016)}]{ragain2016pairwise}
Ragain S, Ugander J (2016) Pairwise choice markov chains. \emph{Advances in Neural Information Processing Systems}, 3198--3206.

\bibitem[{Rajkumar \protect\BIBand{} Agarwal(2014)}]{rajkumar2014statistical}
Rajkumar A, Agarwal S (2014) A statistical convergence perspective of algorithms for rank aggregation from pairwise data. \emph{International Conference on Machine Learning}, 118--126.

\bibitem[{Ruiz(2001)}]{ruiz2001scaling}
Ruiz D (2001) A scaling algorithm to equilibrate both rows and columns norms in matrices. Technical report, CM-P00040415.

\bibitem[{Ruschendorf(1995)}]{ruschendorf1995convergence}
Ruschendorf L (1995) Convergence of the iterative proportional fitting procedure. \emph{The Annals of Statistics} 1160--1174.

\bibitem[{Schneider \protect\BIBand{} Zenios(1990)}]{schneider1990comparative}
Schneider MH, Zenios SA (1990) A comparative study of algorithms for matrix balancing. \emph{Operations Research} 38(3):439--455.

\bibitem[{Schr{\"o}dinger(1931)}]{schro1931uber}
Schr{\"o}dinger E (1931) {\"U}ber die umkehrung der naturgesetze. \emph{Sitzungsberichte der preussischen Akademie der Wissenschaften, physikalisch-mathematische Klasse} 8(9):144--153.

\bibitem[{Seshadri et~al.(2020)Seshadri, Ragain, \protect\BIBand{} Ugander}]{seshadri2020learning}
Seshadri A, Ragain S, Ugander J (2020) Learning rich rankings. \emph{Advances in Neural Information Processing Systems} 33:9435--9446.

\bibitem[{Shah et~al.(2015)Shah, Balakrishnan, Bradley, Parekh, Ramchandran, \protect\BIBand{} Wainwright}]{shah2015estimation}
Shah N, Balakrishnan S, Bradley J, Parekh A, Ramchandran K, Wainwright M (2015) Estimation from pairwise comparisons: Sharp minimax bounds with topology dependence. \emph{Artificial intelligence and statistics}, 856--865 (PMLR).

\bibitem[{Sheffi(1985)}]{sheffi1985urban}
Sheffi Y (1985) \emph{Urban transportation networks}, volume~6 (Prentice-Hall, Englewood Cliffs, NJ).

\bibitem[{Silva \protect\BIBand{} Tenreyro(2006)}]{silva2006log}
Silva JS, Tenreyro S (2006) The log of gravity. \emph{The Review of Economics and Statistics} 88(4):641--658.

\bibitem[{Sinkhorn(1964)}]{sinkhorn1964relationship}
Sinkhorn R (1964) A relationship between arbitrary positive matrices and doubly stochastic matrices. \emph{The Annals of Mathematical Statistics} 35(2):876--879.

\bibitem[{Sinkhorn(1974)}]{sinkhorn1974diagonal}
Sinkhorn R (1974) Diagonal equivalence to matrices with prescribed row and column sums. ii. \emph{Proceedings of the American Mathematical Society} 45(2):195--198.

\bibitem[{Sinkhorn \protect\BIBand{} Knopp(1967)}]{sinkhorn1967concerning}
Sinkhorn R, Knopp P (1967) Concerning nonnegative matrices and doubly stochastic matrices. \emph{Pacific Journal of Mathematics} 21(2):343--348.

\bibitem[{Soufiani et~al.(2013)Soufiani, Chen, Parkes, \protect\BIBand{} Xia}]{soufiani2013generalized}
Soufiani HA, Chen W, Parkes DC, Xia L (2013) Generalized method-of-moments for rank aggregation. \emph{Advances in Neural Information Processing Systems}, 2706--2714.

\bibitem[{Soules(1991)}]{soules1991rate}
Soules GW (1991) The rate of convergence of sinkhorn balancing. \emph{Linear Algebra and its Applications} 150:3--40.

\bibitem[{Spielman(2007)}]{spielman2007spectral}
Spielman DA (2007) Spectral graph theory and its applications. \emph{48th Annual IEEE Symposium on Foundations of Computer Science (FOCS'07)}, 29--38 (IEEE).

\bibitem[{Stellato et~al.(2020)Stellato, Banjac, Goulart, Bemporad, \protect\BIBand{} Boyd}]{stellato2020osqp}
Stellato B, Banjac G, Goulart P, Bemporad A, Boyd S (2020) Osqp: An operator splitting solver for quadratic programs. \emph{Mathematical Programming Computation} 12(4):637--672.

\bibitem[{Stone(1962)}]{stone1962multiple}
Stone R (1962) Multiple classifications in social accounting. \emph{Bulletin de l’Institut International de Statistique} 39(3):215--233.

\bibitem[{Stone \protect\BIBand{} Brown(1971)}]{stone1971computable}
Stone R, Brown A (1971) \emph{A computable model of economic growth}.

\bibitem[{Ta{\c{s}}kesen et~al.(2023)Ta{\c{s}}kesen, Shafieezadeh-Abadeh, \protect\BIBand{} Kuhn}]{tacskesen2023semi}
Ta{\c{s}}kesen B, Shafieezadeh-Abadeh S, Kuhn D (2023) Semi-discrete optimal transport: Hardness, regularization and numerical solution. \emph{Mathematical Programming} 199(1):1033--1106.

\bibitem[{Theil(1967)}]{theil1967economics}
Theil H (1967) \emph{Economics and Information Theory} (North-Holland Publishing Company).

\bibitem[{Theil \protect\BIBand{} Rey(1966)}]{theil1966quadratic}
Theil H, Rey G (1966) A quadratic programming approach to the estimation of transition probabilities. \emph{Management Science} 12(9):714--721.

\bibitem[{Thionet(1964)}]{thionet1964note}
Thionet P (1964) Note sur le remplissage d'un tableau {\`a} double entr{\'e}e. \emph{Journal de la Soci{\'e}t{\'e} Fran{\c{c}}aise de Statistique} 105:228--247.

\bibitem[{Thurstone(1927)}]{thurstone1927method}
Thurstone LL (1927) The method of paired comparisons for social values. \emph{The Journal of Abnormal and Social Psychology} 21(4):384.

\bibitem[{Tomlin(2003)}]{tomlin2003new}
Tomlin JA (2003) A new paradigm for ranking pages on the world wide web. \emph{Proceedings of the 12th international Conference on World Wide Web}, 350--355.

\bibitem[{Tseng(2001)}]{tseng2001convergence}
Tseng P (2001) Convergence of a block coordinate descent method for nondifferentiable minimization. \emph{Journal of Optimization Theory and Applications} 109:475--494.

\bibitem[{Tseng \protect\BIBand{} Bertsekas(1987)}]{tseng1987relaxation}
Tseng P, Bertsekas DP (1987) Relaxation methods for problems with strictly convex separable costs and linear constraints. \emph{Mathematical Programming} 38(3):303--321.

\bibitem[{Tversky(1972)}]{tversky1972elimination}
Tversky A (1972) Elimination by aspects: A theory of choice. \emph{Psychological Review} 79(4):281.

\bibitem[{Uribe et~al.(1966)Uribe, De~Leeuw, \protect\BIBand{} Theil}]{uribe1966information}
Uribe P, De~Leeuw C, Theil H (1966) The information approach to the prediction of interregional trade flows. \emph{The Review of Economic Studies} 33(3):209--220.

\bibitem[{Villani et~al.(2009)}]{villani2009optimal}
Villani C, et~al. (2009) \emph{Optimal transport: old and new}, volume 338 (Springer).

\bibitem[{Vojnovic \protect\BIBand{} Yun(2016)}]{vojnovic2016parameter}
Vojnovic M, Yun S (2016) Parameter estimation for generalized thurstone choice models. \emph{International Conference on Machine Learning}, 498--506 (PMLR).

\bibitem[{Vojnovic et~al.(2020)Vojnovic, Yun, \protect\BIBand{} Zhou}]{vojnovic2020convergence}
Vojnovic M, Yun SY, Zhou K (2020) Convergence rates of gradient descent and mm algorithms for bradley-terry models. \emph{International Conference on Artificial Intelligence and Statistics}, 1254--1264 (PMLR).

\bibitem[{Vojnovic et~al.(2023)Vojnovic, Yun, \protect\BIBand{} Zhou}]{vojnovic2023accelerated}
Vojnovic M, Yun SY, Zhou K (2023) Accelerated mm algorithms for inference of ranking scores from comparison data. \emph{Operations Research} 71(4):1318--1342.

\bibitem[{Wilson(1969)}]{wilson1969use}
Wilson AG (1969) The use of entropy maximising models, in the theory of trip distribution, mode split and route split. \emph{Journal of Transport Economics and Policy} 108--126.

\bibitem[{Xiao et~al.(2007)Xiao, Boyd, \protect\BIBand{} Kim}]{xiao2007distributed}
Xiao L, Boyd S, Kim SJ (2007) Distributed average consensus with least-mean-square deviation. \emph{Journal of Parallel and Distributed Computing} 67(1):33--46.

\bibitem[{Yule(1912)}]{yule1912methods}
Yule GU (1912) On the methods of measuring association between two attributes. \emph{Journal of the Royal Statistical Society} 75(6):579--652.

\bibitem[{Zermelo(1929)}]{zermelo1929berechnung}
Zermelo E (1929) Die berechnung der turnier-ergebnisse als ein maximumproblem der wahrscheinlichkeitsrechnung. \emph{Mathematische Zeitschrift} 29(1):436--460.

\end{thebibliography}
